\documentclass{article}

\usepackage{PRIMEarxiv}

\usepackage[utf8]{inputenc} % allow utf-8 input
\usepackage[T1]{fontenc}    % use 8-bit T1 fonts
\usepackage{url}            % simple URL typesetting
\usepackage{booktabs}       % professional-quality tables
\usepackage{amsfonts}       % blackboard math symbols
\usepackage{nicefrac}       % compact symbols for 1/2, etc.
\usepackage{microtype}      % microtypography
\usepackage{lipsum}
\usepackage{fancyhdr}       % header
\usepackage{graphicx}       % graphics% Required for inserting images
\usepackage{chemfig}
\usepackage{wrapfig}        % to have text wrapping the figure
\newif\ifincludegraphs % Change this to \includegraphsfalse to exclude graphs
\includegraphsfalse
\usepackage{float}  % For [H] specifier to be able to specify the position of a figure
\usepackage{amsmath} % For \text{} in captions
\usepackage[ruled, linesnumbered, vlined]{algorithm2e}
\usepackage{todonotes}
\usepackage{longtable}
\usepackage{caption}
\usepackage{multicol} %to plot images side by side
\usepackage{subcaption} %to give multiple subplots an individual caption
\usepackage{bm} %to use bold letters in math style
\usepackage{amssymb}
\graphicspath{{media/}}     % organize your images and other figures under media/ folder
\usepackage[sorting=none]{biblatex}
\bibliography{My_Library}
\usepackage{blkarray} %to be able to give a matrix a caption
\usepackage{booktabs}
\usepackage{nicematrix}
\usepackage{tikz}% added 
\usetikzlibrary{decorations.pathreplacing, calligraphy}% added

\title{Surrogate-Based Optimization Techniques for Process Systems Engineering
}

\author{
  \begin{tabular}{ccc}
    \textbf{Mathias Neufang} & \hspace{0.5cm} & \textbf{Emma Pajak} \\
    \textnormal{Imperial College London} & \hspace{0.5cm} & \textnormal{Imperial College London} \\
    \texttt{mathias.neufang22@imperial.ac.uk} &  \hspace{0.5cm} & \\
    \\
    \textbf{Damien van de Berg} & \hspace{0.5cm} & \textbf{Ye Seol Lee} \\
    \textnormal{Imperial College London} & \hspace{0.5cm} & \textnormal{University College London} \\
  \end{tabular} 
  \AND
  \textbf{Ehecatl Antonio del Rio Chanona} \\
  \textnormal{Imperial College London} \\
  \texttt{a.del-rio-chanona@imperial.ac.uk}
}
\begin{document}
\maketitle

\begin{abstract}
Optimization plays an important role in chemical engineering, impacting cost-effectiveness, resource utilization, product quality, and process sustainability metrics. This chapter broadly focuses on data-driven optimization, particularly, on model-based derivative-free techniques, also known as surrogate-based optimization. The chapter introduces readers to the theory and practical considerations of various algorithms, complemented by a performance assessment across multiple dimensions, test functions, and two chemical engineering case studies: a stochastic high-dimensional reactor control study and a low-dimensional constrained stochastic reactor optimization study. This assessment sheds light on each algorithm's performance and suitability for diverse applications. Additionally, each algorithm is accompanied by background information, mathematical foundations, and algorithm descriptions. Among the discussed algorithms are Bayesian Optimization (BO), including state-of-the-art TuRBO, Constrained Optimization by Linear Approximation (COBYLA), the Ensemble Tree Model Optimization Tool (ENTMOOT) which uses decision trees as surrogates, Stable Noisy Optimization by Branch and Fit (SNOBFIT), methods that use radial basis functions such as DYCORS and SRBFStrategy, Constrained Optimization by Quadratic Approximations (COBYQA), as well as a few others recognized for their effectiveness in surrogate-based optimization. By combining theory with practice, this chapter equips readers with the knowledge to integrate surrogate-based optimization techniques into chemical engineering. The overarching aim is to highlight the advantages of surrogate-based optimization, introduce state-of-the-art algorithms, and provide guidance for successful implementation within process systems engineering.
\end{abstract}

\keywords{Data-Driven Optimization, \and Model-Based Optimization, \and Derivative-Free Optimization, \and Black-Box Optimization, \and Chemical Process Engineering}
% \classification{65C05, 62M20, 93E11, 62F15, 86A22}

%\affil[1]{Department of Chemical Engineering, Imperial College London, Exhibition Rd, South Kensington, London SW7 2AZ, United Kingdom, e-mail: mathias.neufang22@imperial.ac.uk}
%\affil[2]{Department of Chemical Engineering, Imperial College London}
%\affil[3]{Department of Chemical Engineering, Imperial College London}
%\affil[4]{Department of Chemical Engineering, University College London}

%\makecontributiontitle
%  \DOI{10.1515/futur-2012-0001}

%\setlength{\parindent}{0pt} %switch-off autoindent by latex 
%\vspace{1em}
%\hrule
%\vspace{1em}

\newpage
\section{Introduction}\label{sec:intro}
This chapter aims to equip process systems and chemical engineers with a comprehensive understanding of data-driven optimization, focusing specifically on model-based (surrogate) techniques. It explores both the theoretical foundations and the practical performance assessment of various algorithms, fostering an intuitive grasp of their behavior and efficacy. This understanding will enable readers to utilize these methods effectively. The remaining chapter takes the following structure:

\begin{itemize}
    \item \textbf{Section \ref{sec:intro} - Introduction:} This section sets the scene for the chapter, defining derivative-free optimization (DFO), introducing the three types of DFO methods, and exploring their application in chemical engineering.
    
    \item \textbf{Section \ref{sec:mb dfo} - Model-Based Derivative-Free Optimization:} Here, readers are formally introduced to the concept of model-based DFO, where popular algorithms like Bayesian Optimization take center stage. Background information, mathematical formulations, and algorithm descriptions are presented for a holistic overview of each method.
    
    \item \textbf{Section \ref{sec:b} - Unconstrained Performance Assessment:} The performance assessment focuses on benchmarking model-based algorithms for unconstrained problems. The objective is to equip the reader with sufficient intuition and background information to make informed decisions when selecting algorithms.

    \item \textbf{Section \ref{sec:c} - Constrained Performance Assessment:} This part of the performance assessment focuses on benchmarking model-based algorithms for constrained problems. The objectives mirror those in Section 3 but in the case of constrained black-box problems.
    
    \item \textbf{Section \ref{sec:d} - Chemical Engineering Case Studies:} 
    In this section, two chemical engineering optimization case studies are presented, an unconstrained, stochastic high-dimensional problem, and a stochastic low-dimensional constrained problem. This section serves to demonstrate the applicability of the surrogate-based algorithms to the process systems engineering field.

    \item \textbf{Section \ref{sec:e} - Code:}
    Alongside this study a GitHub repository is published. Here, the interested reader can find all algorithm implementations, as well as the benchmarking library developed for this work.

    \item \textbf{Section \ref{sec:f} - Concluding Remarks:}
    This section summarizes the outcomes of the book chapter's benchmarking section. It emphasizes the most effective algorithms, as well as their drawbacks. Furthermore, future outlook on surrogate-based optimization is provided. 
\end{itemize}

\subsection{Background}
Traditionally, optimization within chemical engineering relies on algebraic expressions or established knowledge-based models which can be optimized by leveraging derivative information from their analytical expressions. However, the rise of digitalization, such as smart measuring devices, process analytical technology, sensor technologies, cloud platforms, and the \textit{Industrial Internet of Things} (IIoT) has called for the need for optimization algorithms guided purely by the collected data, and therefore the term \textit{data-driven optimization} has emerged. In complex chemical systems, it is not unusual for data collection to be feasible only through the evaluation of an expensive black-box function. This function may represent an in-vitro chemical experiment with undetermined mechanisms or a costly process reconfiguration as well as an in-silico simulation in the form of computational fluid dynamics, or quantum mechanical calculations. Evaluating such deterministic models often relies on complex and expensive simulations that are corrupted by computational noise as unintended variations and inaccuracies happen, which makes even numerical derivatives difficult and unreliable \cite{moreBenchmarkingDerivativeFreeOptimization2009}. In such instances, data-driven algorithms emerge as a solution, enabling the optimization of these systems \cite{vandebergDatadrivenOptimizationProcess2022}.
\\[8pt]
Many engineering optimization challenges can be framed as 'costly' black box problems, which are constrained by the number of function evaluations. Engineers often construct precise models of physical systems that are either differentiable or economical to evaluate. These models can be resolved efficiently, and their solutions can be applied to the actual system. However, when gradient information or cost-effective models are unavailable, it becomes necessary to utilize efficient optimization methods that depend solely on function evaluations. The process of developing a model can be viewed as an integral part of the expensive black-box optimization procedure itself \cite{vandebergDatadrivenOptimizationProcess2022}. Within this context, algorithms must rise to the challenge, with optimization, \textit{Artificial Intelligence} (AI), and machine learning playing pivotal roles in enabling advancements in automated control and decision-making \cite{boglePerspectiveSmartProcess2017}. In fact, data-driven optimization boasts a rich historical background within the realm of chemical engineering \cite{bieglerMultiscaleOptimizationProcess2014, wilsonALAMOApproachMachine2017}.
%In situations where observations entail significant costs, there is a preference for optimization algorithms that adjust their behavior based on evolving information. This adaptability poses a primary challenge in the field of optimization \cite{garnettBayesianOptimization2023}.  
\\[8pt]
% The necessity to base algorithmic decisions on data and the absence of analytical expressions for objective functions have led to various terminology to describe this family of optimization algorithms as several communities have developed solution methods for the same problem.
Similar challenges have been the topic of research across many communities.
This chapter refers to "Data-Driven Optimization". However, literature similarly refers to this subject area as "derivative-free optimization", "gradient-free optimization", "zeroth-order optimization", "simulation-based optimization", or, more specifically within the context of process systems engineering, "black-box optimization"  \cite{bieglerMultiscaleOptimizationProcess2014, wilsonALAMOApproachMachine2017}. The term "zeroth-order" alludes to algorithms that do not utilize a function's first or second-order derivatives \cite{shiNumericalPerformanceFinitedifferencebased2023}.

\subsection{Unconstrained Optimization Formulation}
Problem \ref{eq:basic_opt} represents the generic, unconstrained optimization formulation of interest throughout the chapter, where $f:\mathbb{R}^{n_x} \longrightarrow \mathbb{R}$ describes the objective function:
\begin{equation}\label{eq:basic_opt}
\begin{aligned}
     & \min_{\mathbf{x}} && f(\mathbf{x}) \\
     &&& \mathbf{x} \in \mathcal{X} \subseteq \mathbb{R}^{n_x} \\
\end{aligned}
\end{equation}

In traditional numerical optimization, given an analytical expression of $f$, the necessary conditions for a (local) optimal solution of Problem (\ref{eq:basic_opt}) would be any point or region where the gradient is zero $\nabla_x f(\bold{x}) = \bold{0}$, and the sufficient condition would, in addition, have the Hessian matrix be positive semi-definite $\nabla_{xx}^2 f(\bold{x}) \succeq 0$. Therefore, optimization algorithms that use derivatives (e.g. Newton's Method, gradient descent) seek to find regions where these conditions (particularly the necessary conditions) are met.
\\[8 pt]
In \textit{Data-Driven Optimization} (DDO), things are not as 'simple'. Given that there is no analytical expression, algorithms must seek to explore the space in the hope of gathering adequate information about the function at hand, while at the same time using this information to optimize the function. The general assumption for practical applications in DDO is that ensuring convergence to an optimum is hard given that the function itself is unknown, and therefore termination criteria are generally set in the number of evaluations or runtime \cite{neumaierCompleteSearchContinuous2004a}.
\\[8 pt]
\textit{Derivative-Free Optimization} (DFO), a term that will be used in this chapter interchangeably with DDO, encompasses algorithms designed to optimize functions without explicitly using derivative information. DFO methods can be broadly classified into two main categories: \textit{model-based} derivative-free methods, also known as \textit{surrogate-based} optimization, and \textit{direct derivative-free} methods. Situated in-between these categories are the \textit{finite-difference} methods, which approximate derivatives using function evaluations \cite{shiNumericalPerformanceFinitedifferencebased2023}, and for second or higher order methods, it can be viewed as constructing a local approximation (a model) of the system. In this chapter, the focus is on surrogate-based optimization methods. A brief overview of the different methods is presented next, before focusing on model-based derivative-free (also known as surrogate-based) methods.

\subsection{Direct Methods}
Direct derivative-free methods directly utilize sampled points to determine the next sampling location to approach the optimum without relying on the intermediate construction and optimization of surrogates. Many early DFO algorithms, such as the Simplex (Nelder-Mead) algorithm, evolutionary (particle swarm) algorithms (which are a type of metaheuristics), as well as direct global optimization algorithms, simulated annealing, and pattern search, fall into this category \cite{vandebergDatadrivenOptimizationProcess2022}.

\subsection{Finite Difference Methods}
Finite-difference methods operate on the principle that even in cases where the overall function is unknown, sampled data can be utilized to approximate derivatives. By approximating the gradient, typical gradient descent methods can be applied. These methods fall between direct derivative-free methods and surrogate-based DFO methods; they base sampling decisions directly on obtained data, similar to direct methods. However, when second-order or quasi-Newton methods are employed, they can be thought of as quadratic surrogates of the underlying true function, in some cases with iterative updates and 'memory'. Well-known finite difference methods include the Gradient Descent + Momentum Method, Adam, RSM Prob, and (L)BFGS. The interested reader is referred to the work on \cite{shiNumericalPerformanceFinitedifferencebased2023} for further treatment on the subject.

\subsection{Model-Based Methods}\label{sec:mbintro}
In process engineering, there has been a notable surge in interest in model-based techniques, drawing upon the field's expertise in surrogate modeling, meta-modeling, and reduced-order modeling \cite{bieglerMultiscaleOptimizationProcess2014, wilsonALAMOApproachMachine2017}. Traditionally, within the process systems engineering community, model-based DFO methods are also referred to as surrogate optimization, highlighting the creation of a 'surrogate' function \cite{bhosekarAdvancesSurrogateBased2018, kimMachineLearningbasedSurrogate2020, schweidtmannDeterministicGlobalProcess2019}. As these terms are interchangeable, this chapter will refer to both throughout.
\\[8pt]
The fundamental concept of model-based DFO involves creating a model of $f$ using data, known as a surrogate function $\hat{f}$, once this surrogate model is built it can be used to determine candidate solutions. The general framework of model-based DFO seeks to explore the function $f$ to refine its surrogate $\hat{f}$ while at the same time finding points that optimize $f$. The underlying shape of $f$ remains unknown, and the choice of model type (e.g., \textit{Radial Basis Functions} (RBFs), Gaussian Processes, quadratic surrogates, neural networks, etc.) used to fit $\hat{f}$ makes some assumptions on how $f$ behaves. Section \ref{sec:mb dfo} is devoted to the comprehensive treatment of model-based derivative-free methods.

\subsection{Process Systems Engineering Applications of Model-Based Methods}\label{sec:mbapps}
Model-based DFO has been extensively used in chemical engineering, and the community has a rich history of developing state-of-the-art algorithms for various applications across a variety of fields. For instance, in quantum chemistry, using a classical state vector simulator as a surrogate model to approximate the optimization landscape in variational quantum eigensolvers, enhancing the convergence and efficiency of quantum circuit optimization \cite{gustafsonSurrogateOptimizationVariational2024}. In a different vein, model-based DDO has also found applications in the pharmaceutical manufacturing sector, as described in a recent paper where feasibility-driven optimization incorporates additional stages to improve both local exploitation and global exploration. This approach leads to lower costs, improved product quality, and greater process flexibility and robustness \cite{tianSurrogateBasedFramework2024}. This section highlights some of the recent research and implementation of model-based DDO in chemical engineering, grouped by application areas to showcase the versatility and efficacy of these approaches in addressing complex challenges across the chemical engineering domain.
\\[8pt]
\textbf{Process Design \& Flowsheeting:} Process design and flowsheeting involve strategic planning and optimization of chemical processes and systems. Advanced methodologies and tools are employed to create optimized, efficient, and robust processes, emphasizing operational efficiency and environmental sustainability. For example, surrogate-based optimization strategies have been proposed for the global optimization of process flowsheets, leveraging algebraic surrogates constructed from rigorous simulations via Bayesian symbolic regression \cite{forsterAlgebraicSurrogateBased2023}, other approaches have suggested the use of Gaussian Processes (kriging) to optimize nonlinear problems that include noisy implicit black box functions, such as modular process simulators, to manage noise and compensate for unavailable derivatives \cite{caballeroAlgorithmUseSurrogate2008a}. Other related works have focused on the use of Gaussian Processes as versatile surrogate modeling tools to address feasibility constraints from non-converged simulations based on performance-risk trade-offs \cite{durkinSurrogatebasedOptimisationProcess2024}. Another study introduces a multi-fidelity Bayesian optimization framework for reactor design, reducing design time for highly parameterized reactors while ensuring optimal geometry and operating conditions with experimental validation of 3D-printed reactor geometries \cite{savageMultifidelityDatadrivenDesign2023b}. An overview of process design via Bayesian optimization can be found in the following work \cite{paulsonBayesianOptimizationFlexible2024}. 
\\[8pt]
Moving away from Bayesian approaches, other studies on specific processes demonstrate the potential of surrogate-based methods in process design, such as the optimization of a hybrid polycrystalline silicon production route, where surrogate models of key unit operations are constructed to enable entire process optimization exploring various scenarios, e.g., maximizing silicon production, minimizing operating costs, and maximizing total profit \cite{ramirez-marquezSurrogateBasedOptimization2020}. Similarly, a trust region filter framework for heat exchanger network synthesis integrates detailed shell-and-tube heat exchanger models to optimize network topology and exchanger design, including parameters such as pressure drops, shell configurations, and tube arrangements \cite{kaziTrustRegionFramework2021}.  Surrogate-based optimization has also been applied to process systems for resource recovery from wastewater, integrating DFO modeling tools that incorporate classification surrogate models and address uncertainties, thereby offering holistic solutions to reduce environmental impacts in food and beverage production \cite{durkinSurrogatebasedOptimisationProcess2024}. In addition to the above, many other surrogates have been used, from \textit{Graph Neural Networks} (GNNs) for Granular Flows \cite{jiangIntegratingGraphNeural2024}, to symbolic regression-based surrogates for flexibility analysis \cite{forsterAlgebraicSurrogatebasedFlexibility2024}, to Quantile Neural Networks for Two-stage Stochastic Optimization \cite{alcantaraQuantileNeuralNetwork2024}. An overview of data-driven and hybrid models for subsequent optimization of separation processes can be found in \cite{mcbrideHybridSemiParametric2020}.
 \\[8pt]
\textbf{Supply Chains \& Planning, Scheduling, and Operations:} Surrogate-based DDO is also a powerful tool utilized within supply chains, exemplified by a study where historical data models optimize demand response scheduling of air separation units within a nonlinear dynamic framework, ensuring dynamic feasibility and computational efficiency \cite{tsayDataDrivenModelsAlgorithms2018}. Another example is the extension of the DOMINO framework to tackle mixed-integer bi-level multi-follower stochastic optimization problems for integrated planning and scheduling problems under demand uncertainty  \cite{beykalDatadrivenOptimizationMixedinteger2022a}. In a different study on scheduling optimization in integrated chemical plants, \textit{Convex Region Surrogate} (CRS) models are used in mixed-integer programming frameworks, of which the effectiveness is demonstrated through its application to an industrial test case from a Praxair plant \cite{zhangDatadrivenConstructionConvex2016}. Furthermore, the automation and simulation of plant-level surrogate construction, as well as a propagation error mitigation strategy, has enabled the investigation of various levels of abstraction in surrogate modeling to enhance site-level optimization accuracy and efficiency \cite{maDatadrivenStrategiesOptimization2022}. Additionally, techniques such as optimality surrogates and DFO have been employed to address tractability challenges in large-scale, multi-level formulations, integrating hierarchical planning, scheduling, and control decisions within chemical companies \cite{vandebergHierarchicalPlanningschedulingcontrolOptimality2023}.
\\[8pt]
\textbf{Design of Experiments:} A recent paper illustrates high-throughput screening for \textit{Catalytically Active Inclusion Bodies} (CatIBs), utilizing a semi-automated cloning workflow and Bayesian optimization to efficiently generate and screen 63 glucose dehydrogenase variants from Bacillus subtilis, reducing manual effort and enhancing reproducibility \cite{helleckesHighthroughputScreeningCatalytically2024}. Similarly, Bayesian optimization has been utilized for the development of computation-driven materials discovery workflows, focusing on the exploration of unchartered material space \cite{mrozUnknownHowComputation2022}. Another example explores multi-objective Bayesian optimization in flow reactor experiments to identify the Pareto front of the optimal solution – using the \textit{q-Noisy Expected Hypervolume Improvement} (qNEVI) acquisition function \cite{zhangMultiobjectiveBayesianOptimisation2024}.  Additionally, \textit{Piecewise Affine Surrogate-based optimization} (PWAS) has been applied to tackle experimental planning challenges, such as optimizing Suzuki–Miyaura cross-coupling reaction conditions – benchmarked against genetic algorithms and Bayesian optimization variants \cite{zhuDiscreteMixedvariableExperimental2024}.
\\[8pt]
\textbf{Process Dynamics \& Control:} Surrogate DDO finds diverse applications within process control systems, for example, to create surrogates of the dynamic system by time-series modeling \cite{bradfordStochasticDatadrivenModel2020}. In a different work, artificial neural networks with rectifier units accurately represent piecewise affine functions for linear time-invariant systems within \textit{Model Predictive Control} (MPC) \cite{kargEfficientRepresentationApproximation2020}. Another study enhances computational efficiency, noise resilience, and control action smoothness in optimal dynamic product transitions using a data-driven Bayesian approach \cite{flores-tlacuahuacDataDrivenBayesian2024}. Additionally, an innovative application of Online Feedback Optimization with Gaussian Process regression mitigates plant-model mismatch in compressor station operations, resulting in reduced power consumption despite incomplete plant knowledge \cite{zagorowskaOnlineFeedbackOptimization2023}.
\\[8pt]
\textbf{Methodological Developments:}
Finally, looking at methods-based applications of surrogate-based optimization, a recent publication introduces a new data-driven optimization algorithm using \textit{Support Vector Machines} (SVMs) to tackle numerical infeasibilities within \textit{Differential Algebraic Equations} (DAEs), showcasing its effectiveness across diverse case studies, including complex scenarios in reaction engineering, such as the thermal cracking of natural gas liquids \cite{beykalDataDrivenOptimization2020a}. Another approach uses decision-focused surrogate modeling, which aims to address computationally challenging nonlinear optimization problems in real-time settings, validated through nonlinear process case studies such as reactors and heat exchangers \cite{guptaDataDrivenDecision2024}. Surrogate models also find applications in global optimization, exemplified by a novel algorithm tailored for solving linearly constrained mixed-variable problems, employing a piecewise affine surrogate of the objective function alongside an exploration function utilizing \textit{Mixed-Integer Linear Program} (MILP) solvers to search the feasible domain \cite{zhuGlobalPreferencebasedOptimization2023}. 
\\[8pt]
Four works explore optimization methodologies. One employs surrogate-based branch-and-bound algorithms for simulation-based optimization, ensuring consistent convergence to optimal solutions despite variability in initialization, sampling, and surrogate model selection \cite{zhaiSurrogatebasedBranchandboundAlgorithms2023}. \cite{rebelloAssuringOptimalitySurrogatebased2024} develop a robustness test, which ensures the optimality of results derived from surrogate models by guaranteeing they adhere to the universal approximation theorem. Another develops a method for Bayesian optimization of mixed-integer nonlinear programming problems \cite{morlet-espinosaBayesianOptimizationApproach2024}. A final example is a data-driven coordination framework for enterprise-wide optimization which maintains organizational autonomy while outperforming conventional distributed optimization methods across various case studies \cite{vandebergDataDrivenCoordination2023}.
\\[8pt]
As it is clear from this section, the process systems engineering community has a rich history of utilizing and advancing model-based (surrogate-based) DFO techniques.
\\[8pt]
This book chapter has two main objectives. The first objective is to serve as an introduction to newcomers to the field, to introduce the cornerstone concepts behind model-based derivative-free methods, commonly termed surrogate-based optimization, as well as to present the main advantages, shortcomings, and general framework. A second objective of this book chapter is to offer a framework for objectively comparing different model-based derivative-free methods, and particularly focus on process systems engineering applications. While an exhaustive and all-encompassing benchmark is beyond the scope of this book chapter, benchmarking results have been provided for a selection of functions, both in the constrained and unconstrained case, as well as process systems engineering examples. A particular effort has been made to share all the associated code from this benchmarking, outlining its use, in the appended GitHub repository. The authors hope this will inspire practitioners to use the code base to compare different algorithms to meet their specific goals.

%%%%%%%%%%%%%%%%%
% Model-Based Derivative-Free Optimization
%%%%%%%%%%%%%%%%%
\newpage
\section{Model-Based Derivative-Free Optimization}\label{sec:mb dfo}
As introduced in Section \ref{sec:mbintro}, model-based methods utilize surrogate models of the objective function to inform their updates, thereby guiding the optimization process by only leveraging sampled data without relying on explicit derivative information\cite{larsonDerivativefreeOptimizationMethods2019}. However, sampling in many cases does not yield the true function value $f(\mathbf{x})$, because of measurement noise. The feedback of the sampled system is represented by $y$, which consists of the objective function value $f(\mathbf{x})$ and may be blurred by measurement noise $\varepsilon$ as shown in Equation \ref{eq:basic_sample_eq}.
\begin{equation} \label{eq:basic_sample_eq}
    y = f(\mathbf{x}) + \varepsilon
\end{equation}
It is generally assumed that $\varepsilon$ is a random variable, and in many cases, this assumption is extended to be normally distributed i.e., $\varepsilon \sim \mathcal{N}(0,\sigma_\varepsilon)$.
The sampled data set $\mathcal{D}$ containing $n_d$ sampled positions and corresponding system feedback $(\mathbf{x}, y)$ are depicted in this chapter as $\mathcal{D} = \{(\mathbf{x}^{(1)}, y^{(1)}), (\mathbf{x}^{(2)}, y^{(2)}), ..., (\mathbf{x}^{(n_d)}, y^{(n_d)}) \}$. This leaves $X \in \mathbb{R}^{n_x \times n_d}$ as the matrix containing the data from the inputs (decision variables), and $\mathbf{y} \in \mathbb{R}^{n_d}$ as the vector containing the data from the output (sampled objective function values possibly corrupted by noise):
\begin{equation*}
    \mathbf{x} = 
        \begin{bmatrix}
            x_1 \\
            x_2 \\
            \vdots \\
            x_{n_x}
        \end{bmatrix}, 
    \quad
    X = \begin{bmatrix}
                        (\mathbf{x}^{(1)})^\intercal \\
                        (\mathbf{x}^{(2)})^\intercal \\
                        \vdots \\
                        (\mathbf{x}^{(n_d)})^\intercal
                    \end{bmatrix}
                    =
                    \begin{bmatrix}
                        x_1^{(1)} & \dots & x_{n_x}^{(1)} \\
                        x_1^{(2)} & \dots & x_{n_x}^{(2)} \\
                        \vdots & \ddots & \vdots \\
                        x_1^{(n_d)} & \dots & x_{n_x}^{(n_d)}
                    \end{bmatrix}
                    ,
    \quad
    \mathbf{y} = \begin{bmatrix}
                    y^{(1)}\\
                    y^{(2)}\\
                    \vdots
                    \\
                    y^{(n_d)}
                    \end{bmatrix}
\end{equation*} 

For simplicity, the output dimension is set as $n_y = 1$. However, this is revisited later on in the case of constrained problems.

% \\[8pt]
In most model-based optimization algorithms, two sequential optimizations are conducted. The first one is responsible for the surrogate (model) building step, which generally relies on building a surrogate (model) function that minimizes a likelihood function between the data and the surrogate function (commonly a least squares). This creates an approximation of $f$ denoted as $\hat{f}$. The second optimization optimizes $\hat{f}$ to find the next best candidate. Solely optimizing the surrogate function $\hat{f}$ might lead to over-exploitation, and not enough exploration, and therefore algorithms include some exploration components into their routine. Different algorithms propose different ways to add this exploration element, for example, Bayesian optimization, which uses Gaussian Processes as surrogates, leverages the prediction of the uncertainty to sample points that are promising but also unknown. Other algorithms, directly sample points that are meant to be used in the surrogate building step to obtain more information from the function $f$. The subsequent sections describe the different algorithms and how they handle this exploration-exploration dilemma. 

%The primary focus is on optimizing the given objective function. However, alongside this, there's a secondary optimization of the parameter set for the surrogate model. This model parameter optimization aims to fit the surrogate model with the data obtained from function evaluations at each iteration. 
In summary, at each iteration, a surrogate function is created through some optimization procedure, and then, this surrogate function is optimized to find the best candidate point(s) for the next iteration. It should also be noted that DFO does not preclude the use of conventional derivative-based solvers for the optimization of surrogate models \cite{vandebergDatadrivenOptimizationProcess2022}. Notably, for model-based DFO to prove effective, the building of the surrogate and its subsequent optimization should be computationally less expensive than sampling the true function \cite{larsonDerivativefreeOptimizationMethods2019b}. 
\\[8pt]

The choice of the most appropriate surrogate model is contingent upon the specific characteristics of the system at hand, and a variety of surrogates can achieve state-of-the-art performance, some notable examples:
\begin{itemize}
    \item Algorithms that construct quadratic surrogates, such as COBYQA \cite{ragonneauModelBasedDerivativeFreeOptimization2022} and CUATRO \cite{vandebergDatadrivenOptimizationProcess2022}.
    \item Approaches utilizing Gaussian Processes, as seen in Bayesian optimization \cite{thegpyoptauthorsGPyOptBayesianOptimization2016}.
    \item Techniques relying on Radial Basis Functions such as DYCORS and SRBFStrategy \cite{krityakierneSOPParallelSurrogate2016, regisStochasticRadialBasis2007a}.
    \item The use of decision trees as surrogates such as in ENTMOOT \cite{thebeltENTMOOTFrameworkOptimization2021}.
    \item Strategies that leverage basis functions to construct surrogates \cite{huyerSNOBFITStableNoisy2008,wilsonALAMOApproachMachine2017, boukouvalaARGONAUTAlgoRithmsGlobal2017}.
\end{itemize}

In the next section, the general workflow for surrogate-based optimization is presented.

\subsection{Workflow for Model-Based DFO}
\begin{figure}[htbp]
    \centering
    \includegraphics[scale=0.75]{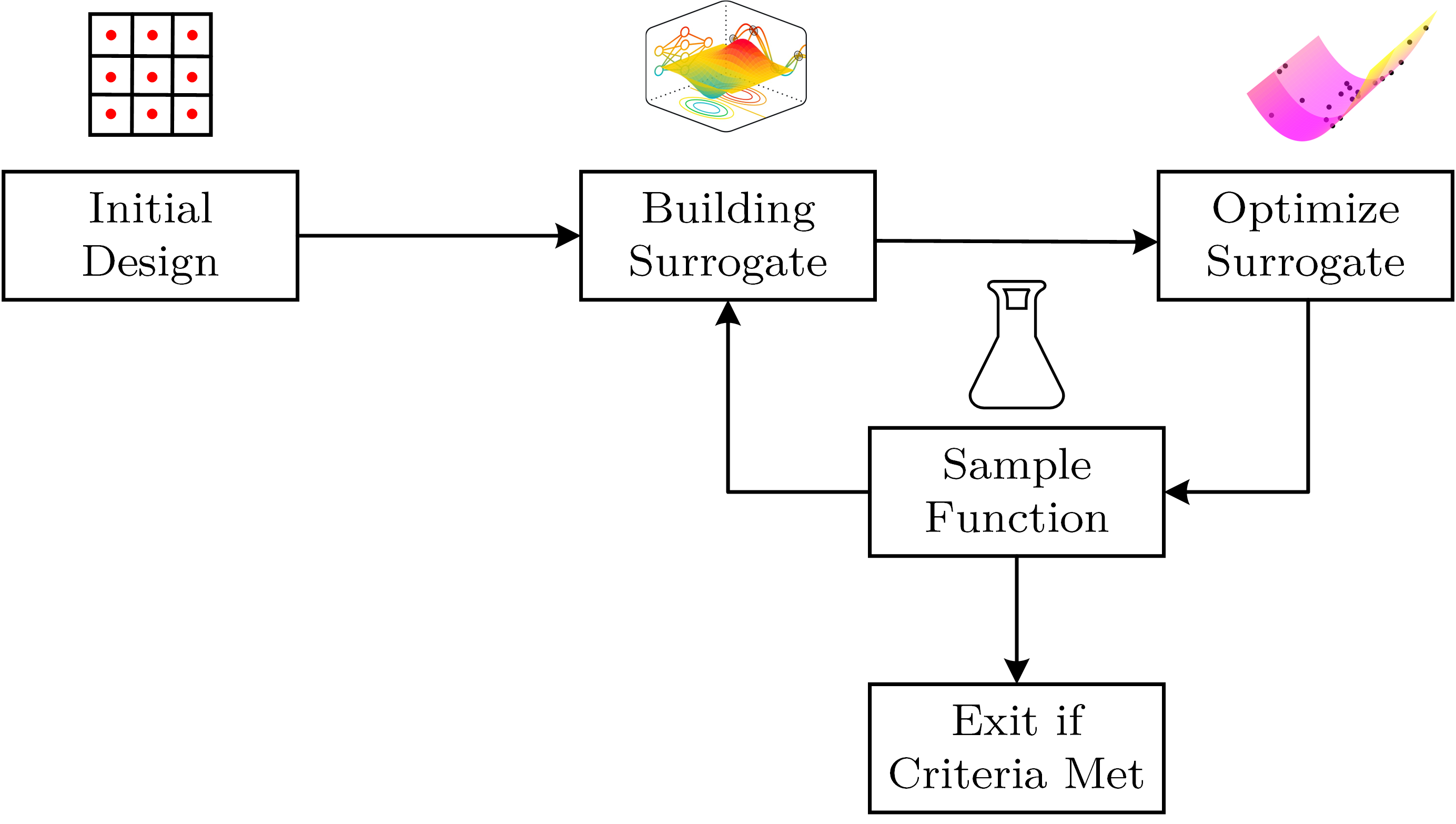}
    \caption{Fundamental workflow of model-based DFO.}
    \label{fig:Image_opti_flow}
\end{figure}

The fundamental workflow of model-based DFO is illustrated in Figure \ref{fig:Image_opti_flow}, with further details on each stage below:

\begin{itemize}
    \item \textbf{Initial Design:} This stage entails the initial set of experiments or samples aimed at gathering essential information about the objective function. Techniques commonly employed include DoE, \textit{Latin hypercube sampling} (LHS), factorial designs, or space-filling designs in general which explore the input space to capture its behavior and characteristics.
    
    \item \textbf{Build Surrogate:} This stage uses the available data to build a surrogate model of the objective function. Surrogates can take various forms, including machine learning techniques such as Gaussian Processes, decision trees, neural networks, etc. It is noteworthy that surrogates do not have to be machine learning models, for example, quadratic surrogates, which approximate the objective function using quadratic polynomials can be used. For instance, algorithms like COBYQA utilize quadratic surrogates to iteratively optimize the objective function. 
    
    \item \textbf{Optimize Surrogate:} This step uses the surrogate to find the next best point(s) to sample. One straightforward option is to optimize the surrogate (a model of the objective function) and sample there. However, more often than not the surrogate is not simply optimized and sampled, but some exploration component is added. This involves a trade-off between exploration (sampling in unexplored regions) and exploitation (sampling in regions that are likely to contain the optimal solution).

    \item \textbf{Sample Function:} This phase involves querying the 'true' system and evaluating the objective function to obtain new data. These samples are later used to update the surrogate, and at the same time represent optimal (or near optimal) solutions at termination.
    
    \item \textbf{Termination:} Termination occurs when the optimization problem-specific criterion is met, such as reaching a predefined number of iterations, achieving a satisfactory solution, or exhausting computational resources.
\end{itemize}

\subsection{Local vs Global Surrogates}
Surrogate optimization approaches can be broadly classified into local and global approaches. Global approaches proceed by constructing a surrogate model based on all their samples and constructing a single (flexible) surrogate for the whole decision space. The optimization of this surrogate and the next sample is also allowed to be anywhere within the problem constraints, without it having to be close to the current (or any previous) sampled point. Several practical implementations rely on neural networks \cite{henaoSurrogateBasedSuperstructure2011, savageMultifidelityDatadrivenDesign2023}, GPs \cite{caballeroAlgorithmUseSurrogate2008, chanonaRealtimeOptimizationMeets2021}, RBFs \cite{gutmannRadialBasisFunction2001, costaRBFOptOpensourceLibrary2018}, or a combination of basis functions \cite{wilsonALAMOApproachMachine2017, boukouvalaARGONAUTAlgoRithmsGlobal2017} to create the surrogate model.
\\[8pt]
By contrast, local approaches seek to maintain an accurate approximation of the original optimization problem within a trust region, whose position and size are adapted iteratively. This procedure entails updating or reconstructing the surrogate model as the trust region moves around, and it benefits from a well-developed convergence theory providing sufficient conditions for local optimality in unconstrained and bound-constrained problems \cite{connGlobalConvergenceGeneral2009}. However, it is important to highlight that to guarantee convergence these methods rely on a substantial number of function evaluations, which might be too expensive for costly optimization problems, and therefore are rarely used in expensive black-box optimization problems.

Handling constraints for DFO problems is still an active field of research \cite{audetProgressiveBarrierDerivativefree2018, gardnerBayesianOptimizationInequality,boukouvalaDerivativeFreeOptimization2014,easonAdvancedTrustRegion2018}, and in some instances trust regions have been found beneficial, particularly in safe exploration \cite{petsagkourakisSafeRealTimeOptimization2021,sorekDimensionalityReductionProduction2017}.
\\[8pt]
In the following sections, the main ideas behind different surrogate-based optimization algorithms are presented. These are intended to give a high-level understanding to the audience, and the interested reader is encouraged to go into additional material and references for a complete treatment of the different methodologies.

\subsection{Bayesian Optimization (BO)}

\textit{Bayesian Optimization} (BO) is a surrogate-based optimization strategy that relies on a probabilistic model as the surrogate. This probabilistic surrogate is used to predict the outcome of the objective function at unobserved points and then uses an acquisition function to determine the next point to evaluate. This approach balances exploration (trying points where little is known about the function) and exploitation (trying points where the model predicts a high value).
\\[8pt]
The process begins with initial observations of the objective function, which are used to build a probabilistic model, often a \textit{Gaussian Process} (GP). This model is then used to predict the outcome of the function at any untried point. The acquisition function, such as Expected Improvement, Probability of Improvement, or Upper Confidence Bound, uses these predictions to decide which point to evaluate next. This point is then evaluated on the actual objective function, and the result is used to update the probabilistic model. This cycle continues until a satisfactory solution is found or the computational budget is exhausted. As the reader can appreciate, BO follows the standard model-based DFO framework that is seen in Figure \ref{fig:Image_opti_flow}.

\subsubsection{BO - Surrogate Model}
BO relies on a probabilistic model, usually a GP. A GP is a stochastic process that, instead of defining a specific function, specifies a distribution over functions, where any finite set of function values have a joint Gaussian distribution. This makes GPs particularly useful for tasks like regression, interpolation, uncertainty quantification, and Bayesian optimization.
\\[8pt]
A GP is fully characterized by its mean function $m(\mathbf{x})$ and covariance function (or kernel) $k(\mathbf{x}, \mathbf{x}')$. For a given set of input points $X$, the corresponding GP is defined as:

\begin{itemize}
    \item Mean function: $m(\mathbf{x})$ describes the expected value of the process at any point $\mathbf{x}$. Typically, it is assumed to be zero, assuming that the dataset is standardized. This is the value that the GP will predict as an expected value far away from any data point.
    \item Covariance function: $k(\mathbf{x}, \mathbf{x}')$ specifies how the process covaries between any two points $\mathbf{x}$ and $\mathbf{x}'$. It captures the similarity between inputs and determines how much influence nearby points have on each other. Common choices include the squared exponential kernel, Matérn kernel, and rational quadratic kernel.
\end{itemize}

Given that a GP is a probabilistic model of the objective function $f \sim \mathcal{GP}(m(\cdot), k(\cdot,\cdot))$, once a set of data points, $\mathcal{D}=\{X,\mathbf{y} \}=\{(\mathbf{x}^{(i)},y^{(i)})\}_{i=1}^{n_d}$, has been sampled, including input points $X$ and their corresponding objective function values $\mathbf{y}$, the joint distribution between the sampled points and the objective function value to be predicted at a new point $\mathbf{x}^{(\text{new})}$ can be defined by Equation \ref{eq:BO_jd}.
\begin{equation}\label{eq:BO_jd}
p(f,\mathbf{y})= \mathcal{GP} \left( \begin{bmatrix} f \\ \mathbf{y} \end{bmatrix};\begin{bmatrix} m(\mathbf{x}^{(\text{new})}) \\ \mathbf{m}(X) \end{bmatrix}, \begin{bmatrix} K(\mathbf{x}^{(\text{new})},\mathbf{x}^{(\text{new})}) & k^\intercal(X,\mathbf{x}^{(\text{new})}) \\ k(X,\mathbf{x}^{(\text{new})}) & K(X,X)+ \sigma_n^2I \end{bmatrix} \right)
\end{equation}

Where $\mathbf{m}(X)=[m(\mathbf{x}^{(1)}),...,m(\mathbf{x}^{(n_d)})]^\intercal$ is the mean function evaluated at each datapoint, $\sigma_n^2I$ is the observation noise covariance, and $K$ is the covariance of the samples. It is therefore possible to find the probability of $f$ given the data by marginalization:
\begin{equation}\label{eq:BO_prob}
p(f|\mathbf{y})=\mathcal{GP}(f;\mu_\mathcal{D}(\cdot), \sigma^2_\mathcal{D}(\cdot,\cdot))
\end{equation}

Where $\mu_\mathcal{D}$ is the expected value and $\sigma^2_\mathcal{D}$ is the variance of the prediction of $f$ at a new point $\mathbf{x}^{(\text{new})}$:
\begin{equation}
\mu_\mathcal{D} = m(\mathbf{x}^{(\text{new})}) + k^\intercal(X,\mathbf{x}^{(\text{new})}) [K(X,X)+ \sigma_n^2I]^{-1}(\mathbf{y}-\mathbf{m}(X))   
\end{equation}
\begin{equation}
\sigma^2_\mathcal{D} = K(\mathbf{x}^{(\text{new})},\mathbf{x}^{(\text{new})}) - k^\intercal(X,\mathbf{x}^{(\text{new})})[K(X,X)+ \sigma_n^2I]^{-1}k(X,\mathbf{x}^{(\text{new})})  
\end{equation}

\subsubsection{BO - Surrogate Optimization} \label{sec:BO_so}
Given a GP that models the objective function, BO uses this model to pick the next point. One straightforward way to conduct this would be to minimize the expected value of the predicted function, i.e., $\min_{\mathbf{x}} \mu_\mathcal{D}(\mathbf{x})$. However, given that the objective function is modeled by a probabilistic model (e.g., a GP) it is possible to use the uncertainty predicted by the model to search for the next point to sample. 
\\[8pt]
In BO the \emph{acquisition function}  $\mathcal{A}$ is used to determine the next point to evaluate in the search space. It quantifies the utility or potential of evaluating a point $\mathbf{x}$ based on the current probabilistic model. Popular acquisition functions include \textit{Probability of Improvement} (PI), \textit{Expected Improvement} (EI), and \textit{Lower Confidence Bound} (LCB). The choice of acquisition function depends on the trade-off between exploration (sampling in regions of uncertainty) and exploitation (sampling where the surrogate model predicts high values).
\\[8pt]
One easy-to-understand acquisition function is the Lower Confidence Bound, where at each iteration, after the GP is built the surrogate optimization step solves the following problem:
\begin{equation}\label{eq:GP_LCB}
\min_\mathbf{x} \mathcal{A}^{lcb}(\mathbf{x}) = \min_\mathbf{x} \mu_\mathcal{D}(\mathbf{x}) - \gamma \sigma_\mathcal{D}(\mathbf{x})   
\end{equation}

where $\gamma$ is a hyperparameter of the algorithm that determines how much exploration versus exploitation is preferred. From an empirical perspective, it can be observed that this acquisition function minimizes the expected value (i.e., $\mu_\mathcal{D}(\mathbf{x})$) while at the same time encouraging the exploration of regions with uncertainty that could yield an even better solution (i.e., the term $\gamma \sigma_\mathcal{D}(\mathbf{x})$ ).
\\[8pt]
In this way by iteratively updating the surrogate model (e.g., a GP) and selecting points to evaluate based on minimizing the acquisition function, BO explores the search space and converges to optimal or near-optimal solutions while minimizing the number of evaluations of the true objective function.

\subsubsection{BO - Algorithms}
This book chapter uses various BO algorithms in its benchmarking procedure. Two of these are standard variants of BO, the first, GPyOpt \cite{thegpyoptauthorsGPyOptBayesianOptimization2016}, is an off-the-shelf implementation based on \cite{gonzalezBatchBayesianOptimization, gonzalezBayesianOptimizationSynthetic2015, gonzalezGLASSESRelievingMyopia} and the second is an in-house implementation. Additionally, the state-of-the-art TuRBO (Trust Region BO) algorithm was included in the high-dimensional case study \cite{erikssonScalableGlobalOptimization2019}. One key challenge of the standard version of BO with GPs is its efficacy in high-dimensional settings. This shortcoming arises because standard BO uses a single, uniform model to represent the entire search space. Such a model fails to capture the diverse characteristics of different regions and overly focuses on exploring new areas across the entire search space, rather than efficiently targeting the most promising regions. The TuRBO algorithm addresses these issues by adopting a local probabilistic approach for global optimization in large-scale, high-dimensional problems. Instead of relying on a single global model, TuRBO fits multiple local models within dynamically adjusted trust regions, focusing on promising areas of the search space. This method leverages an implicit bandit approach to allocate samples efficiently across these local models.

\subsection{Ensemble Tree MOdel Optimization Tool (ENTMOOT)}
The \textit{Ensemble Tree MOdel Optimization Tool} (ENTMOOT) is a framework designed to perform BO using \textit{Decision Tree} (DT) based surrogate models \cite{thebeltENTMOOTFrameworkOptimization2021}. In each iteration, ENTMOOT approximates the black-box function using a \textit{Gradient-Boosted Decision Tree} (GBDT) model from LightGBM \cite{keLightGBMHighlyEfficient}. GBDTs are captured in a basic fashion in the subsequent section, followed by a description of their embedding into the BO framework.

\subsubsection{ENTMOOT - Surrogate Model}
DTs are supervised learning algorithms that can be interpreted as a piece-wise constant approximation of the underlying function \cite{breimanClassificationRegressionTrees1984}, in this case, the objective function. As the name suggests DTs have a tree-like structure composed of nodes, leaves, and branches. The topmost node in a DT is the root node representing the entire dataset. The root node is split into decision nodes, representing the split of the data based on certain conditions. Terminal nodes are the leaf nodes, representing the final output and no further splitting occurs at these nodes. All nodes are connected through branches, which represent the outcome of a decision rule applied to a former node.

When adapted to a regression task as required for this study, a DT can be depicted as an optimization problem that embeds the search for optimal splitting parameters that minimize the mean-squared error loss function $\mathcal{L}^{MSE}$ at each node. Thereby the overall optimization problem for a DT yields
\begin{equation} \label{eq:dt_optim}
    \min_{T} \sum_{n \in \text{nodes}} \left( \frac{N_L(n)}{N_n} \mathcal{L}^{MSE}_L(n) + \frac{N_R(n)}{N_n} \mathcal{L}^{MSE}_R(n) \right)
\end{equation}
with $T$ representing the structure of the tree composed of $n$ nodes. Each node $n$ governs the split of $N_n$ samples with $N_L(n)$ and $N_R(n)$ being the numbers of samples in the left and right subsets for node $n$, respectively. $\mathcal{L}^{MSE}_L(n)$ and $\mathcal{L}^{MSE}_R(n)$ are the mean squared errors for the left and right subsets for node $n$ with $\mathcal{L}^{MSE}$ as
\begin{equation} \label{eq:dt_mse}
    \mathcal{L}^{MSE}_S = \frac{1}{|S|} \sum_{i \in S} (y_i - \bar{y}_S)^2    
\end{equation}
Here, $|S|$ is the number of samples in subset $S$, $y_i$ is the actual target value of the $i$-th sample, and $\bar{y}_S$ is the mean target value for subset $S$.

Gradient boosting in the context of GBDTs refers to the construction of a series of trees with a new tree added to the ensemble after each iteration. Each tree attempts to improve upon the performance of the ensemble so far, minimizing the error between objective function observation $y_i$ and prediction $\hat{f}_{\text{entmoot},i}$. 

The residuals $r_i$ are defined as the negative gradients of the loss function \ref{eq:dt_mse} 
\begin{equation}
    r_i = - \frac{\delta \mathcal{L}^{MSE}(y_i, \hat{f}_{\text{entmoot},i})}{\delta \hat{f}_{\text{entmoot},i}}
\end{equation}
The new tree is trained to predict these residuals. 

For each observation $i$ the residuals $r_i = y_i - \hat{f}_{\text{entmoot},i}$ are calculated and a new DT $h_m(x)$ is trained to predict the residuals from the previous step. The new tree's predictions are added to the ensemble:
\begin{equation} \label{eq:pred_dt}
\hat{f}_{\text{entmoot},i} = \hat{f}_{\text{entmoot},i}^{previous} + \nu h_m(x_i)
\end{equation}
with $\nu$ being the learning rate controlling the contribution of each tree. The final prediction is the sum of all contributions from trees $m \in M$:
\begin{equation}\label{eq:entmoot_pred}
\hat{f}_{\text{entmoot}} = \hat{f}_{\text{entmoot},0} + \sum_{m=1}^{M} \nu h_m(x)    
\end{equation}

%There is some information in the paper regarding how the author's managed to overcome the shortcomings of traditional implementations of GBDT: LightGBM, which contains two novel techniques: Gradient-based One-Side Sampling and Exclusive Feature Bundling

\subsubsection{ENTMOOT - Surrogate Optimization}
Following the BO paradigm, ENTMOOT leverages an acquisition that balances exploration and exploitation to determine the next sampled point. To do this, the surrogate model's prediction is combined with an uncertainty measure, reflecting varying degrees of trust in the data points. Since DTs do not come with an uncertainty quantification as GPs do, \cite{thebeltENTMOOTFrameworkOptimization2021} use the uncertainty measure $\alpha(x)$ to quantify the model uncertainty using the distance to the closest point $\mathbf{x}_d$ in data set $\mathcal{D}$ to quantify the confidence of model predictions.
\begin{equation} \label{eq:entmoot_alpha}
    \alpha(\mathbf{x})=\min_{d \in \mathcal{D}}||\mathbf{x}-\mathbf{x}_d||_p
\end{equation}
where $p$ determines the metric used to determine the distance. \cite{thebeltENTMOOTFrameworkOptimization2021} use implementations of Euclidean or Manhattan distance metrics and discuss advantages and shortcomings in their paper.

This combination forms the acquisition function, which is then optimized to identify the best candidate for the optimal point in the current iteration. The following acquisition function represents the lower confidence bound acquisition function as seen in Equation \ref{eq:GP_LCB} adjusted to the elements of GBDTs. It is a simplified version as presented in \cite{thebeltENTMOOTFrameworkOptimization2021} that captures the main elements:
\begin{equation} \label{eq:entmoot_aq}
    \min_{\mathbf{x}} ~ \mathcal{A}^{\text{lcb, entmoot}}(\mathbf{x}) = \min_{\mathbf{x}} ~ \hat{f}_{\text{entmoot}}(\mathbf{x}) - \gamma \alpha(\mathbf{x})
\end{equation}
$\hat{f}_{\text{entmoot}}$ refers to the tree model prediction in Equation \ref{eq:entmoot_pred}, capturing how Equation \ref{eq:entmoot_aq} exploits the underlying surrogate model to find promising areas in the search space. Variable $\alpha(\mathbf{x})$ as introduced in Equation \ref{eq:entmoot_alpha} handles exploration and quantifies the degree of uncertainty expected from prediction $\hat{f}_{\text{entmoot}}$. $\gamma \in \mathbb{R}$ balances exploitation and exploration and is a hyperparameter depending on the application as described previously in Section \ref{sec:BO_so}.

Full details of this algorithm can be found in \cite{thebeltENTMOOTFrameworkOptimization2021}. 

\subsection{Constrained Optimization by Linear Approximation (COBYLA)}\label{sec:cobyla}

First introduced in 1994 by M.J.D. Powell, a pioneer of computational mathematics \cite{buhmannMichaelPowell292018}, the Constrained Optimization BY Linear Approximations (COBYLA) algorithm iteratively refines solutions for nonlinearly constrained optimization problems by leveraging linear approximations of the objective (and any constraint) functions \cite{powellDirectSearchOptimization1994a}. COBYLA employs a trust region bound to limit changes to the variables. The algorithm dynamically adjusts the trust region radius based on predicted improvements to both the objective function and feasibility conditions. Additionally, it utilizes a merit function to compare the effectiveness of different variable vectors in improving the shape of the simplex, while ensuring adherence to constraints. 

\subsubsection{COBYLA - Surrogate Model}
This algorithm solves the following problem at every iteration to determine the next best point to sample:

\begin{equation}\label{eq:cobyla_1}
\begin{aligned}
& \min_{\mathbf{x} \in \mathbb{R}^{n_x} } && \hat{f}_\text{cobyla}(\mathbf{x})\\
& \text{s.t.} && \hat{c}_{_\text{cobyla},i}(\mathbf{x}) \geq 0, &&& \quad i=1,2,\ldots,m\\
\end{aligned}
\end{equation}

As mentioned, COBYLA utilizes a linear approximation to construct a surrogate of the objective function. Hence, Equation \ref{eq:cobyla_1}, presents the linear programming problem which is iteratively solved, yielding a vector of variables, $\mathbf{x}$. Where $\hat{f}_\text{cobyla}(\mathbf{x})$ represents the surrogate, and $f(\mathbf{x})$ is the true objective function. Similarly, $\hat{c}_{\text{cobyla},i}(\mathbf{x})$ is a set of unique linear functions that serve as surrogates to the true constraints, ${c}_{\text{cobyla},i}(\mathbf{x})$.

At every iteration, COBYLA models the objective and the constraint functions with linear interpolants, which consists of $n_x + 1$ points that
are updated along the iterations.

\subsubsection{COBYLA - Surrogate Optimization}
After constructing the initial surrogate, $\hat{f}_\text{cobyla}(\mathbf{x})$, and a given center point $\mathbf{x}_k$, the algorithm iteratively enhances the linear approximation by adjusting the trust region radius, $\rho$, and evaluating potential solutions to determine the next variable vector, $\mathbf{x}_*$. Initially, it is verified whether $\mathbf{x}_k$ is the optimal vertex and ensures the simplex is acceptable. Subsequently, the trust region condition on $\mathbf{x}_*$ is given by Equation \ref{eq:cobyla_tr}.
\begin{equation}\label{eq:cobyla_tr}
\| \mathbf{{x}_*} - \mathbf{x}_k \|_2 < \frac{1}{2} \rho
\end{equation}
If feasible, the surrogate is now minimized by $\mathbf{x}_*$, subject to the trust region condition and linear constraints ${c}_{\text{cobyla},i}(\mathbf{x})$. In cases where multiple possible $\mathbf{x}_*$ exist, the vector that yields the lowest value of $\| \mathbf{x}_* - \mathbf{x}_k \|_2$ is selected.
\\[8pt]
To deal with constraints the COBYLA algorithm uses a merit function. This merit function, denoted as $\Phi(\mathbf{x})$, combines both the objective function $f(\mathbf{x})$ and the constraint functions $c_i(\mathbf{x})$ into a single scalar value:
\begin{equation}
\Phi_\text{cobyla}(\mathbf{x}) = f(\mathbf{x}) + \left[ \max_{i= 1, \ldots , m} c_i(\mathbf{x})\right]_+
\end{equation}
The merit function serves a dual purpose: the first term, $f(\mathbf{x})$, ensures that solutions are evaluated based on their ability to optimize the objective function. The second term, $\left[ \max {c_i(\mathbf{x}) : i = 1, 2, \ldots , m} \right]_+$, captures the magnitude of constraint violations. By considering the maximum violation among all constraints, the merit function guides the optimization towards solutions that not only optimize the objective function but also satisfy the constraints as closely as possible.
\\[8pt]
During the optimization process, COBYLA aims to minimize the merit function $\Phi_\text{cobyla}(\mathbf{x})$ by iteratively adjusting the variable vector $\mathbf{x}$ within the trust region bounds and evaluating potential solutions. By minimizing the merit function, COBYLA effectively balances the trade-off between optimizing the objective function and satisfying the constraints, ultimately guiding the search toward feasible and optimal solutions.

\subsection{Constrained Optimization by Quadratic Approximation (COBYQA)}\label{sec:cobyqa}

\textit{Constrained Optimization by Quadratic Approximation} (COBYQA) builds on Powell's COBYLA, presenting a similar model-based DFO method, but instead utilizing quadratic approximations. COBYQA constructs quadratic models of objective and constraint functions using derivative-free symmetric Broyden updates, enabling efficient optimization without explicit derivatives. It dynamically adjusts its trust-region radius and incorporates a geometry-improving procedure to enhance numerical stability. Importantly, COBYQA strictly adheres to bound constraints, ensuring robustness in various engineering and industrial applications. Unlike alternative methods like SQPDFO, COBYQA directly handles inequality constraints without introducing additional slack variables, maintaining efficiency and consistency throughout optimization iterations.

\subsubsection{COBYQA - Surrogate Model}
In COBYQA, the quadratic surrogate is constructed by approximating the objective \( f(\mathbf{x}) \) and constraints \( c_i(\mathbf{x}) \) functions within a trust region defined by a radius \( \Delta_k \) around the current iteration \( \mathbf{x}_k \) - where \( k \) represents the current iteration. The quadratic model \( \hat{f}_\text{cobyqa}(\mathbf{x}) \) and \( \hat{c}_{\text{cobyqa},i}(\mathbf{x}) \), defined in Equation \ref{eq:cobyqa:fsurr} and \ref{eq:cobyqa:csurr} respectively, are formed using function evaluations at selected points within this region, capturing the curvature of the functions and providing a refined estimate of their behavior \cite{ragonneauModelBasedDerivativeFreeOptimization2022}.

\begin{equation}\label{eq:cobyqa:fsurr}
\begin{aligned}
\hat{f}_\text{cobyqa}(\mathbf{x}) &= \hat{f}(\mathbf{x}_k) + \hat{\mathbf{g}}_k^\intercal(\mathbf{x} - \mathbf{x}_k) + \frac{1}{2} (\mathbf{x} - \mathbf{x}_k)^\intercal \hat{B}_k (\mathbf{x} - \mathbf{x}_k) 
\end{aligned}
\end{equation}
\begin{equation}\label{eq:cobyqa:csurr}
\begin{aligned}
\hat{c}_{\text{cobyqa},i}(\mathbf{x}) &= \hat{c}_i(\mathbf{x}_k) + \nabla \hat{c}_i(\mathbf{x}_k)^\intercal (\mathbf{x} - \mathbf{x}_k) + \frac{1}{2} (\mathbf{x} - \mathbf{x}_k)^\intercal \nabla^2 \hat{c}_i(\mathbf{x}_k) (\mathbf{x} - \mathbf{x}_k), \quad i = 1, 2, \dots, m 
\end{aligned}
\end{equation}

Where \( \hat{\mathbf{g}}_k \in \mathbb{R}^{n_x} \) is the surrogate of the gradient of the objective function, \( \hat{B}_k \in \mathbb{R}^{n_x \times n_x} \) is a positive definite symmetric matrix representing the surrogate Hessian of the  objective function. Similarly, \( \nabla \hat{c}_i(\mathbf{x}_k) \) is the surrogate of the gradient of the \( i \)-th constraint function, and \( \nabla^2 \hat{c}_i(\mathbf{x}_k) \) is the surrogate of its Hessian matrix.

\subsubsection{COBYQA - Surrogate Optimization}
COBYQA optimizes the quadratic surrogate to identify the next iterate \( \mathbf{x}_{k+1} \) by minimizing a merit function \( \Phi_{\text{cobyqa},k}(\mathbf{x}) \) within the trust region. The merit function, defined by Equation \ref{eq:cobyqa_merit} combines the objective and constraint violations, weighted by penalty parameters, guiding the optimization process. The trust region, denoted by \( \mathcal{B}(\mathbf{x}_k, \Delta_k) \), is given by Equation \ref{eq:cobyqa_trust}.
\begin{equation}\label{eq:cobyqa_merit}
\Phi_{\text{cobyqa},k}(\mathbf{x}) = \hat{f}_{\text{cobyqa},k}(\mathbf{x}) + \sum_{i=1}^{m} \rho_i [\hat{c}_{i,\text{cobyqa},k}(\mathbf{x}) - \Delta_k]_+
\end{equation}
\begin{equation}\label{eq:cobyqa_trust}
\mathcal{B}(\mathbf{x}_k, \Delta_k) = \{ \mathbf{x} : \| \mathbf{x} - \mathbf{x}_k \| \leq \Delta_k \}
\end{equation}
Here, \( \rho_i \) represents penalty parameters associated with the constraint functions.

\subsection{Local Search with Quadratic Models (LSQM)}
\textit{Local Search with Quadratic Models} (LSQM) is a naive surrogate-based optimization method studied in this chapter. It constructs a quadratic surrogate model based on sampled data and optimizes it within a trust region. LSQM is advantageous as both the model construction and optimization are convex problems. However, the algorithm is inherently exploitative, hence it lacks any exploration element.

\subsubsection{LSQM - Surrogate Model}
LSQM builds surrogates of the form:
\begin{equation}
\hat{f}_\text{lsqm}(\mathbf{x};Q,c,b) = \mathbf{x}^\intercal Q\mathbf{x} + c^\intercal\mathbf{x} + b  
\end{equation}
where $Q \in \mathbb{R}^{n_x \times n_x}$ is symmetric, $c \in \mathbb{R}^{n_x}$, and $b \in \mathbb{R}$. Given the dataset $\mathcal{D}= \{ (\mathbf{x}^{(i)},y^{(i)}) \}_{i=1}^{n_d}$, where values obtained from the objective function are denoted as $y^{(i)} \leftarrow f(\mathbf{x}^{(i)})$. The following least squares problem is formulated to estimate $Q$, $c$ and $b$:
\begin{equation}
\min_{Q,c,b} \quad \sum_{i=1}^{n_d} \left( \hat{f}_\text{lsqm}(Q,c,b; \mathbf{x}^{(i)})  - y^{(i)}  \right)^2
\end{equation}
This is a convex optimization problem which can be easily solved. Furthermore, if semi-positive-definiteness is imposed on $Q$ (i.e., $Q \succeq 0$) then the ensuing surrogate optimization problem is also convex.

\subsubsection{LSQM - Surrogate Optimization} \label{sec:lsqm_surr_opt}
Once the quadratic surrogate model is constructed, the following optimization problem is formulated to find the next point to sample and evaluate:
\begin{equation}
\min_\mathbf{x} \quad \hat{f}_\text{lsqm}(\mathbf{x};Q,c,b) = \min_\mathbf{x} \quad \mathbf{x}^\intercal Q\mathbf{x} + c^\intercal\mathbf{x} + b
\end{equation}
N.B. The bias term $b$ can be omitted from the optimization as it will not influence the optimal point. In summary, LSQM presents a simple surrogate-based optimization algorithm that is used as a baseline in this chapter.

\subsection{Convex qUAdratic Trust-Region Optimizer (CUATRO)}\label{sec:cuatro}

CUATRO is a quadratic model-based trust region method first introduced in \cite{vandebergDatadrivenOptimizationProcess2022}. It works similarly to COBYQA and LSQM with some key distinctions. CUATRO is implemented in the convex programming framework CVXPY \cite{diamondCVXPYPythonEmbeddedModeling2016}, which allows explicit handling of surrogate constraints. Additionally, CUATRO is developed with the peculiarities of chemical engineering applications in mind: it includes heuristic routines for sample-efficient and safe exploration, as well as dimensionality reduction techniques to exploit latent structure in the solution space to scale to thousands of variables \cite{vandebergHighdimensionalDerivativefreeOptimization2024}.

\subsubsection{CUATRO - Surrogate Model}

Similarly to LSQM, CUATRO builds objective and constraint surrogates of the following form:
\begin{equation}
\hat{f}_\text{cuatro}(\mathbf{x};Q,\mathbf{p},r) = \mathbf{x}^\intercal Q\mathbf{x} + \mathbf{p}^\intercal\mathbf{x} + r
\end{equation}
\begin{equation}
\hat{c}_{\text{cuatro},i}(\mathbf{x}; Q_i,\mathbf{p}_i,r_i) = \mathbf{x}^\intercal Q_i\mathbf{x} + \mathbf{p}_i^\intercal\mathbf{x} + r_i 
\end{equation}

The objective $\hat{f}(\cdot)$ is again trained via regular least-squares regression. This time, however, convexity is enforced by constraining $Q$ to be positive semi-definite:

\begin{equation} \label{eq:cuatro_lsqr}
\min_{Q \succeq 0,\mathbf{p},r} \quad \sum_{i=1}^{n_d} \left( \hat{f}_\text{cuatro}(\mathbf{z}^{(i)}; Q,\mathbf{p},r)  - y^{(i)}  \right)^2
\end{equation}

Instead of training separate constraint surrogates using least squares regression on each set of constraint evaluations, the default implementation of CUATRO performs convex quadratic discrimination: it finds the ellipsoid that minimizes the total distance to the discrimination boundary of all falsely classified samples. The interested reader is referred to \cite{vandebergDataDrivenCoordination2023} for the convex formulation.

\subsubsection{CUATRO - Surrogate Optimization}

The surrogate optimization step follows that of COBYQA, with the exception that box bounds ($[\mathbf{x}_L, \mathbf{x}_U]$) and the constraint surrogate $\hat{c}(\cdot)$ are included as explicit constraints resulting in a convex semidefinite program:

\begin{equation}
\label{CUATRO minimisation}
\begin{aligned}
     & \underset{\mathbf{x} \in \mathcal{B}(\mathbf{x}_k, \delta_k) }{\text{min.}}
     &&  \hat{f}_\text{cuatro}(\mathbf{x}; \cdot) \\
     & \text{s.t.} 
     && \hat{c}_\text{cuatro}(\mathbf{x}; \cdot) \leq 0 \\
     &&& \mathcal{B}(\mathbf{x}_k, \delta_k) = \{ \mathbf{x} \in [\mathbf{x}_L, \mathbf{x}_U] \text{ } | \text{ } ||\mathbf{x} - \mathbf{x}_k||_F^2 \leq \delta_k^2 \} \\
\end{aligned}
\end{equation}

\subsubsection{CUATRO - Improving High-Dimensional Performance}

\cite{vandebergHighdimensionalDerivativefreeOptimization2024} introduces CUATRO-pls, an extension to CUATRO that performs dimensionality reduction to find a lower-dimensional subspace over which to perform surrogate updates. At each iteration, \textit{Partial Least Squares regression} (PLS) identifies the linear projection matrix $M \in \mathbb{R}^{n_z \times n_x}$ and reconstruction matrix $R=(M^\intercal M) M ^\intercal$ that best predict the outputs $\mathbf{y}$ of all samples within the original trust region $\mathcal{B}(\cdot)$ after projection to their linear embedding $Z$ such that $Z = M X$ and $X \approx R Z$. The surrogates are then trained in the reduced dimensional dataset $Z$:

\begin{equation}
\min_{Q \succeq 0,\mathbf{p},r} \quad \sum_{i=1}^{n_d} \left( \hat{f}_\text{cuatro}(\mathbf{z}^{(i)}; Q,\mathbf{p},r)  - y^{(i)}  \right)^2
\end{equation}

Surrogate fitting and minimization are also performed in the subspace. To conserve convexity, the equivalent trust region constraint in the original space $(R \mathbf{z} - \mathbf{x}_k)^\intercal (R \mathbf{z} - \mathbf{x}_k) \leq r^2$ is replaced with a heuristic trust region in embedding space centered around $\mathbf{z}_k=M \mathbf{x}_k$, where an effective radius is defined as the distance from $\mathbf{z}_k$ to the furthest projected sample as $\hat{\delta}^2 =  \max_{\mathbf{z} \in Z} \text{ } (\mathbf{z} - \mathbf{z}_k)^\intercal (\mathbf{z} - \mathbf{z}_k) $:

\begin{equation}
\label{eq: cuatro_min_pls}
\begin{aligned}
     & \underset{\mathbf{z} \in \mathcal{B}(\mathbf{z}_k, \hat{\delta}_k) }{\text{min.}}
     &&  \hat{f}_\text{cuatro}(\mathbf{z}; \cdot) \\
     & \text{s.t.} 
     && \hat{c}_\text{cuatro}(\mathbf{z}; \cdot) \leq 0 \\
\end{aligned}
\end{equation}

The minimization candidate is evaluated and the trust region is updated after reconstruction to the original space. This methodology introduces a crucial choice in the form of the embedding dimensionality $n_{pls}$. \cite{vandebergHighdimensionalDerivativefreeOptimization2024} show that a default underestimation of $n_{pls}=2$ works quite well in the absence of information about the intrinsic embedding. This makes intuitive sense, as the method reduces to a search in the most informative linear combination of dimensions, which would have merit as an idea on its own.

\subsection{Stable Noisy Optimization by Branch and Fit (SNOBFIT)}
\textit{Stable Noisy Optimization by Branch and Fit} (SNOBFIT) by \cite{huyerSNOBFITStableNoisy2008} is a surrogate-based optimization algorithm that combines global and local search. To do so SNOBFIT branches the search space to create smaller sub-spaces. Subsequently, surrogate models are fitted within the sub-spaces to obtain information about promising areas of the objective function. Candidate points generated in these sub-spaces are sorted into classes indicating, whether a local or global aspect of SNOBFIT has determined the respective subspace. There are three classes for local candidates and two classes for global candidates. Technically, "local" aspects refer to an exploitative strategy leveraging so-called \textit{save guarded nearest neighbors} (SNNs) to promising function evaluations. Accordingly, "global" aspects refer to an explorative strategy applied to gain information about unexplored sub-spaces. The subsequent part of this section starts by explaining how surrogate models are built to then delve into the usage of local and global information to navigate the optimization of the objective function. 

\subsubsection{SNOBFIT - Surrogate Model}
For fast local convergence, SNOBFIT handles local search from the location $\mathbf{x}^{best}$ of the best-so-far objective function value with a full quadratic model. 
To create such a model, the number of points in a local area $N$ must exceed the dimensions of the problem $n$ by $\Delta n$, $N \geq n + \delta n$. $\Delta n$ are called the previously introduced SNNs and can go up to $n$. The SNNs are technically a set of previous function evaluation points that are consulted as support points for the intended local surrogate model fit. The procedure to determine the SNNs for a point uses coordinate-wise comparison of previous function evaluation positions guided by a threshold to ensure diversity and to promote exploration. After each iteration of SNOBFIT, the SNNs for the locations of new evaluations are determined and the SNNs for previous evaluated locations are updated. For an in-depth description of the SNN determination the interested reader is referred to Section 3 in \cite{huyerSNOBFITStableNoisy2008}.

Subsequently, a local model around each point $\mathbf{x}$ is fitted with the aid of the SNNs:
\begin{equation}
    q(x) := f_{best} + \mathbf{g}^\intercal (\mathbf{x} - \mathbf{x}^{best}) + \frac{1}{2} (\mathbf{x} - \mathbf{x}^{best})^\intercal G(\mathbf{x} - \mathbf{x}^{best})
\end{equation}
which corresponds to a Taylor approximation around the best observed objective function value $f_{best}$. Gradient $\mathbf{g}$ and Hessian $G$ are determined based on the SNNs as described previously. For a more in-depth description of this approach, the reader is referred to \cite{huyerSNOBFITStableNoisy2008}. Depending on the user's preferences regarding the number of local optimization candidates, the previously described approach can be repeated, with adaptations described in \cite{huyerSNOBFITStableNoisy2008}.

\subsubsection{SNOBFIT - Surrogate Optimization}
Different from other algorithms previously discussed, SNOBFIT generates a batch of new candidate locations rather than only one. These candidate locations are generated by optimizing the surrogate models around the best-so-far evaluation from the previous section. The optimization of the quadratic surrogate model around this point is described in Section \ref{sec:lsqm_surr_opt}. Depending on the user's preference, the number of such exploitative generated candidates can be increased to include more candidates from SNN-supported local surrogates. 

The global aspect of SNOBFIT fills the remaining spots in the batch of candidate points for the succeeding iteration. Such candidates are generated from so far unexplored regions of the search space. For a box $[x_l,x_u]$ with corresponding point previously observed point $x$ the candidate point $z$ is generated by 
\begin{equation} \label{eq:snobfit_branch}
    z := 
        \begin{cases} 
            \frac{1}{2} (x_l + x) & \text{if } x - x_l > x_u - x \\ 
            \frac{1}{2} (x + x_u) & \text{otherwise},
        \end{cases}
\end{equation}
If after applying Equation \ref{eq:snobfit_branch} there are still spots in the batch to be filled by global candidates, random sampling within the bounds $[x_l,x_u]$ is used. Typically this happens during the initialization of the algorithm or the early iterations of the algorithm when not enough points are yet observed to generate the local quadratic models.

\subsection{DYnamic COordinate search using Response Surface models (DYCORS)}
This chapter considers three algorithms that are based on RBF surrogates by Shoemaker et al.: DYCORS, SOPStrategy, and SRBFStrategy \cite{regisCombiningRadialBasis2013}. RBFs are versatile mathematical tools used in surrogate modeling, characterized by their radial symmetry around a center point and employed to approximate complex functions. These functions are flexible, capable of capturing non-linearities, and adapt well to irregularly sampled data points, making them suitable for various applications in optimization and machine learning \cite{RadialBasisFunction}. While this section focuses extensively on detailing DYCORS, the interested reader is referred to the literature for comprehensive insights into SOPStrategy \cite{krityakierneSOPParallelSurrogate2016} and SRBFStrategy \cite{regisStochasticRadialBasis2007a, regisParallelStochasticGlobal2009}.
\\[8pt]
The DYCORS framework represents a sophisticated approach for surrogate-based optimization, tailored for high-dimensional, expensive, and black-box functions. DYCORS integrates elements from the \textit{Dynamically Dimensioned Search} (DDS) algorithm \cite{tolsonDynamicallyDimensionedSearch2007} into a surrogate-based optimization context, specifically leveraging RBF surrogates. DYCORS is particularly effective for high-dimensional optimization problems due to its dynamic search strategy, which balances exploration and exploitation. By progressively focusing the search and using a surrogate model to guide the selection of trial points, DYCORS can efficiently navigate the complex landscape of black-box functions.

\subsubsection{DYCORS - Surrogate Model}
In the context of DYCORS and other RBF algorithms discussed, the interpolation model employs RBFs to approximate the objective function based on known data points. Given \( n_d \) distinct points \( \mathbf{x}^{(1)}, \ldots, \mathbf{x}^{(n_d)} \in \mathbb{R}^{n_x} \) with corresponding function values \( f(\mathbf{x}^{(1)}), \ldots, f(\mathbf{x}^{(n_d)}) \), the RBF interpolant is formulated as:
\begin{equation}
s_n(\mathbf{x}) = \sum_{i=1}^{n_d} \lambda_i \phi(\|\mathbf{x} - \mathbf{x}^{(i)}\|) + p(\mathbf{x}), \quad \mathbf{x} \in \mathbb{R}^{n_x},
\end{equation}
where \( \phi \) is an RBF kernel, such as the cubic form \( \phi(r) = r^3 \). Other kernels like the thin plate spline and Gaussian can also be used. The coefficients \( \lambda_i \in \mathbb{R} \) and \( p(\mathbf{x}) \) represent a linear polynomial in \( n_x \) variables. 

The matrix \( \Phi \in \mathbb{R}^{n_d \times n_d} \) is defined by \( \Phi_{ij} = \phi(\|\mathbf{x}^{(i)} - \mathbf{x}^{(j)}\|) \), and \( P \in \mathbb{R}^{n_d \times (n_x+1)} \) is constructed such that its \( i \)-th row is \( [1, (\mathbf{x}^{(i)})^T] \). The cubic RBF model that interpolates the points \( (\mathbf{x}^{(1)}, f(\mathbf{x}^{(1)})), \ldots, (\mathbf{x}^{(n_d)}, f(\mathbf{x}^{(n_d)})) \) is obtained by solving the system:
\begin{equation}
\begin{bmatrix}
\Phi & P \\
P^T & 0
\end{bmatrix}
\begin{bmatrix}
\lambda \\
c
\end{bmatrix}
=
\begin{bmatrix}
F \\
0
\end{bmatrix},
\end{equation}
where \( F = (f(\mathbf{x}^{(1)}), \ldots, f(\mathbf{x}^{(n_d)}))^T \), \( \lambda = (\lambda_1, \ldots, \lambda_{n_d})^T \in \mathbb{R}^{n_d} \), and \( c = (c_1, \ldots, c_{n_x+1})^T \in \mathbb{R}^{n_x+1} \) are the coefficients for the RBF and linear polynomial \( p(\mathbf{x}) \), respectively. This coefficient matrix is invertible if and only if \( \text{rank}(P) = n_x + 1 \).

\subsubsection{DYCORS - Surrogate Optimization}

The optimization of the surrogate model within the DYCORS framework involves the \textit{Dynamic Coordinate Search} (DCS), which integrates several key steps to efficiently navigate \textit{High-dimensional, Expensive, and Black-box} (HEB) objective functions:
\\[8pt]
\textbf{Perturbation Probability:} A subset of the coordinates of the current best solution \( x_{\text{best}} \) is perturbed to generate trial points. The probability \( p_{\text{select}} \) of perturbing a coordinate is given by a strictly decreasing function \( \phi(n_x) \). This function ensures that the number of perturbed coordinates decreases as the algorithm progresses, leading to a more localized search.
\\[8pt]
\textbf{Generating Trial Points:} For each trial point, the selected coordinates are perturbed by adding normally distributed random variables with mean zero and a standard deviation \( \sigma_n \), known as the step size. The set of coordinates to be perturbed, \( I_{\text{perturb}} \), is chosen randomly in each iteration, ensuring a diverse exploration of the search directions.
\\[8pt]
\textbf{Selection of the Next Iterate:} From the set of generated trial points, the next iterate is selected using criteria that balance exploration and exploitation. In the DYCORS-LMSRBF algorithm, for instance, the selection is based on a weighted score combining the estimated function value from the RBF surrogate and the distance from previously evaluated points. This balance encourages both the discovery of new regions and the refinement of promising areas.
\\[8pt]
\textbf{Evaluation and Update:} The objective function \( f \) is evaluated at the selected trial point. This new data point is then incorporated into the surrogate model, updating it to reflect the most recent information. This iterative updating process ensures that the surrogate model continuously improves in accuracy.
\\[8pt]
\textbf{Adaptive Step Size Adjustment:} The step size \( \sigma_n \) is adjusted based on the success of the iterations. Parameters such as the number of consecutive successful iterations \( C_{\text{success}} \) and the number of consecutive failed iterations \( C_{\text{fail}} \) are monitored. Depending on these parameters and optional thresholds, the step size is either increased or decreased, facilitating an adaptive search strategy.
\\[8pt]
\textbf{Convergence Criteria:} The optimization process continues until a predefined convergence criterion is met. Common criteria include a maximum number of function evaluations, a maximum number of iterations, or a tolerance threshold for changes in the objective function value. The best solution found during the search is then reported as the optimal solution.
\\[8pt]
\textbf{Specific Implementations in DYCORS:} The DYCORS algorithm has two specific variants: DYCORS-LMSRBF and DYCORS-DDSRBF. In the performance assessment carried out in this chapter, the DYCORS-LMSRBF variant is used. This variant has been shown to effectively balance exploration and exploitation by selecting the iterate based on a combination of the RBF surrogate value and the distance from previously evaluated points. The nuances of the variants are briefly detailed below:

\begin{itemize}
    \item \textbf{DYCORS-LMSRBF:} This variant is a modification of the LMSRBF algorithm. It uses the dynamic coordinate search strategy and selects the iterate based on a weighted score of the RBF surrogate value and the distance criterion. The perturbation probability function \( \phi(n_x) \) is chosen to ensure an initial high probability of perturbing coordinates, which decreases logarithmically as the number of iterations increases.
    \item \textbf{DYCORS-DDSRBF:} This variant incorporates the RBF surrogate into the DDS algorithm. It maintains the dynamic coordinate search strategy but uses the surrogate model to enhance the efficiency and accuracy of the search.
\end{itemize}

\newpage 
\section{Surrogate Optimization Methods Performance Assessment}\label{sec:b}
In addition to the theoretical presentation, a performance assessment was conducted to evaluate the different model-based algorithms highlighted in this chapter. This assessment is supplemented with the associated code, such that the interested reader can further explore this performance assessment and the algorithms studied. 

\subsection{Performance Assessment Procedure}\label{sec:bp}
For this comparative assessment, the algorithms to be assessed along with the test functions have been defined. The set containing the different algorithms is called $\mathbb{A}$, and the set containing the different test functions as $\mathbb{F}$. This book chapter conducts a performance assessment for unconstrained black-box problems as well as for constrained back-box problems. In the unconstrained case the functions and algorithms are as follows:
\begin{align*}
    a &\in \mathbb{A}, & \mathbb{A} &= \left\{ \text{LSQM, SNOBFIT, SRBF, DYCORS, SOP, COBYLA, COBYQA, CUATRO, BO, ENTMOOT} \right\} \\
    f &\in \mathbb{F}, & \mathbb{F} &= \left\{ \text{Ackley, Levy, Rosenbrock, Quadratic} \right\}.
\end{align*}

Therefore every algorithm $a \in \mathbb{A}$ is assessed on every objective function $f \in \mathbb{F}$ and its performance is compared relative to the other optimization algorithms on the same function. The procedure allocates each algorithm a budget of $n_{e}$ function evaluations, and the optimization of $f$ is conducted five times per algorithm to account for algorithm and function evaluation stochastic factors. For a given algorithm, the trajectories of function values $y_k$ with $k=1, 2, ...n_{e}$ are stored in five vectors each of length $n_{e}$, where only the best-so-far values within a trajectory are stored. The evaluation budget is given proportional to the dimension of the function; where $n_x$ denotes dimensions, the function evaluation budgets are allocated as follows, 20 function evaluations for $n_x=2$, 50 function evaluations for $n_x=5$, and 100 function evaluations for $n_x=10$.
\\[8pt]
In addition to displaying all the performance figures for each algorithm, a quantitative metric is provided for each algorithm's performance for each function. It is important to note, that this performance assessment is not based on the final objective value that the algorithms arrive at. Instead, the assessment is based on the algorithms' respective trajectories; trajectories offer a more robust, and less arbitrary, measure for comparison. Equation \ref{eq:norm_score}, defines the normalized scoring metric. 
\\[8pt]
To account for random initialization resulting in sub-optimal objective function values, as well as the necessity of using initial sampling to build the first surrogate model before starting the optimization (see Figure \ref{fig:Image_opti_flow}) the best point found so far does not start counting from the very first evaluation, but rather after the fist $n_c$ evaluations. 
For $n_x=2$ this is $n_c=5$, for $n_x=5$ this is $n_c=10$ and for $n_x=10$ this is $n_c=15$. By this, the effective length of the trajectory considered for the benchmarking reduces to $n$, yielding $n=15$, $n=40$, and $n=85$ for $d=2$, $d=5$, and $d=10$ respectively. 

For example, when LSQM is benchmarked on the 2-dimensional Levy function LSQM will do a total of 20 function evaluations. Initially, $n_x+1\rightarrow3$ evaluations are done to construct the surrogate (implementation details like these can be found in the GitHub repository). After LSQM builds the initial surrogate it will perform 2 more function evaluations for a total of 5, which corresponds to $n_c=5$, before our benchmarking procedure starts counting. From iterations 6 to 20, the best point found so far will be compared to the other algorithm's performance (noting that averages over a total of 5 runs are what is compared).  

The specific scoring metric used to deliver a quantitative score is the trajectories of length $n$, and they are compared using the quotient $r_{k,a}^{(n_x)}$, which represents the relative performance of algorithm $a$ on function $f$ for iteration $k$ with $n_x$ input dimensions:
\begin{equation}\label{eq:norm_score}
    r_{k,a}^{(n_x)} = \frac{\Tilde{y}_{k}-y^{mean}_{k,a}}{\Tilde{y}_{k}-y^{*}_{k}},  \quad \quad \quad \quad 0 \le r_{k,a}^{(n_x)} \le 1 
\end{equation}
Note, for simplicity, the superscript $(n_x)$ is only shown on the left-hand side, whereas all variables presented account for the respective dimension.
Here, $y^{mean}_{k,a}$ denotes the mean value of the objective function $f$ along the trajectory achieved by algorithm $a$ given $n_x$ input dimensions. This mean is computed at a specific trajectory position $k$ and is derived from the averaging of five optimization runs. 

$\Tilde{y}_{k}$ represents the maximum (worst) function value for iteration $k$ out of all algorithms investigated. In contrast, $y^{*}_{k}$ is the lowest (best) function value achieved and therefore denotes the best performance by any algorithm on iteration $k$. To summarise, the closer $y^{mean}_{k,a}$ is to $y^{*}_{k}$ - and hence the closer $r_{k,a}^{(d)}$ is to 1 - the better $a$ performed in iteration $k$ compared to other algorithms. 

% \\[8pt]
Finally, $p_{a}^{(n_x)}$ as in Equation \ref{eq:pa} represents the overall performance of algorithm $a$ on objective function $f$ with a given budget $n$ for $d$ input dimensions. 
Again, the closer $p_a^{(n_x)}$ is to 1, the better $a$ performed relative to the algorithms in $\mathbb{A}$.
\begin{equation}\label{eq:pa}
    p_{a}^{(n_x)} = %\frac{\| \bm{t}^{rel}_{a} \|_1}{dim(\bm{t}^{rel}_{a})} =
    \frac{\sum_{k=1}^{n} r_{k,a}^{(n_x)}}{n}, \quad \quad \quad \quad 0 \le p_{a}^{(n_x)} \le 1
\end{equation}

Geometrically, this score represents how low (good) an algorithm's trajectory is, normalized by the best and worst performance amongst the algorithms.

\subsection{Unconstrained Performance Assessments}\label{sec:ub}
This section presents the performance assessment of the model-based algorithms presented in Section \ref{sec:mb dfo} for unconstrained optimization problems. The section is structured as follows before a conclusion is drawn for the unconstrained case:
\begin{itemize}
    \item \textbf{Mathematical Objective Functions:} Objective functions for the respective optimization problems are presented.
    \item \textbf{2D Trajectory Plots:} For qualitative evaluation, 2-dimensional plots are presented that help gain intuition on the algorithm's behavior.
    \item \textbf{1D Convergence Plots:} For quantitative analysis, 1-dimensional plots to observe the algorithm's convergence.
    \item \textbf{Performance Tables:} Tables summarizing the performance in benchmarking numbers following the procedure described in Section \ref{sec:bp}.   
\end{itemize}

\subsubsection{Mathematical Objective Functions}\label{sec:umof}
The objective functions used in the unconstrained benchmarking are the Ackley function as proposed by David H. Ackley in 1987 \cite{ackleyConnectionistMachineGenetic1987}, the Levy function by Paul Levy in 1954 \cite{levyMouvementBrownien1954}, the Rosenbrock function by Howard H. Rosenbrock in 1960 \cite{rosenbrockAutomaticMethodFinding1960}, and a quadratic mildly ill-conditioned function. The former two are multi-modal functions meaning they have multiple local minima, whereas the latter two are uni-modal functions, denoting their single global minimum. Mathematical expressions can be found in Equations \ref{eq:Ackley}, \ref{eq:Levy}, \ref{eq:Rosenbrock}, and \ref{eq:Antonio}:
\begin{itemize}
    \item Ackley: 
    \begin{equation} \label{eq:Ackley}
        f(\mathbf{x}) = -a \exp \left( -b \sqrt{\frac{1}{n_x} \sum_{i=1}^{n_x} x_i^2} \right) - \exp \left( \frac{1}{n_x} \sum_{i=1}^{n_x} \cos(c x_i) \right) + a + \exp(1)
    \end{equation}
    \item Levy: 
    \begin{equation} \label{eq:Levy}
        f(\mathbf{x}) = \sin^2(\pi w_1) + \sum_{i=1}^{n_x-1} \left[ (w_i - 1)^2 (1 + 10 \sin^2(\pi w_i+1)) \right] + (w_{n_x} - 1)^2 (1 + \sin^2(2 \pi w_{n_x}))
    \end{equation}
    \item Rosenbrock: 
    \begin{equation} \label{eq:Rosenbrock}
        f(\mathbf{x}) = \sum_{i=1}^{n_x-1} \left[ 100(x_{i+1} - x_i^2)^2 + (1 - x_i)^2 \right]
    \end{equation}
    \item Quadratic: 
    \begin{equation} \label{eq:Antonio}
        f(x) = \sum_{i=1}^{n_x} \left[ (i \cdot x_i)^2 + \left(\frac{ai}{n_x}\right) x_i x_{n_x} \right]
    \end{equation}
\end{itemize}
For all functions $\mathbf{x} = (x_1, x_2, \ldots, x_{n_x})$ is the $n_x$-dimensional input vector. For the Ackley function $a$, $b$, and $c$ are constants set to $a=20$, $b=0.2$, and $c=2\pi$ respectively. For the Levy function $w_i = 1 + \frac{x_i - 1}{4}$ for $i = 1, 2, \ldots, n_x$. In the Quadratic function $a = 1.9$, and $x = (x_1, x_2, \ldots, x_{n_x})$. The respective case studies are summarized in Table \ref{tab:Uncon}.
N.B., that all case studies in this section are deterministic, a study of noise has been conducted in: \cite{vandebergDatadrivenOptimizationProcess2022}. Furthermore, both chemical engineering case studies involved an element of stochasticity.
\begin{table}[h]
    \centering
    \captionof{table}{Case studies for the unconstrained benchmarking. The input dimensions are $n_x = 2,5,7$ on domain $\mathbf{x} \in [-5, 5]^{n_x}$ with global minimum $f(\mathbf{x}^*) = 0$} \label{tab:Uncon} 
    \begin{tabular}{lll}
    \hline
         \textbf {Case Study} &\textbf{Topology} & \textbf{Global Minimum}\\
    \hline
         Ackley Function & non-convex & $\mathbf{x}^* = (0,\dots, 0)$ \\
         Levy Function& non-convex & $\mathbf{x}^* = (1,\dots, 1)$\\
         Rosenbrock Function& non-convex & $\mathbf{x}^* = (1,\dots, 1)$ \\
         Quadratic Function & convex, ill-conditioned & $\mathbf{x}^* = (0,\dots, 0)$ \\
    \hline
    \end{tabular}
\end{table}

\newpage
\subsubsection{Results - Convergence Plots}
%\ifincludegraphs
\begin{figure}[H]
    \begin{subfigure}{0.40\linewidth}
        \centering
            \includegraphics[width=\linewidth]{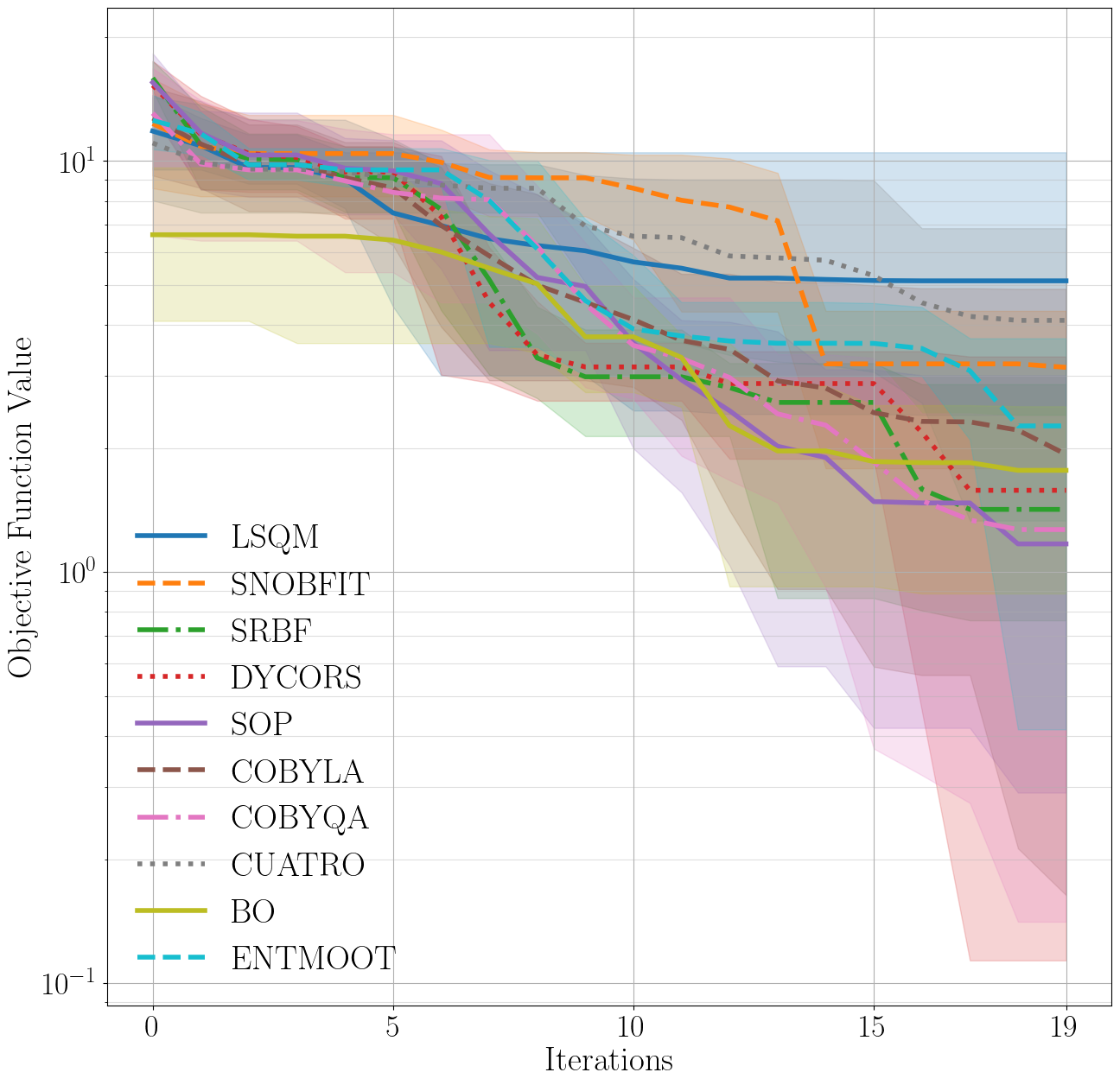}
            \caption{Ackley 2D}
    \end{subfigure}
    \quad
    \quad
    \quad
    \quad
    \quad
    \quad
    \quad
    \begin{subfigure}{0.40\linewidth}
        \centering
            \includegraphics[width=\linewidth]{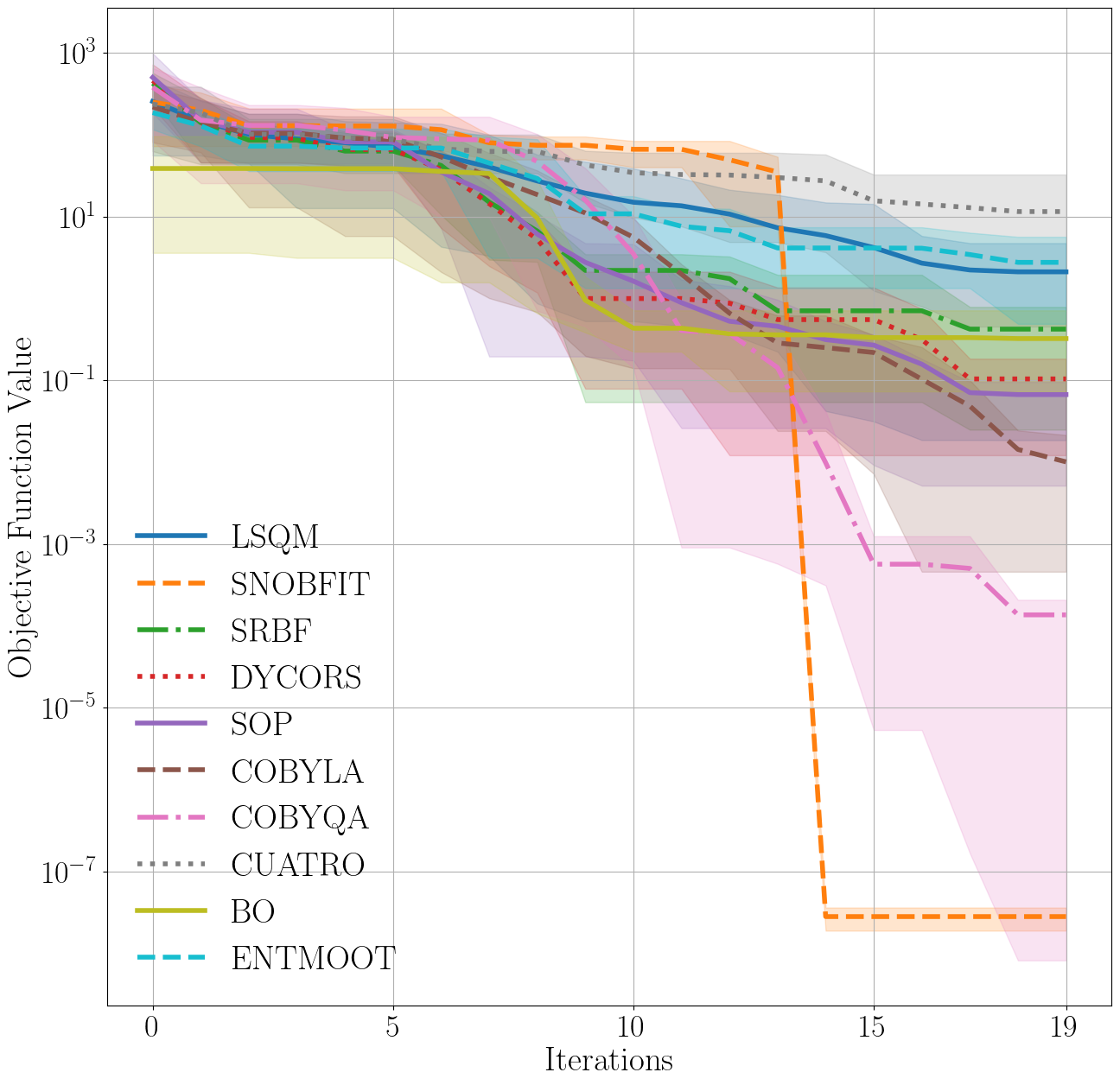}
            \caption{Quadratic 2D}            
    \end{subfigure}        
    \\
    \begin{subfigure}{0.40\linewidth}
        \centering
            \includegraphics[width=\linewidth]{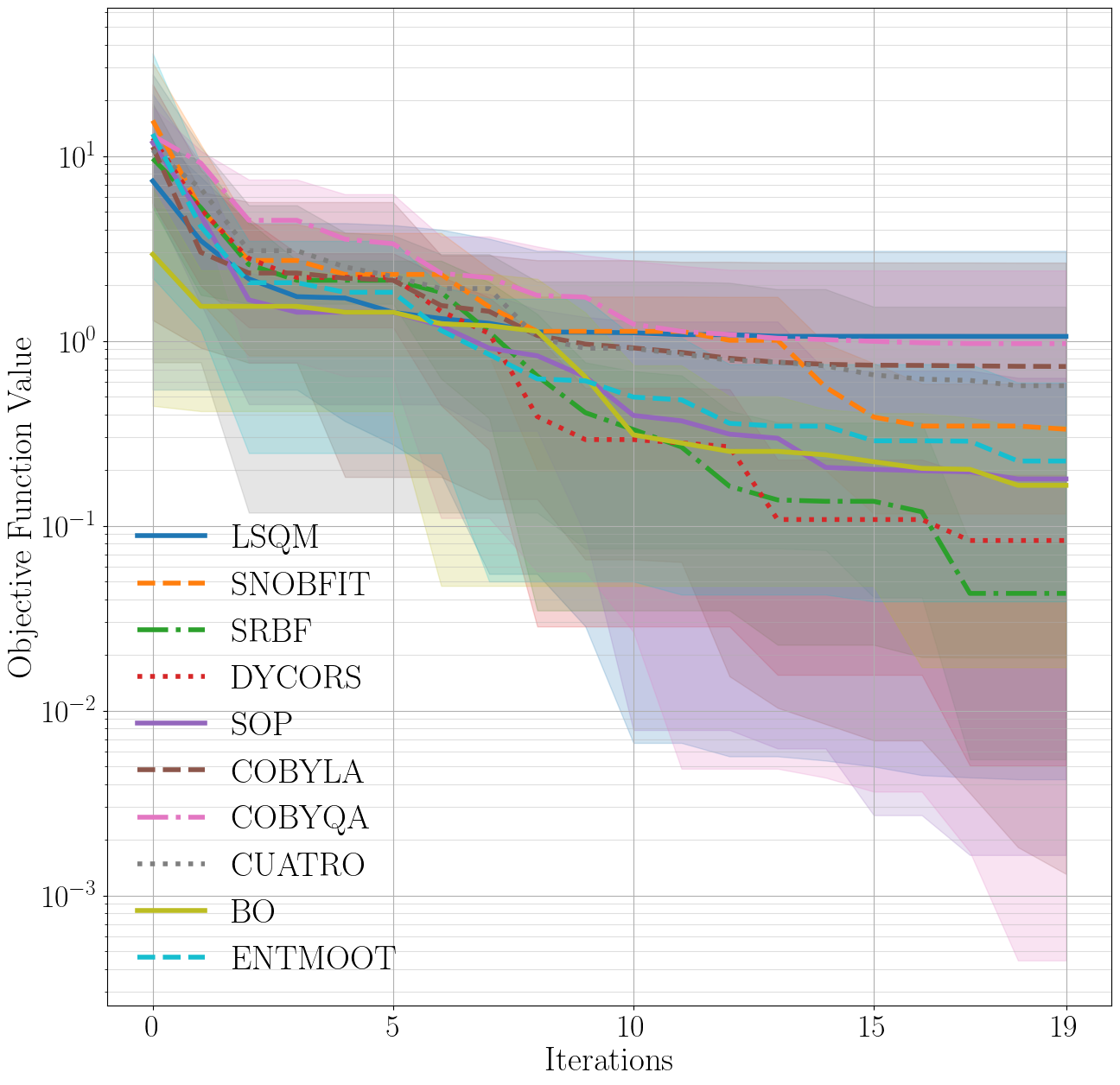}
            \caption{Levy 2D}
    \end{subfigure}   
    \quad
    \quad
    \quad
    \quad
    \quad
    \quad
    \quad
    \begin{subfigure}{0.40\linewidth}
        \centering
            \includegraphics[width=\linewidth]{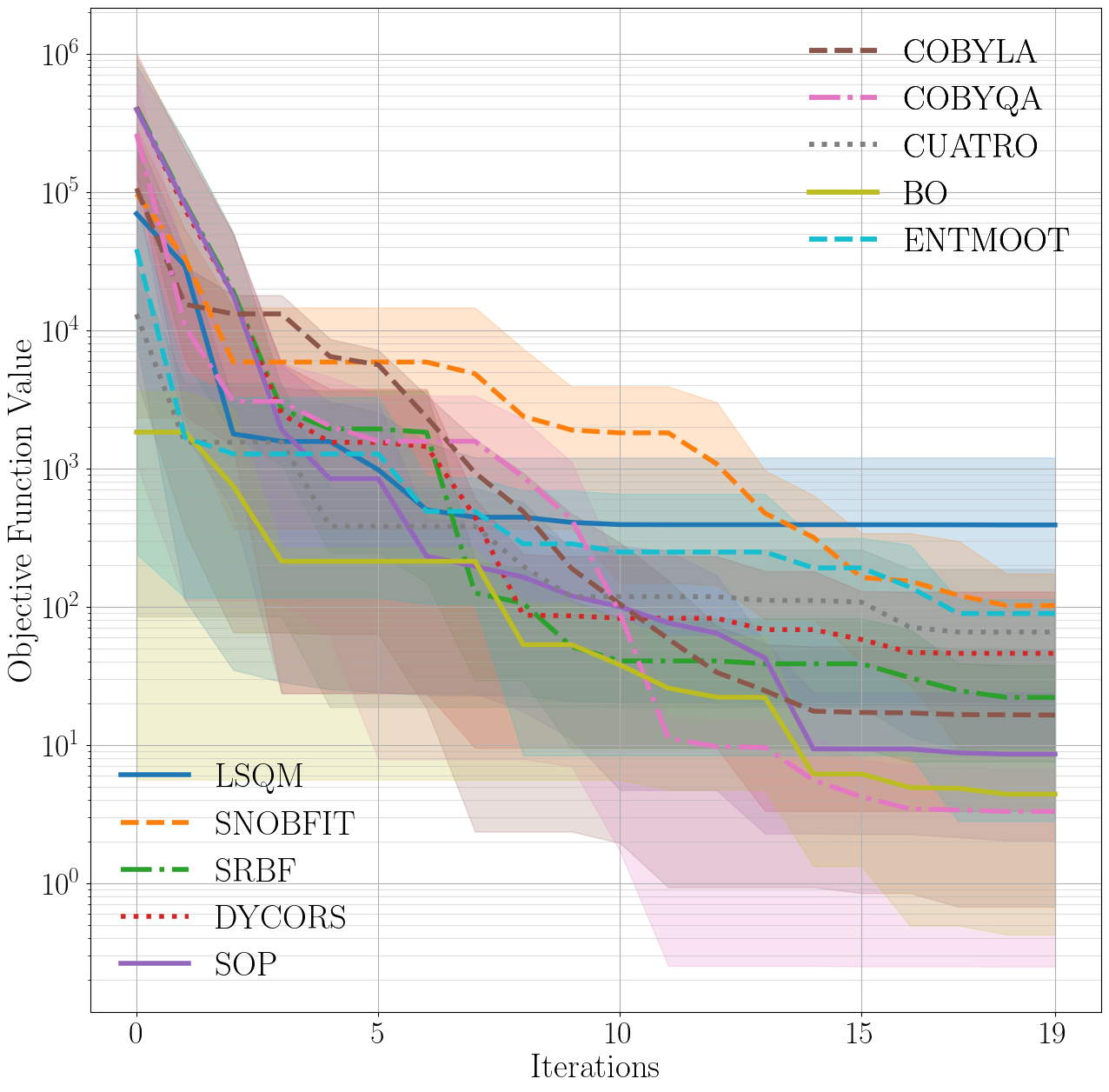}
            \caption{Rosenbrock 2D}
    \end{subfigure}
    \\
    \begin{subfigure}{0.40\linewidth}
        \centering
            \includegraphics[width=\linewidth]{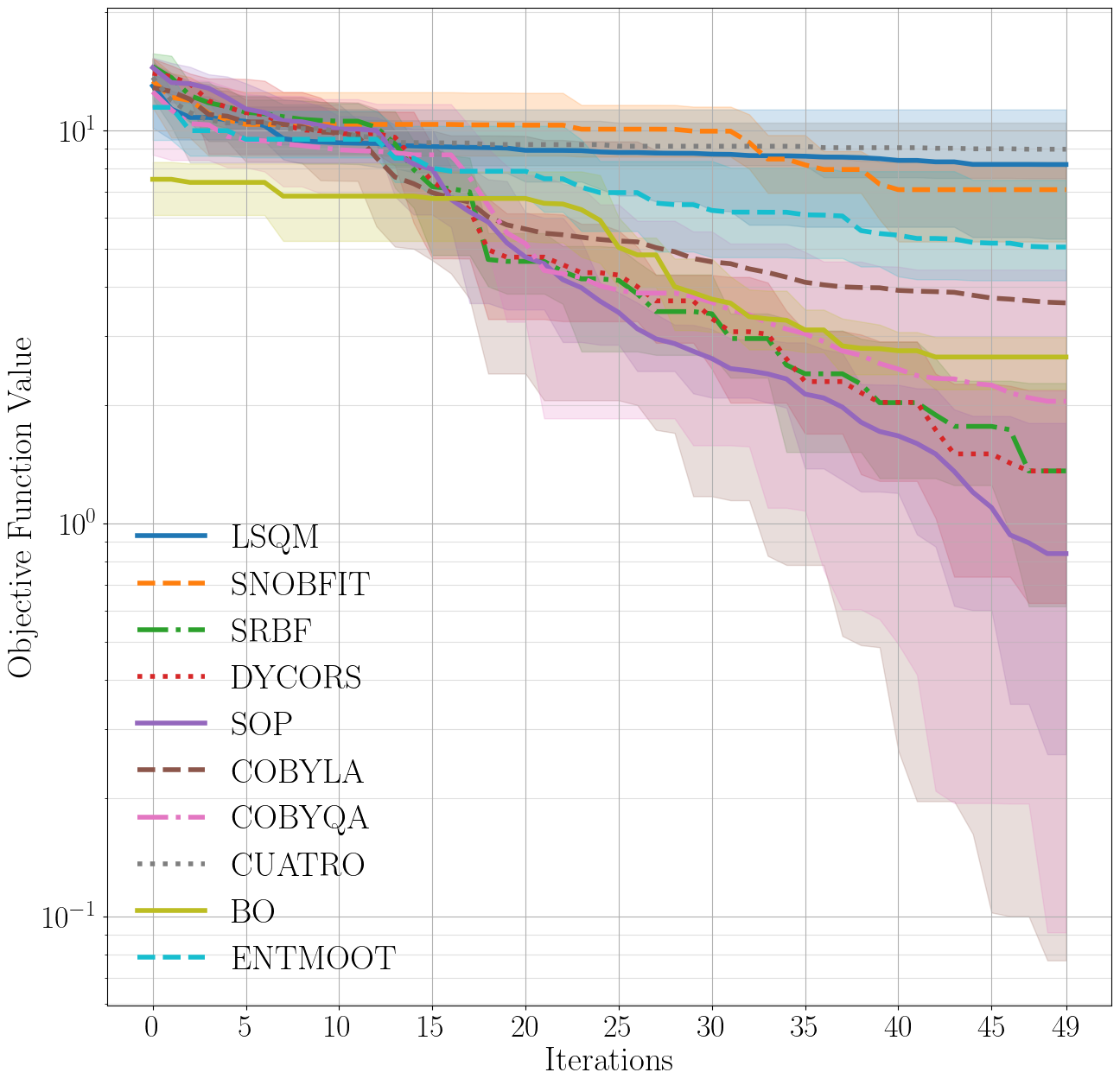}
            \caption{Ackley 5D}
    \end{subfigure}
    \quad
    \quad
    \quad
    \quad
    \quad
    \quad
    \quad
    \begin{subfigure}{0.40\linewidth}
        \centering
            \includegraphics[width=\linewidth]{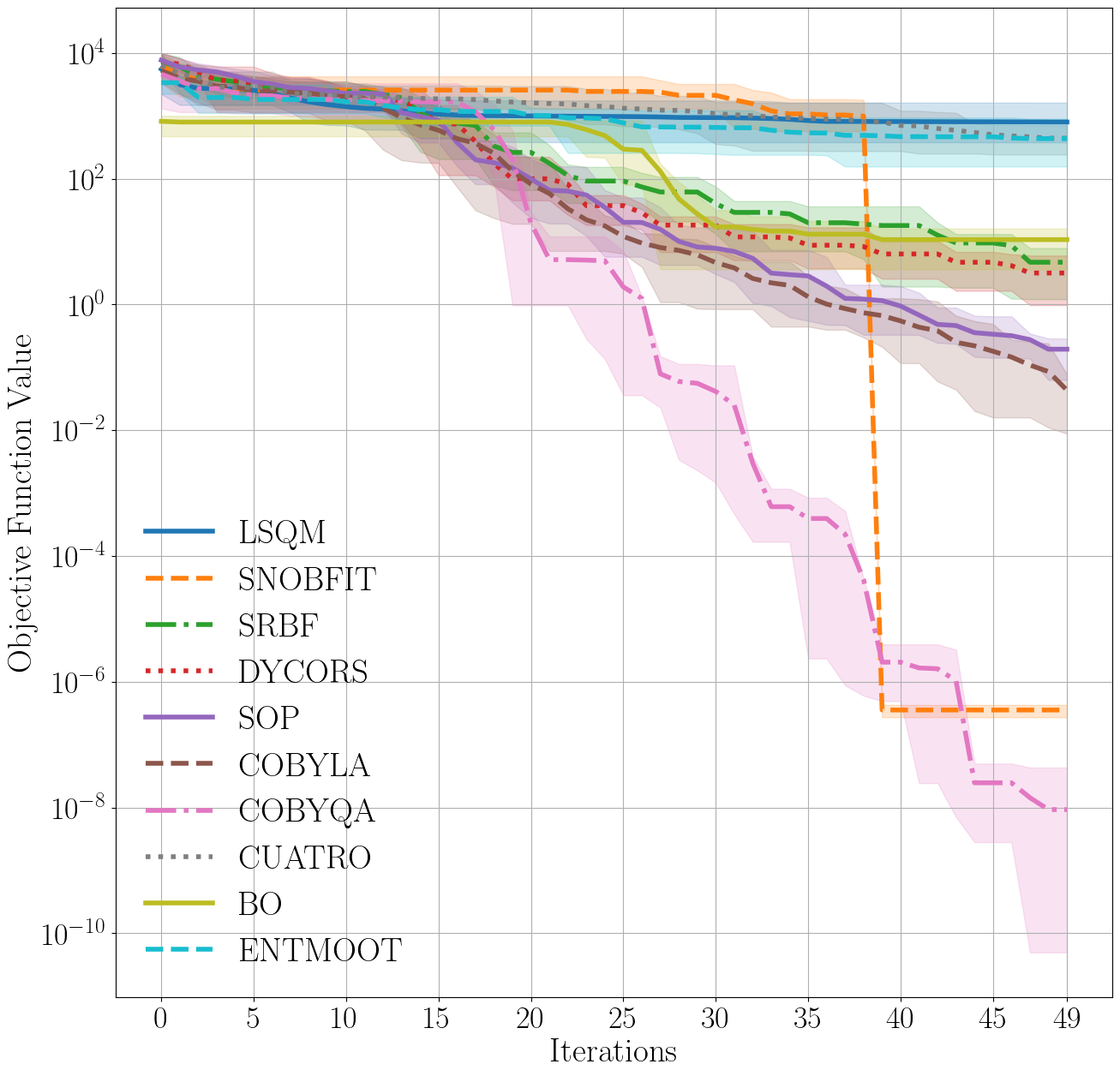}
            \caption{Quadratic 5D}
    \end{subfigure}
    \caption{Convergence plots, showing mean objective function values and $10\%-90\%$-intervals enveloping trajectories over 10 repetitions from 10 different starting points. Sub-captions indicate test function and input dimensionality.}    
\end{figure}

\clearpage
%\ifincludegraphs
\begin{figure}[H]\ContinuedFloat
    \begin{subfigure}{0.40\linewidth}
        \centering
            \includegraphics[width=\linewidth]{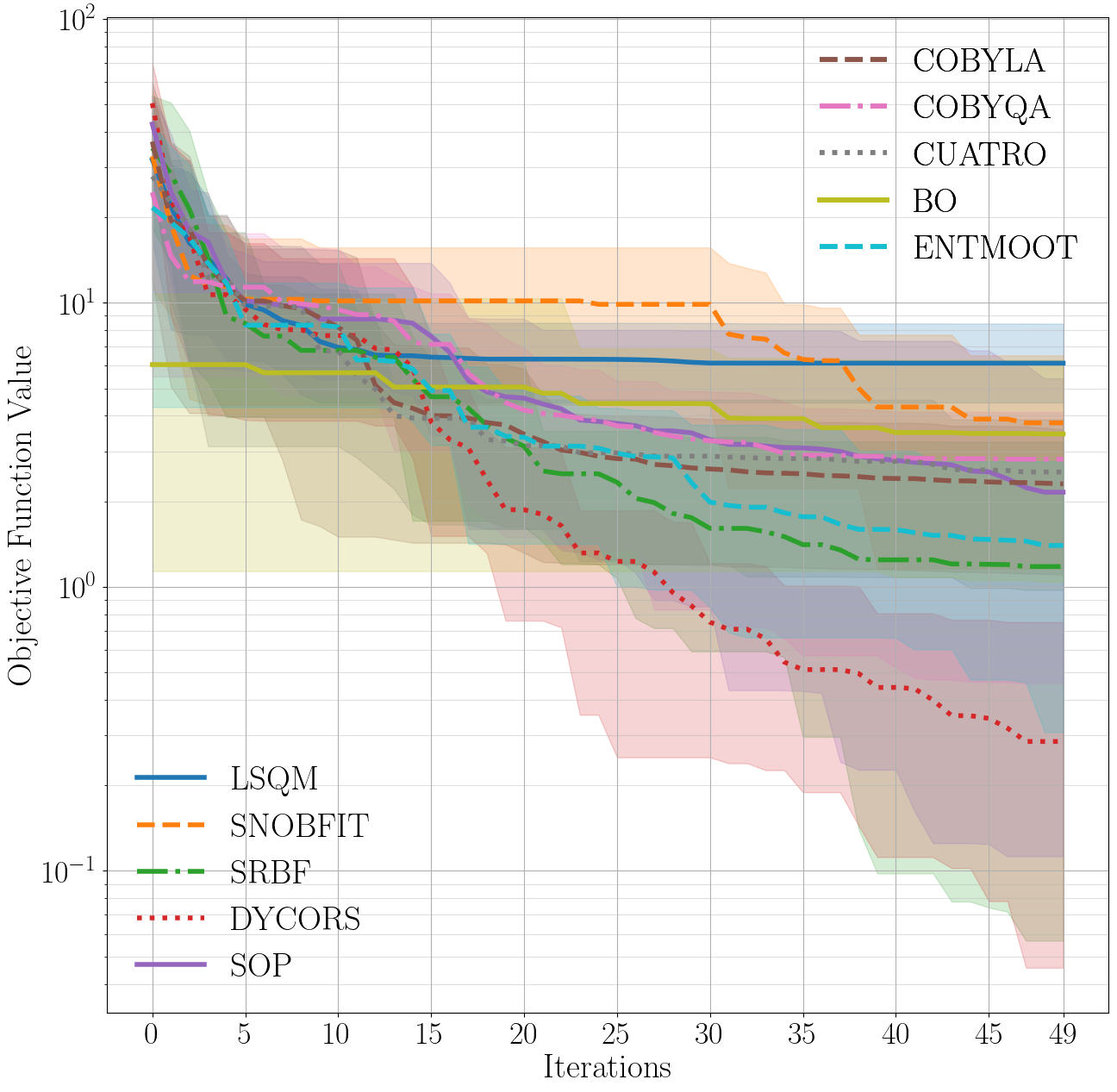}
            \caption{Levy 5D}
    \end{subfigure}
    \quad
    \quad
    \quad
    \quad
    \quad
    \quad
    \quad
    \begin{subfigure}{0.40\linewidth}
        \centering
            \includegraphics[width=\linewidth]{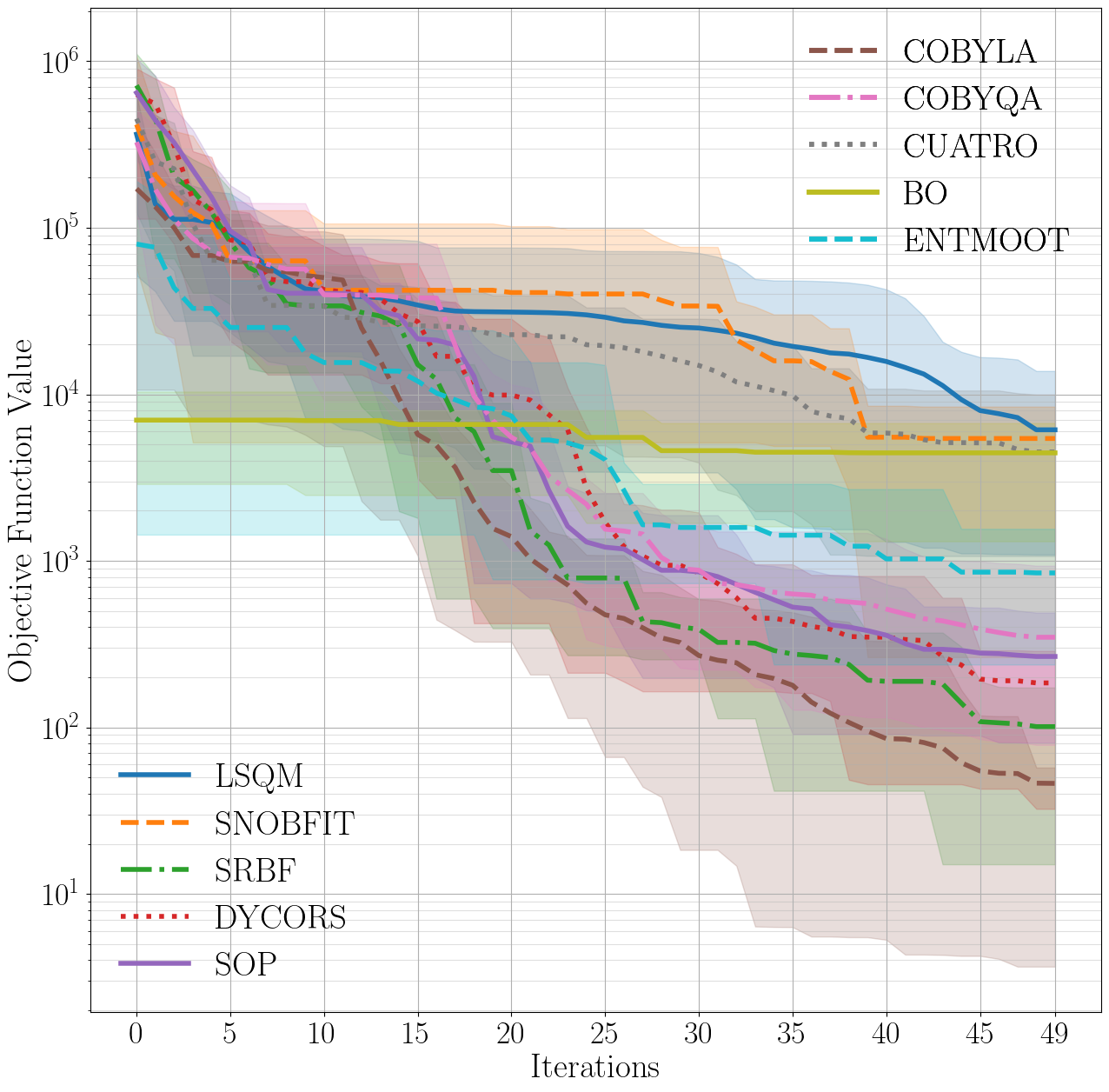}
            \caption{Rosenbrock 5D}            
    \end{subfigure}        
    \\
    \begin{subfigure}{0.40\linewidth}
        \centering
            \includegraphics[width=\linewidth]{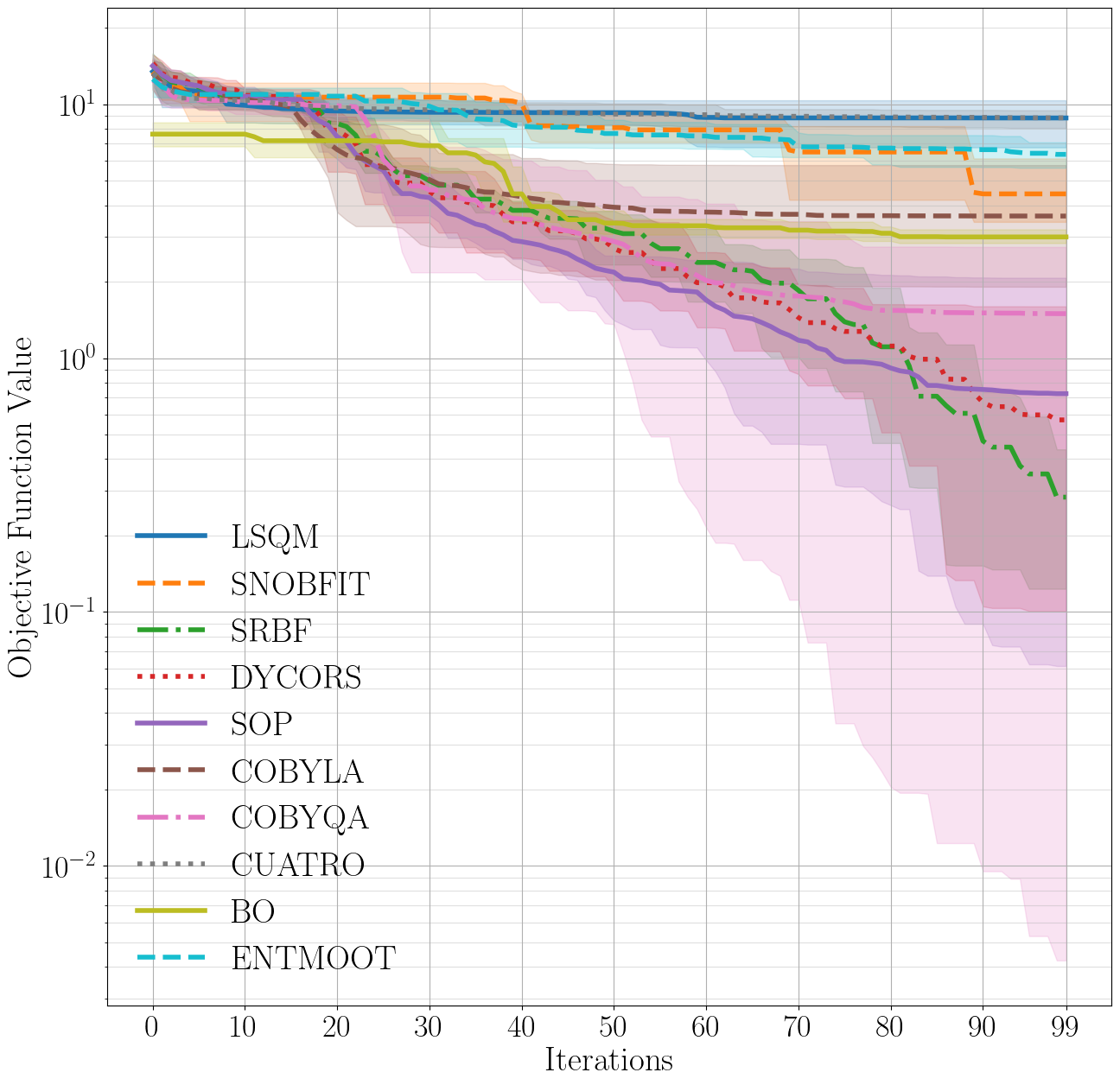}
            \caption{Ackley 7D}
    \end{subfigure}   
    \quad
    \quad
    \quad
    \quad
    \quad
    \quad
    \quad
    \begin{subfigure}{0.40\linewidth}
        \centering
            \includegraphics[width=\linewidth]{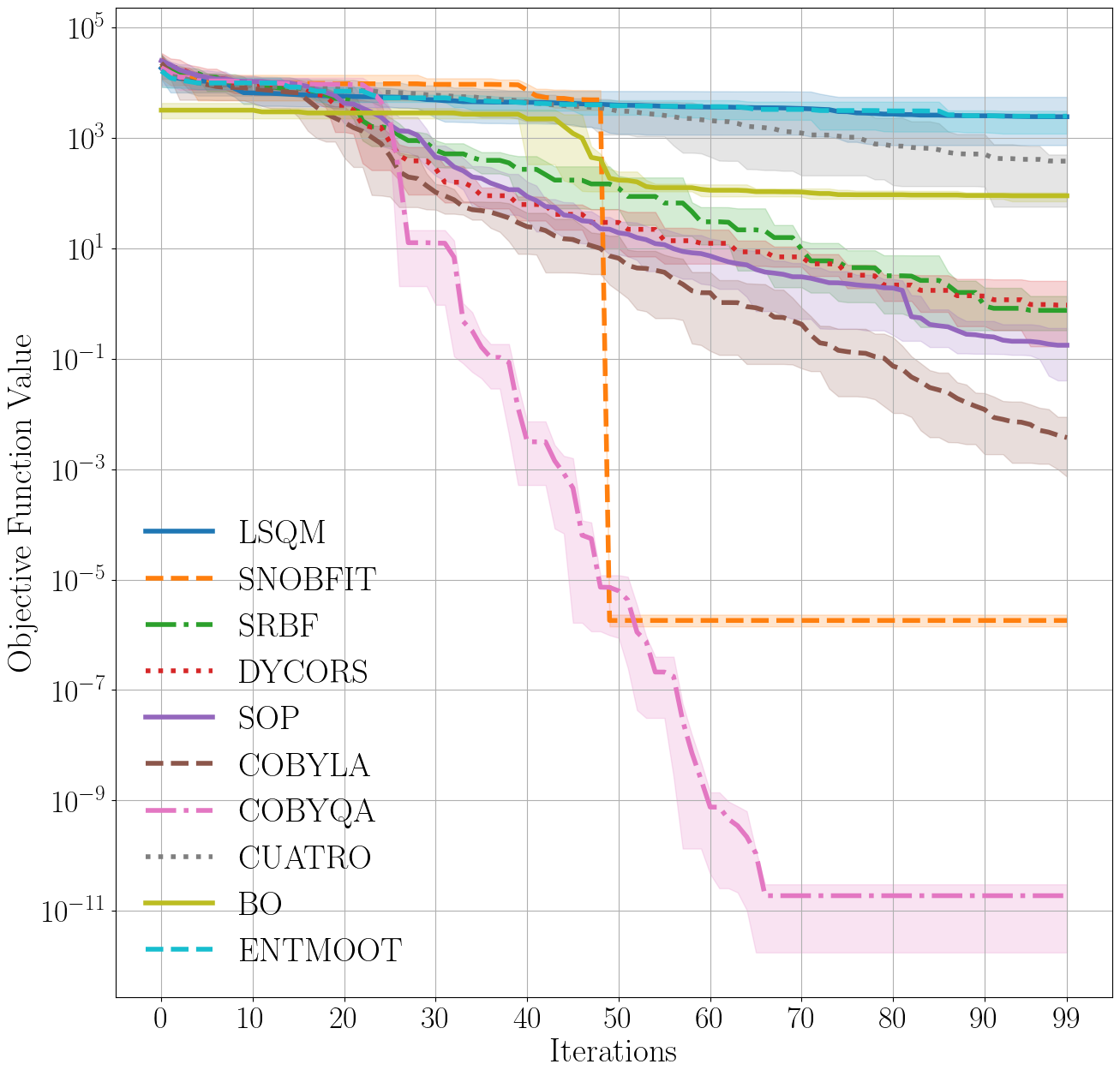}
            \caption{Quadratic 7D}
    \end{subfigure}
    \\
    \begin{subfigure}{0.40\linewidth}
        \centering
            \includegraphics[width=\linewidth]{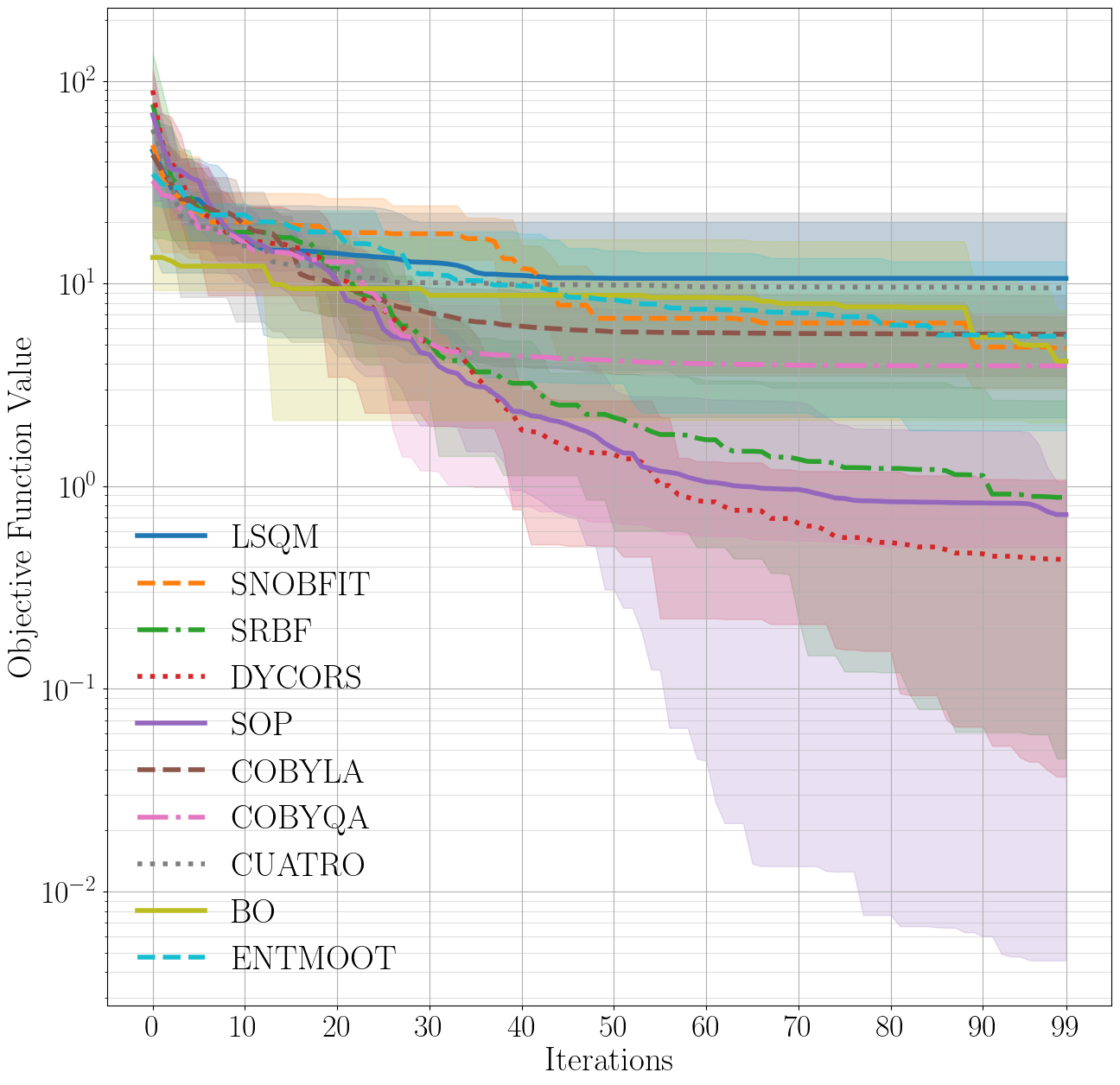}
            \caption{Levy 7D}
    \end{subfigure}
    \quad
    \quad
    \quad
    \quad
    \quad
    \quad
    \quad
    \begin{subfigure}{0.40\linewidth}
        \centering
            \includegraphics[width=\linewidth]{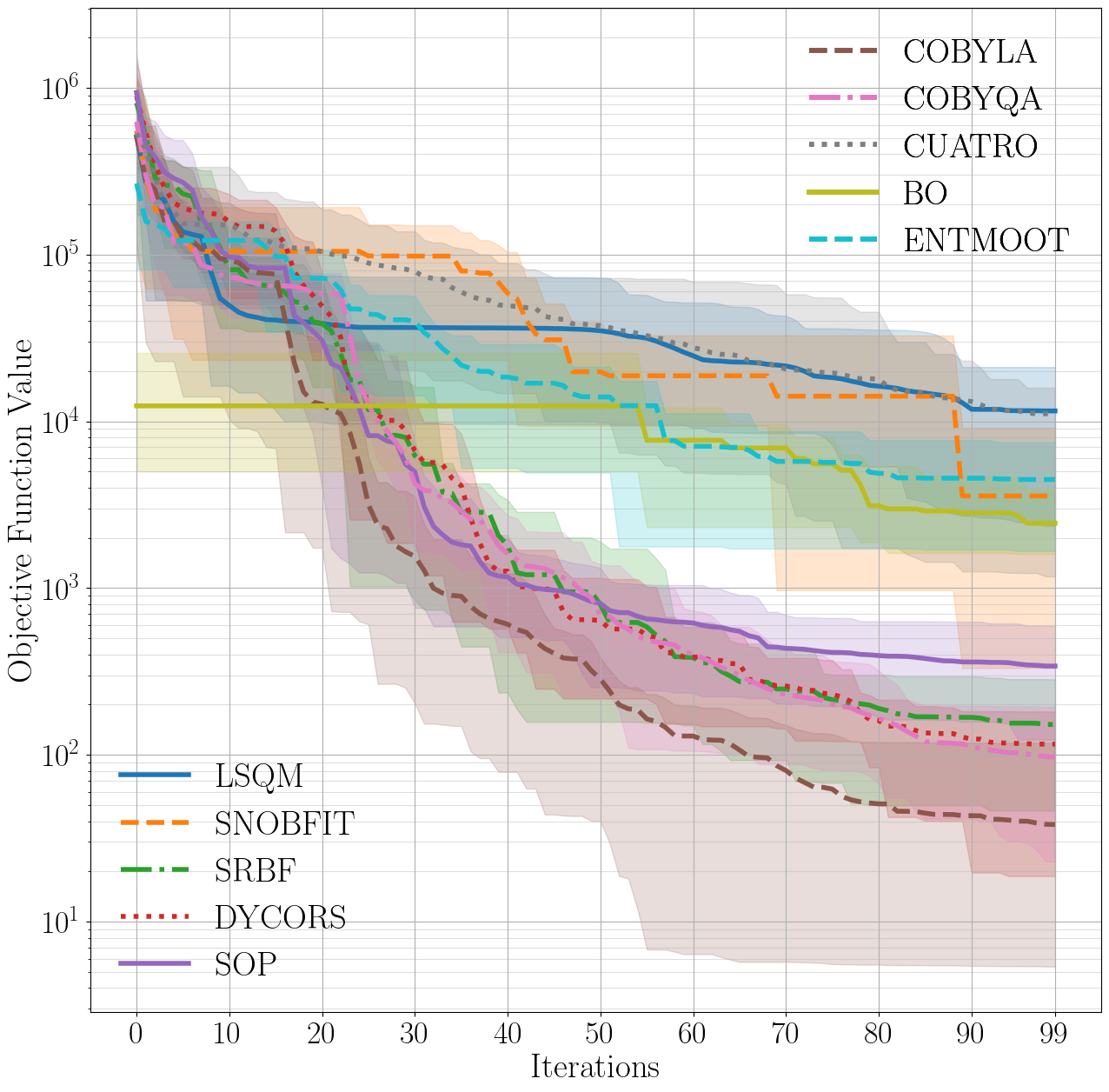}
            \caption{Rosenbrock 7D}
    \end{subfigure}
    \caption{Convergence plots, showing mean objective function values and $10\%-90\%$-intervals enveloping trajectories over 10 repetitions from 10 different starting points. Sub-captions indicate test function and input dimensionality cont.}    
\end{figure}
%\fi
%%%%%%%%%%%%%%%%%%%%%%%%%%%%%%%%%%%%%%%%%%%%%%%%%%%%%%%%%%%%%%%%%%%%%%%%%%%%%%%%%%%%%%%%%%%%%%%%%%%%%%%%%%%%%%%%%%%%%%%%%%%%%%%%%%%%%%%%%%%%%%%

\subsubsection{Results - Trajectory Plots}
%\ifincludegraphs
\begin{figure}[H]
    \begin{subfigure}{0.385\linewidth}
        \centering
            \includegraphics[width=\linewidth]{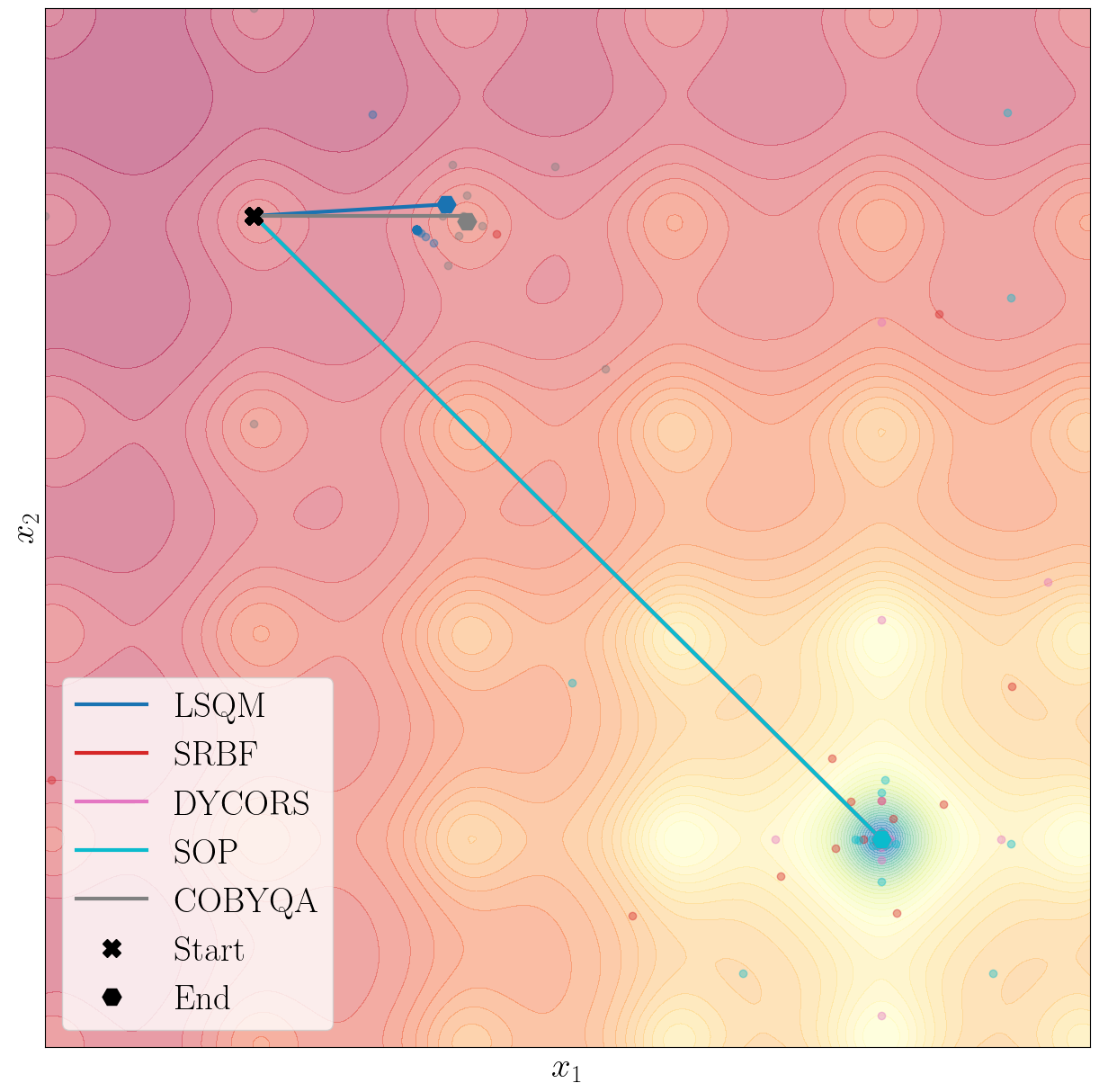}
            \caption{Ackley Contour and Trajectories}
    \end{subfigure}
    \quad
    \quad
    \quad
    \quad
    \quad
    \quad
    \quad
    \quad
    \begin{subfigure}{0.385\linewidth}
        \centering
            \includegraphics[width=\linewidth]{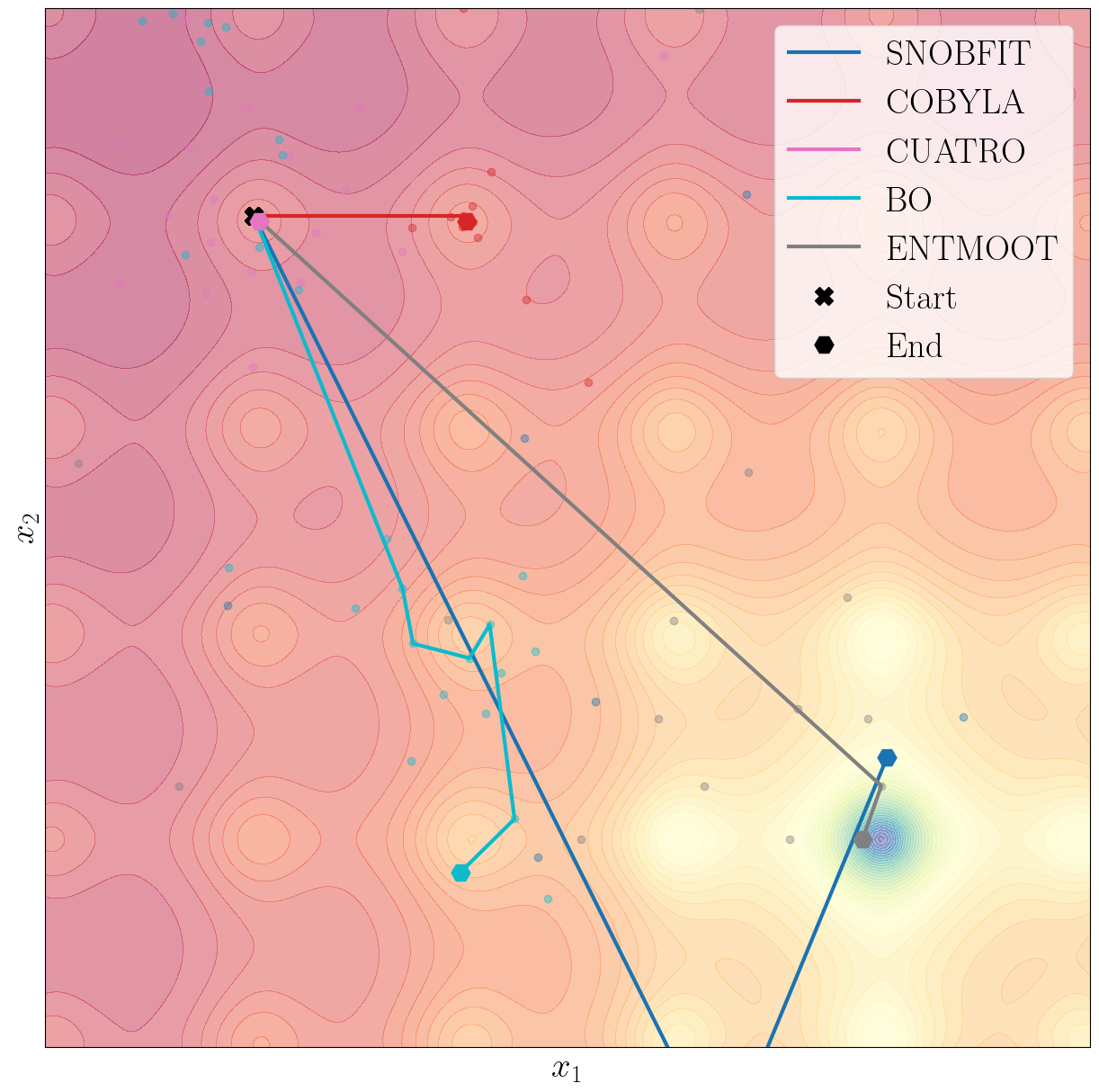}
            \caption{Ackley Contour and Trajectories}            
    \end{subfigure}        
    \\
    \begin{subfigure}{0.385\linewidth}
        \centering
            \includegraphics[width=\linewidth]{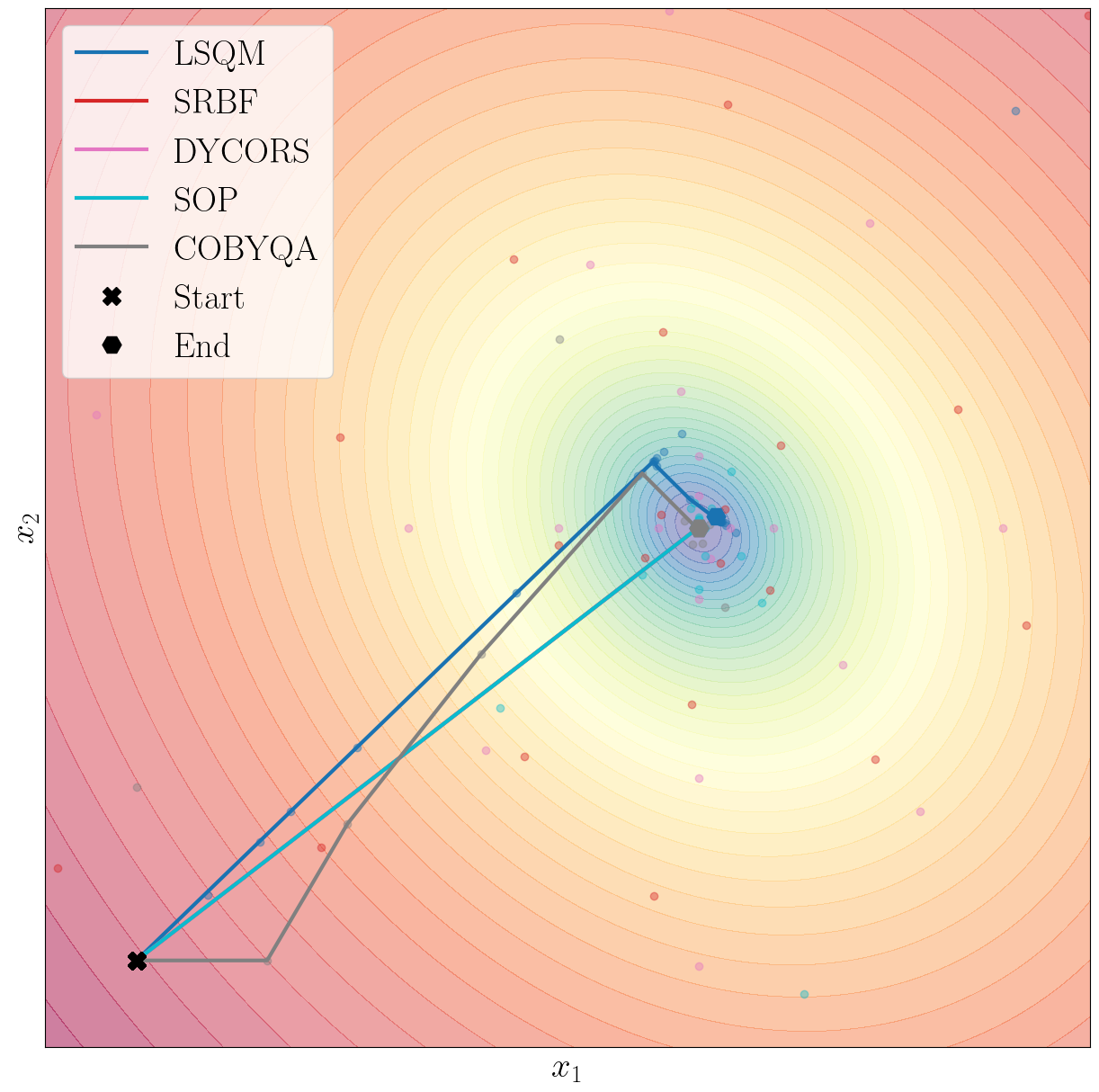}
            \caption{Quadratic Contour and Trajectories}
    \end{subfigure}   
    \quad
    \quad
    \quad
    \quad
    \quad
    \quad
    \quad
    \quad
    \begin{subfigure}{0.385\linewidth}
        \centering
            \includegraphics[width=\linewidth]{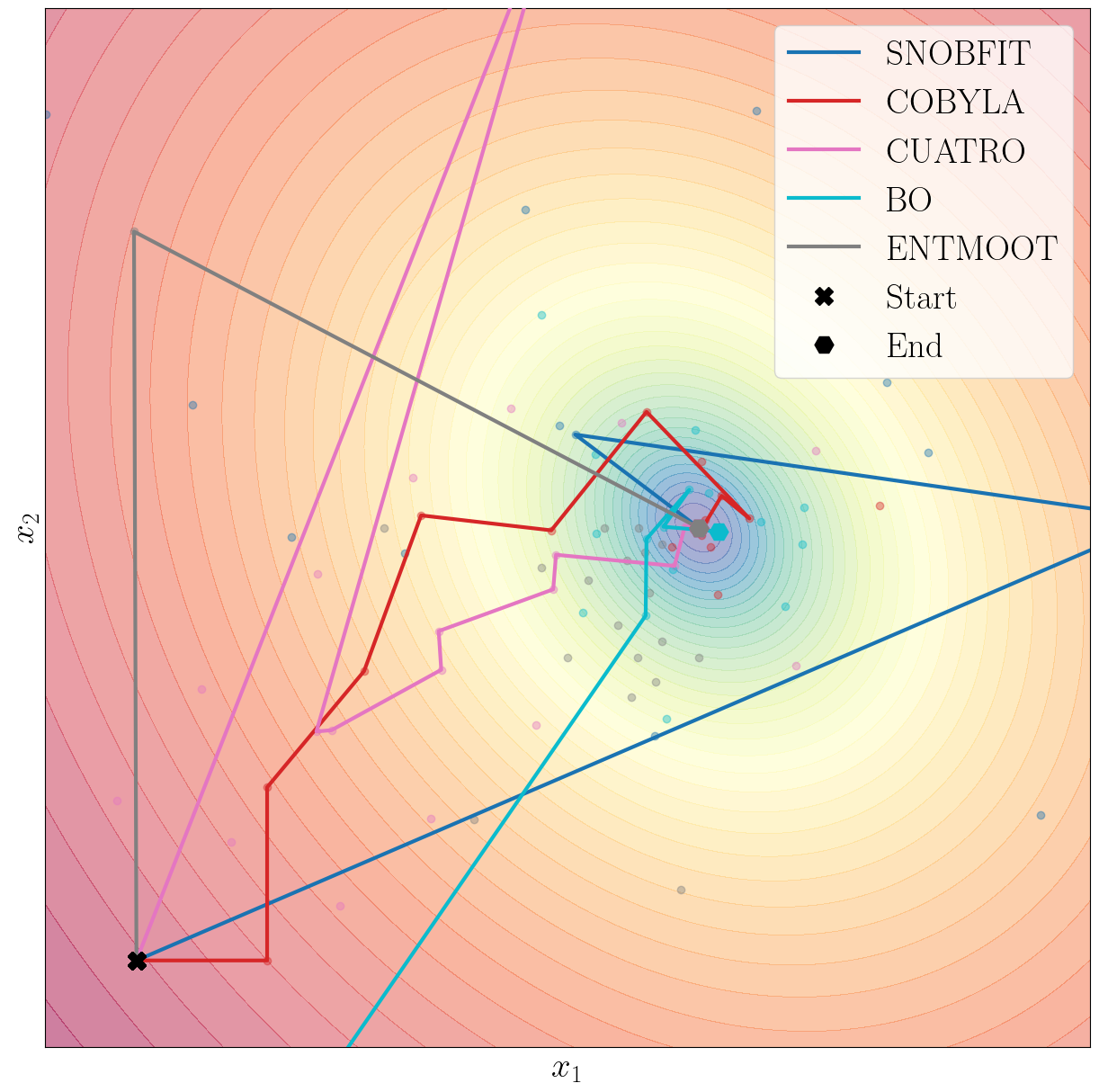}
            \caption{Quadratic Contour and Trajectories}
    \end{subfigure}
    \\
    \begin{subfigure}{0.385\linewidth}
        \centering
            \includegraphics[width=\linewidth]{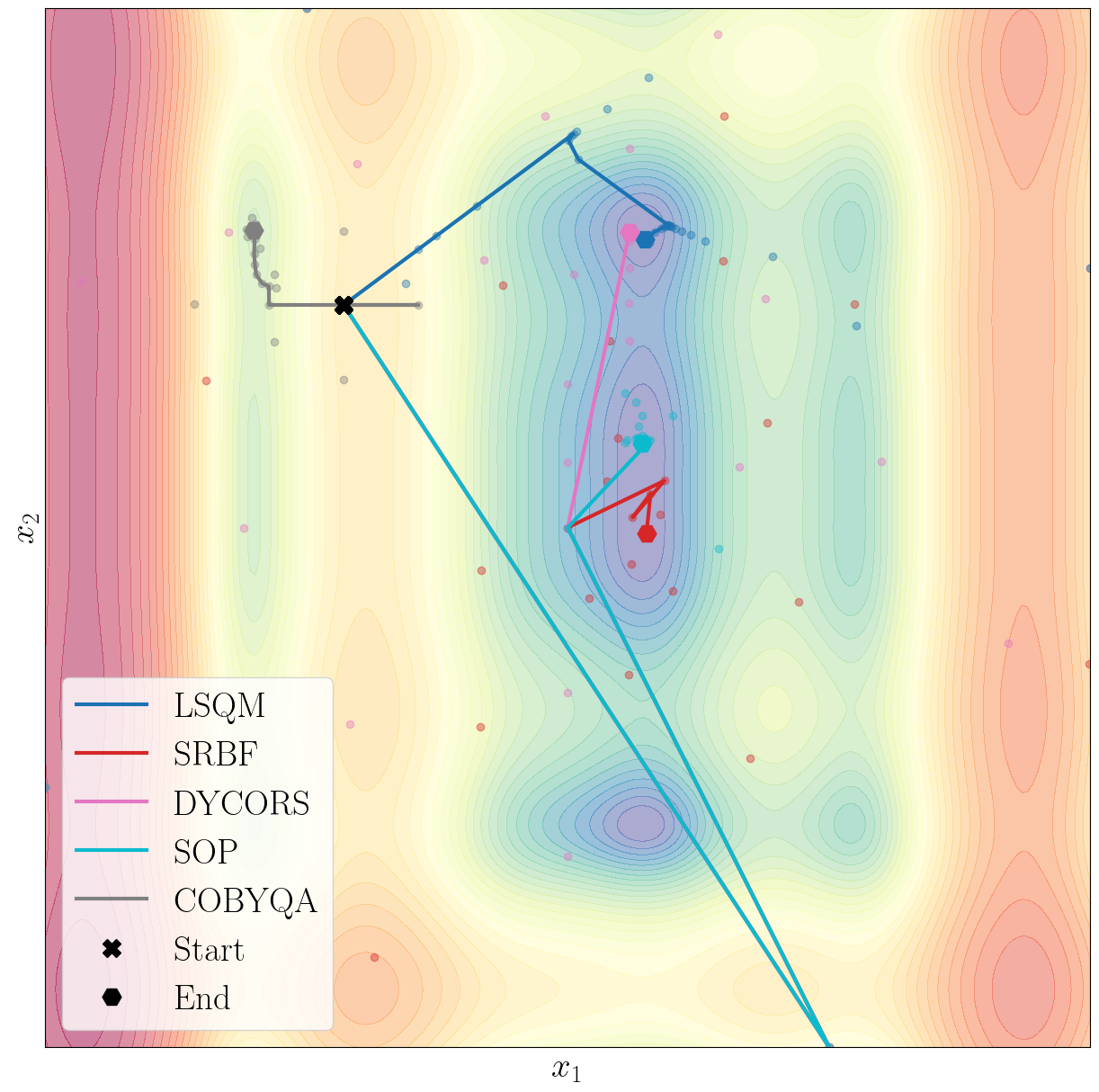}
            \caption{Levy Contour and Trajectories}
    \end{subfigure}
    \quad
    \quad
    \quad
    \quad
    \quad
    \quad
    \quad
    \quad
    \begin{subfigure}{0.385\linewidth}
        \centering
            \includegraphics[width=\linewidth]{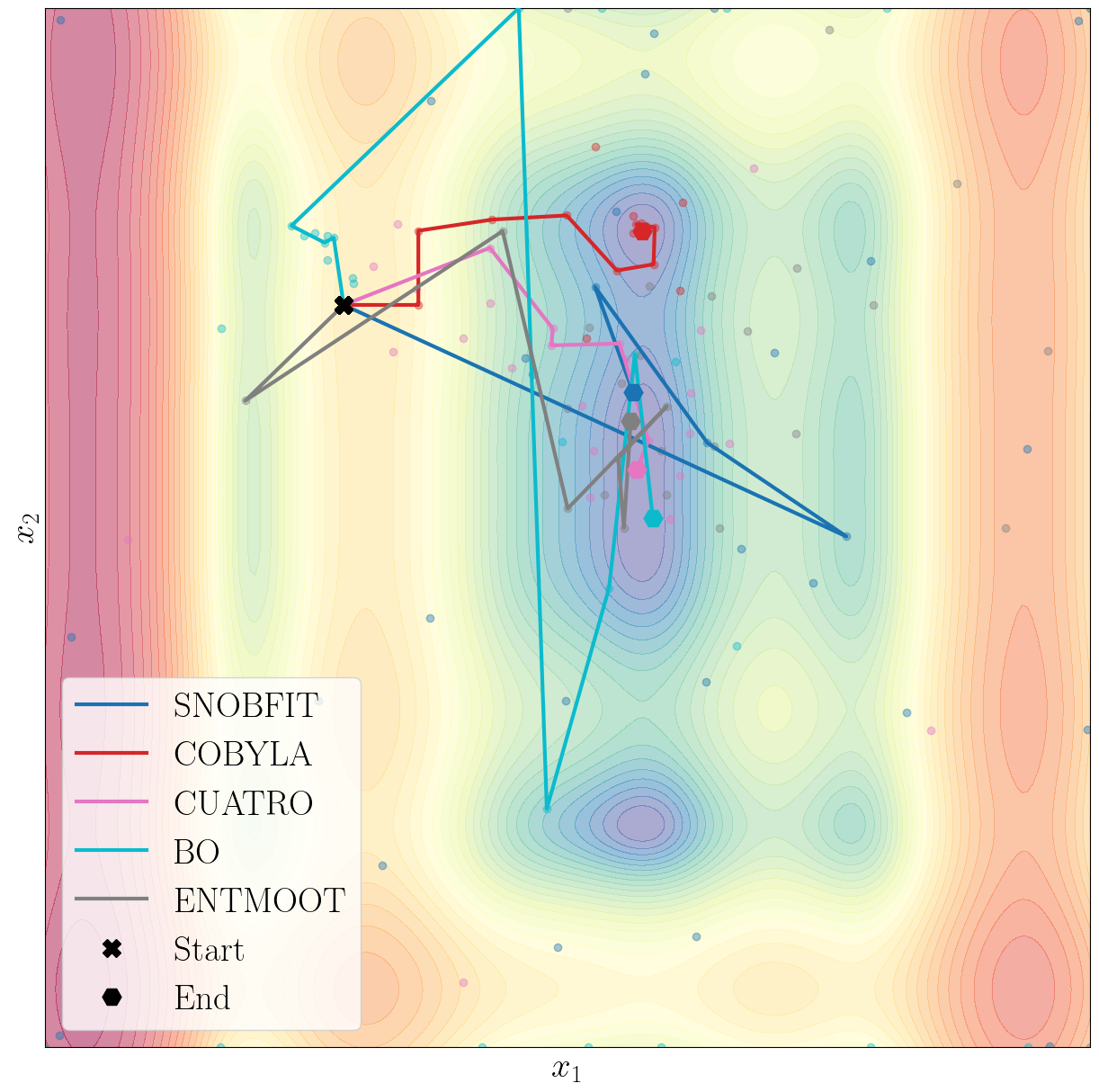}
            \caption{Levy Contour and Trajectories}
    \end{subfigure}
    \caption{2-Dimensional contour plots with exemplary optimization trajectories for the unconstrained case for each test function. The lines connect best-so-far evaluations and the scattered points show remaining evaluations.}
    \label{fig:uncon_2D_contour1}
\end{figure}
%\fi
%%%%%%%%%%%%%%%%%%%%%%%%%%%%%%%%%%%%%%%%%%%%%%%%%%%%%%%%%%%%%%%%%%%%%%%%%%%%%%%%%%%%%%%%%%%%%%%%%%%%%%%%%%%%%%%%%%%%%%%%%%%%%%%%%%%%%%%%%%%%%%%%%%%
\clearpage
%\ifincludegraphs
\begin{figure}[htbp]\ContinuedFloat
    \begin{subfigure}{0.385\linewidth}
        \centering
            \includegraphics[width=\linewidth]{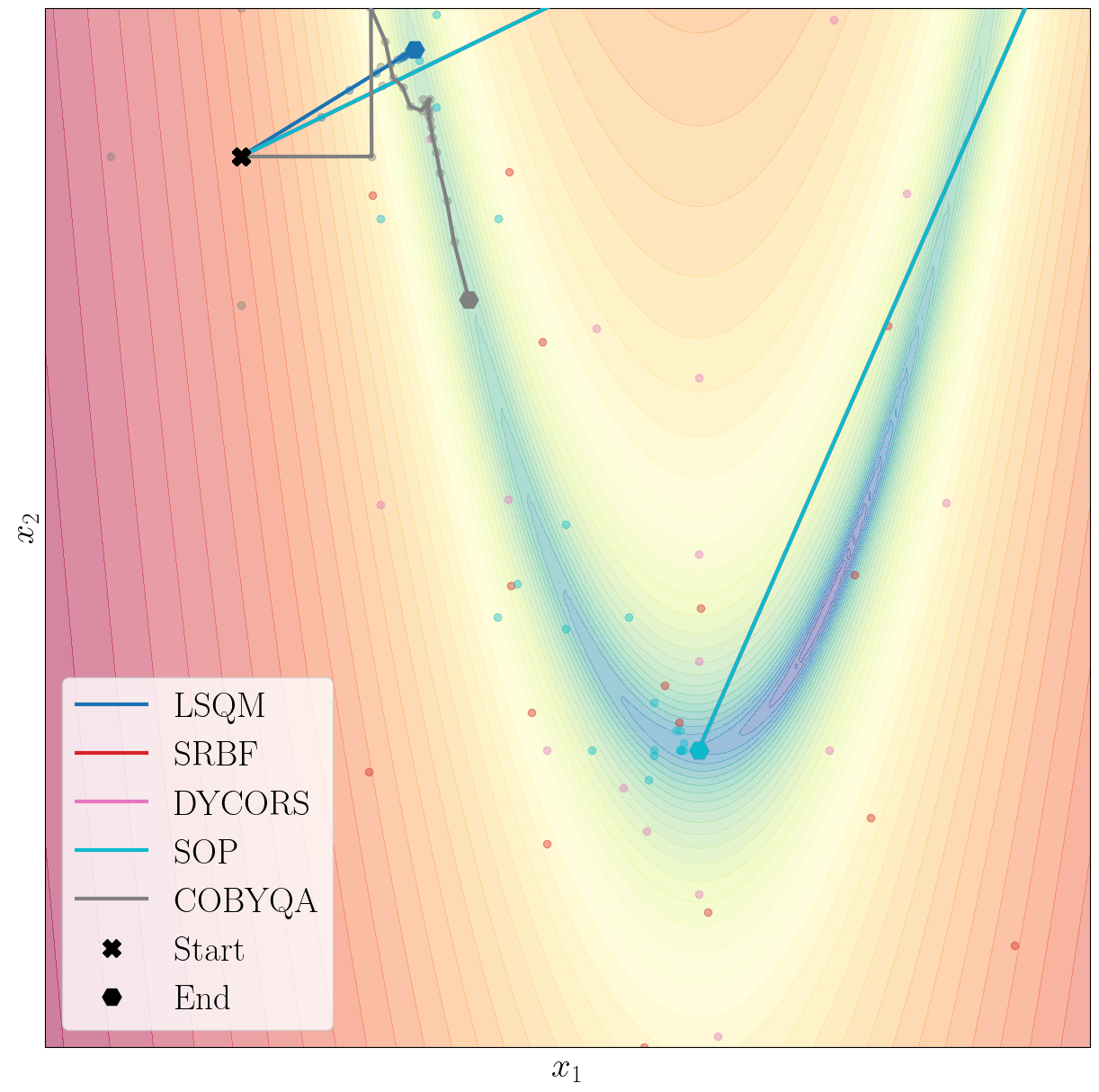}
            \caption{Rosenbrock Contour and Trajectories}
    \end{subfigure}
    \quad
    \quad
    \quad
    \quad
    \quad
    \quad
    \quad
    \quad
    \begin{subfigure}{0.385\linewidth}
        \centering
            \includegraphics[width=\linewidth]{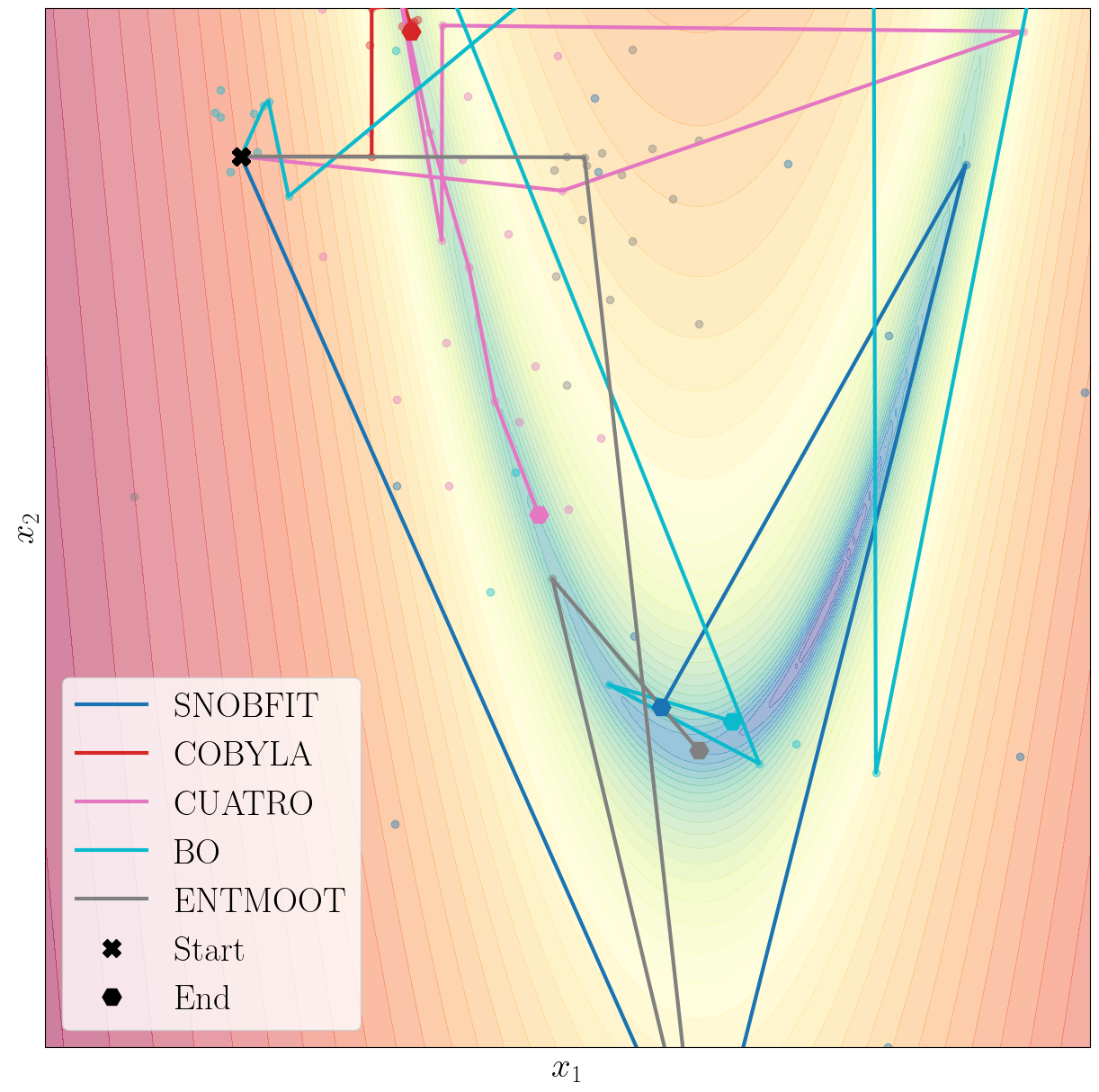}
            \caption{Rosenbrock Contour and Trajectories}
            \label{fig:uncon_2D_contour1_ros}
    \end{subfigure}
    \caption{2-Dimensional contour plots with exemplary optimization trajectories for the unconstrained case for each test function. The lines connect best-so-far evaluations and the scattered points show remaining evaluations (cont.).}    
\end{figure}
%\fi
\subsubsection{Results - Tables}
\begin{longtable}{l|ccc|ccc|c}
\caption{Quantitative algorithm scores in descending order for unconstrained optimization benchmarking} \label{tab:Uncon_BO} \\
    \hline
        \hspace{2cm}       & \textbf{Ackley} & \textbf{Levy} & \textbf{Multimodal} & \textbf{Rosenbrock} & \textbf{Quadratic} & \textbf{Unimodal} & \textbf{All} \\
    \hline
\endfirsthead
\caption[]{Quantitative algorithm scores in descending order for unconstrained optimization benchmarking (cont.)} \\
    \hline
        \hspace{2cm}       & \textbf{Ackley} & \textbf{Levy} & \textbf{Multimodal} & \textbf{Rosenbrock} & \textbf{Quadratic} & \textbf{Unimodal} & \textbf{All} \\
    \hline
\endhead
% Table content here
    \multicolumn{8}{l}{\textbf{DYCORS}} \\
    \hline
         \textbf{D2} & 0.93 & 1.00 & 0.96 & 0.89 & 1.00 & 0.95 & 0.95 \\
         \textbf{D5} & 0.94 & 1.00 & 0.97 & 0.87 & 0.95 & 0.91 & 0.94 \\
         \textbf{D7} & 0.95 & 1.00 & 0.98 & 0.93 & 0.95 & 0.94 & 0.96 \\
    \hline
         \textbf{All} & 0.94 & 1.00 & 0.97 & 0.90 & 0.97 & 0.93 & 0.95 \\
    \hline
         \multicolumn{8}{l}{\textbf{SRBF}} \\
    \hline
         \textbf{D2} & 0.96 & 0.96 & 0.96 & 0.93 & 0.98 & 0.95 & 0.96 \\
         \textbf{D5} & 0.94 & 0.86 & 0.90 & 0.97 & 0.91 & 0.94 & 0.92 \\
         \textbf{D7} & 0.92 & 0.92 & 0.92 & 0.96 & 0.93 & 0.94 & 0.93 \\
    \hline
         \textbf{All} & 0.94 & 0.91 & 0.93 & 0.95 & 0.94 & 0.94 & 0.93 \\
    \hline
             \multicolumn{8}{l}{\textbf{SOP}} \\
    \hline
         \textbf{D2} & 0.87 & 0.93 & 0.90 & 0.97 & 0.99 & 0.98 & 0.94 \\
         \textbf{D5} & 1.00 & 0.58 & 0.79 & 0.89 & 0.97 & 0.93 & 0.86 \\
         \textbf{D7} & 1.00 & 1.00 & 1.00 & 0.96 & 0.92 & 0.94 & 0.97 \\
    \hline
         \textbf{All} & 0.96 & 0.84 & 0.90 & 0.94 & 0.96 & 0.95 & 0.92 \\
    \hline
             \multicolumn{8}{l}{\textbf{BO}} \\
    \hline
         \textbf{D2} & 1.00 & 0.91 & 0.95 & 1.00 & 1.00 & 1.00 & 0.98 \\
         \textbf{D5} & 0.87 & 0.53 & 0.70 & 0.85 & 0.94 & 0.90 & 0.80 \\
         \textbf{D7} & 0.76 & 0.38 & 0.57 & 0.79 & 0.85 & 0.82 & 0.69 \\
    \hline
         \textbf{All} & 0.88 & 0.61 & 0.74 & 0.88 & 0.93 & 0.91 & 0.82 \\
    \hline
         \multicolumn{8}{l}{\textbf{COBYLA}} \\
    \hline
         \textbf{D2} & 0.79 & 0.43 & 0.61 & 0.86 & 0.92 & 0.89 & 0.75 \\
         \textbf{D5} & 0.74 & 0.74 & 0.74 & 1.00 & 1.00 & 1.00 & 0.87 \\
         \textbf{D7} & 0.74 & 0.56 & 0.65 & 1.00 & 1.00 & 1.00 & 0.83 \\
    \hline
         \textbf{All} & 0.76 & 0.58 & 0.67 & 0.95 & 0.97 & 0.96 & 0.81 \\
    \hline
             \multicolumn{8}{l}{\textbf{COBYQA}} \\
    \hline
         \textbf{D2} & 0.83 & 0.00 & 0.41 & 0.88 & 0.81 & 0.84 & 0.63 \\
         \textbf{D5} & 0.88 & 0.57 & 0.72 & 0.85 & 0.92 & 0.89 & 0.81 \\
         \textbf{D7} & 0.89 & 0.72 & 0.80 & 0.94 & 0.86 & 0.90 & 0.85 \\
    \hline
         \textbf{All} & 0.86 & 0.43 & 0.65 & 0.89 & 0.86 & 0.88 & 0.76 \\
    \hline
         \multicolumn{8}{l}{\textbf{ENTMOOT}} \\
    \hline
         \textbf{D2} & 0.55 & 0.88 & 0.71 & 0.74 & 0.73 & 0.74 & 0.73 \\
         \textbf{D5} & 0.46 & 0.80 & 0.63 & 0.96 & 0.44 & 0.70 & 0.66 \\
         \textbf{D7} & 0.18 & 0.27 & 0.22 & 0.63 & 0.00 & 0.31 & 0.27 \\
    \hline
         \textbf{All} & 0.40 & 0.65 & 0.52 & 0.78 & 0.39 & 0.58 & 0.55 \\
    \hline
         \multicolumn{8}{l}{\textbf{CUATRO}} \\
    \hline
         \textbf{D2} & 0.08 & 0.44 & 0.26 & 0.87 & 0.00 & 0.43 & 0.35 \\
         \textbf{D5} & 0.00 & 0.72 & 0.36 & 0.45 & 0.12 & 0.28 & 0.32 \\
         \textbf{D7} & 0.00 & 0.16 & 0.08 & 0.00 & 0.37 & 0.18 & 0.13 \\
    \hline
         \textbf{All} & 0.03 & 0.44 & 0.23 & 0.44 & 0.16 & 0.30 & 0.27 \\
    \hline
             \multicolumn{8}{l}{\textbf{LSQM}} \\
    \hline
         \textbf{D2} & 0.28 & 0.26 & 0.27 & 0.41 & 0.71 & 0.56 & 0.41 \\
         \textbf{D5} & 0.09 & 0.16 & 0.12 & 0.05 & 0.21 & 0.13 & 0.12 \\
         \textbf{D7} & 0.02 & 0.00 & 0.01 & 0.22 & 0.04 & 0.13 & 0.07 \\
    \hline
         \textbf{All} & 0.13 & 0.14 & 0.13 & 0.23 & 0.32 & 0.27 & 0.20 \\
    \hline
         \multicolumn{8}{l}{\textbf{SNOBFIT}} \\
    \hline
         \textbf{D2} & 0.00 & 0.45 & 0.22 & 0.00 & 0.12 & 0.06 & 0.14 \\
         \textbf{D5} & 0.02 & 0.00 & 0.01 & 0.00 & 0.00 & 0.00 & 0.00 \\
         \textbf{D7} & 0.18 & 0.17 & 0.18 & 0.17 & 0.41 & 0.29 & 0.23 \\
    \hline
         \textbf{All} & 0.06 & 0.21 & 0.14 & 0.06 & 0.18 & 0.12 & 0.13 \\
    \hline
\end{longtable}

%For ease of interpretation find attached Table \ref{tab:uncon_ALL} ranking and concluding the overall performance values for each algorithm:
%\begin{table}[H]
%    \centering
%    \captionof{table}{Overall algorithm ranking for unconstrained %optimization benchmarking} \label{tab:uncon_ALL}
%    \begin{tabular}{lcc}
%    \hline
%        \textbf{Algorithm} & \textbf{Overall Performance (All)} & %\textbf{Rank} \\
%    \hline
%         \textbf{DYCORS} & 0.95 & 1  \\
%         \textbf{SRBF} & 0.93 & 2  \\
%         \textbf{SOP} & 0.92 & 3  \\
%         \textbf{BO} & 0.82 & 4  \\
%         \textbf{COBYLA} & 0.81 & 5  \\
%         \textbf{COBYQA} & 0.76 & 6  \\
%         \textbf{ENTMOOT} & 0.55 & 7  \\
%         \textbf{CUATRO} & 0.27 & 8  \\
%         \textbf{LSQM} & 0.20 & 9  \\
%         \textbf{SNOBFIT} & 0.13 & 10  \\
%    \hline
%    \end{tabular}
%\end{table}

\subsubsection{Results and Discussion: Mathematical Unconstrained Functions}
This section presents the conclusions and discussion of results from the unconstrained synthetic benchmarking. The performance assessment considers varying dimensionality and includes both unimodal and multimodal test functions. The interpretations drawn here are based on the quantitative benchmarking ranking and trajectory plots. It is important to note that all conclusions and results are relative to the specific set of algorithms included in this performance assessment. Overall, DYCORS emerged as the best-performing algorithm in the unconstrained synthetic benchmarking. SRBF secured the second-best position, demonstrating consistent and robust performance. On the other hand, SNOBFIT was the poorest-performing algorithm, with LSQM ranking as the second poorest. It is important to highlight that very competitive algorithms were benchmarked, so poor performance does not mean an algorithm is bad. Additionally, 4 test functions across 3 dimensions is not an exhaustive combination, therefore readers are encouraged to draw only directional conclusions, and even make their own assessment if they wish to use these algorithms for a specific application. Furthermore, Figure \ref{fig:uncon_2D_contour1_ros} shows erratic behavior from some of the algorithms. This highlights the difference between gradient-based algorithms, when derivatives are available, and derivative-free methods, which have to explore and "learn" the function as well as optimize it. Algorithm-specific conclusions have been grouped and presented below.
\\[8pt]
\textbf{BO:} The in-house implementation of BO was the best-performing algorithm on the 2D unimodal test and the Ackley function in 2D. However, a decline in performance was observed with increasing dimensionality of the test functions.
\\[8pt]
\textbf{COBYQA:} COBYQA performed the worst on the 2D Levy function. It showed better performance on unimodal test functions than multimodal ones, with performance increasing with dimensionality across all test functions except the Ill-Conditioned Quadratic function.
\\[8pt]
\textbf{LSQM:} LSQM was the second worst-performing algorithm in the unconstrained synthetic category. It performed slightly better on unimodal test functions, with the best performance observed on the Rosenbrock function.
\\[8pt]
\textbf{COBYLA:} COBYLA demonstrated strong overall performance, excelling in unimodal test functions over multimodal ones. It achieved the best results on the 2D and 7D unimodal tests, particularly on the Rosenbrock and Ill-Conditioned Quadratic functions.
\\[8pt]
\textbf{CUATRO:} CUATRO performed the worst in Ackley D5 and D7, Ill-Conditioned Quadratic 2D, and Rosenbrock 7D, but exhibited good performance in the Rosenbrock 2D test.
\\[8pt]
\textbf{DYCORS:} DYCORS emerged as the top-performing algorithm in the unconstrained synthetic benchmarking and was the best among the three RBF-based Shoemaker et al. algorithms. It excelled in the Levy test function for D2, D5, and D7, as well as in the Ill-Conditioned Quadratic D2 test. Increasing dimensionality did not significantly affect its performance within the tested range.
\\[8pt]
\textbf{ENTMOOT:} ENTMOOT had similar performance in both unimodal and multimodal tests, excelling in the Rosenbrock D5 function but performing the worst in the Ill-Conditioned Quadratic D7. Its performance decreased in all test functions except for the Rosenbrock function.
\\[8pt]
\textbf{SNOBFIT:} SNOBFIT was the worst-performing algorithm in the unconstrained synthetic category, with the poorest results in Ill-Conditioned Quadratic D2, Rosenbrock D2, and D5, Levy D5, and being the worst in Unimodal D5.
\\[8pt]
\textbf{SOP:} SOP showed strong performance across both unimodal and multimodal test functions, with performance not significantly impacted by increasing dimensionality. It was the best in Levy D7, Ill-Conditioned Quadratic D5, and D7, and had the best overall performance in Multimodal D7.
\\[8pt]
\textbf{SRBF:} SRBF was a solid all-rounder and the second-best algorithm in the overall unconstrained synthetic benchmarking. Although it did not rank first in any single category, it scored over 0.9 in all but the Levy test function in D5. It had a very small standard deviation in performance, and increasing dimensionality did not significantly impact its results.

\newpage
\section{Surrogate-based Constrained Optimization}\label{sec:c}
Constrained optimization is pivotal in process optimization, as constraints play a critical role in ensuring safety, operational, and environmental constraints. From a safety standpoint, constraints ensure the process operates within safe operating limits, preventing potential hazards such as chemical leaks, fires, or explosions. Operationally, they maintain process efficiency and prevent equipment breakdown by establishing clear boundaries for variables such as flow rates, reaction times, and equipment utilization. By doing so, constraints mitigate equipment damage, minimize downtime, and ensure consistent product quality. Finally, environmental constraints are critical and equally vital, reducing the impact of chemical processes on the environment. Highlighting their significance, in practical terms, every process invariably operates at the limit of at least one constraint - this is driven by the demands of a competitive global economy, which pushes processes to achieve peak performance. Hence, it is imperative to develop algorithms that are capable of optimizing such problems whilst also satisfying constraints.
\\[8pt]
Given the importance of constraints, this section will outline data-driven model-based optimization algorithms that can deal with constrained problems. It will also present the theory behind the algorithms as well as a comparative study for the reader to develop an intuition and understanding of the different algorithms.  

\subsection{Problem Formulation for Constrained Surrogate-Based Optimization}

Recall the general optimization formulation in Equation \ref{eq:basic_opt} with modifications for a constrained system:
\begin{equation}
\begin{aligned}
     & \min_{\mathbf{x}} && f(\mathbf{x}) \\
     & \text{s.t. } && \mathbf{g}(\mathbf{x}) \le \mathbf{0} \\
     &&& \mathbf{x} \in \mathbb{R}^{n_x} \\
     &&& f:\mathbb{R}^{n_x} \longrightarrow \mathbb{R} \\
     &&& \mathbf{g}: \mathbb{R}^{n_x} \longrightarrow \mathbb{R}^{n_g}
\end{aligned}
\end{equation}

The challenge arises from the lack of prior knowledge regarding the locations where the constraint equations $\mathbf{g}(\cdot)$ are satisfied. The handling of constraints varies depending on the safety criticality of the system. For instance, when sampling a simulation, the primary objective is to ensure that the resulting solution complies with the constraints. However, when sampling a physical reactor, the goal might be to minimize both the number and severity of constraint violations during sampling. In scenarios such as simulations without safety considerations, an intuitive way to deal with constraints might be via the penalty method where violations are penalized using the Lp-norm:
\begin{align}
     & \min_{\mathbf{x}} f(\mathbf{x}) + \rho \sum_{j=1}^{n_g} (\max(0,g_j(\mathbf{x})))^p
\end{align}
where $\rho$ constitutes the penalty parameter. However, the choice of penalty parameter $\rho$ is not trivial: if too small, the returned solution might accept small constraint violations in favor of increases in $f(\cdot)$; if too big, the problem might become ill-conditioned, and difficult to solve with model-based methods. 
\\[8 pt]
Typically, it is assumed that both $f(\mathbf{x})$ and $\mathbf{g}(\mathbf{x})$ can only be sampled. Regarding $\mathbf{g}(\mathbf{x})$, there are two scenarios: in the first, constraints may be violated during optimization as long as the optimal candidate meets the constraints. In the second scenario, constraints cannot be violated during the optimization sequence, nor by the optimal candidate. The first is referred to as constrained optimization, and the second is safe optimization. Safe optimization is particularly relevant in real-world chemical engineering problems where, for example, samples could represent experiments conducted in a reactor. For instance, if the constraint $T \leq 500$ [K] represents the safe operating region that avoids runaway reactions, safe optimization is most suitable for ensuring process safety and stability.
\\[8pt]
In what follows, the focus shifts to the more interesting "safe" scenario of minimizing both the number and severity of constraint violations. This necessitates the adoption of a distinct category of explicit constraint-handling methods.

%regarding the constraint handling: Damien mentions in his paper, that CUATRO and BO are the only methods that handle constraints explicitly, by estimating the constraints using convex quadratic surrogates and by directly interpolating the mean respectively. For all other methods, constraints are handled using the penalty method of the Lp norm, usually the L2 norm 

\subsection{Explicit Constraint Handling Methods}
All the constrained DFO solvers studied in this performance assessment handle constraints \textit{explicitly}. These methods not only construct surrogates of the objective function $\hat{f}$ but also for the constraint evaluations: $\hat{g}_j(\mathbf{x}) \approx g_j(\mathbf{x})$. These surrogate models are then incorporated into the minimization step 
%[LINK TO MINIMIZATION STEP IN EARLIER WRITING] 
as constraints, ensuring that the next sample is expected to satisfy the constraints before the evaluation. To address uncertainties that arise from inexact surrogates and noisy constraint evaluations, some considered methods also employ trust regions $T(\cdot)$ to restrict the update step to the vicinity of the current (feasible) candidate, where the algorithms are more confident in the accuracy of its surrogates.

\begin{equation}
    \begin{aligned}
     & \min_{\mathbf{x} \in \mathbb{R}^{n_x} } && \hat{f}(\mathbf{x}) \\
     & \text{s.t. } && \hat{\mathbf{g}}(\mathbf{x}) \le \mathbf{0} \\
     &&& \mathbf{x} \in T(\mathbf{x}_c) 
     \end{aligned}
\end{equation}

where $\mathbf{x}_c$ is a "safe" center point, and $\hat{\mathbf{g}}$ could either map a different surrogate for each constraint or map a single surrogate on the maximum constraint violation. The methods benchmarked mostly differ in the type of surrogate $\hat{\mathbf{g}}$ employed. 

In the following sections, the algorithms for constrained model-based DFO are outlined. Specifically, how the existing algorithms are modified (if in any way) to account for constraints. 

%\todo[inline]{check that the outline sections for constrained and unconstrained match ... or are at least as close as possible}

\subsection{Constrained Surrogate Optimization Methods}
CUATRO, COBYLA and COBYQA account for constraints by construction, see sections \ref{sec:cuatro}, \ref{sec:cobyla}, and \ref{sec:cobyqa}, respectively for more details. However, not all model-based derivative free surrogates account for constraint by construction, for example, BO. In the following section, we will describe how constraints can be incorporated into the BO framework. This approach can also serve as a template for incorporating constraints in other surrogate-based optimization methodologies.

%\subsubsection{CUATRO}
%CUATRO, short for Convex qUadratic Trust Region Optimizer, is similar to COBYQA, in that it combines trust region with quadratic surrogates with the main difference being that CUATRO uses \textit{convex} quadratic surrogates within a convex programming framework. Additionally, instead of training separate surrogates for each constraint, we use CUATRO with the option to construct a single constraint surrogate via quadratic discrimination where we only on the binary (in)feasibility labels of all samples in the trust region. It also incorporates a safe exploration routine by sampling a point in the trust region that is furthest away from all other samples while still expected to be feasible. A detailed description can be found in \cite{vandebergDatadrivenOptimizationProcess2022}.

\subsubsection{Constraint handling in Bayesian optimization}

In \textit{Constrained Bayesian Optimization} (CBO), constraints can be handled through various methods. A common approach involves constructing surrogates of the constraints, typically in the form of GPs. Since the constraints are unknown but can be sampled like the objective function, surrogates can be built following a similar procedure to that of the objective function.

There are two usual ways to utilize these surrogate constraints. The first approach involves penalizing constraint violations in the objective function. This can be done by either multiplying the objective function by one minus the probability of violation or adding the constraint violation as a penalty. The second common approach involves using the surrogate constraints to formulate a constrained optimization problem and solving it using constrained optimization algorithms, such as interior point methods or sequential quadratic programming.
Specifically, given sampled inputs $X = \left[ \mathbf{x}^{(1)}, ..., \mathbf{x}^{(n_d)} \right]^\intercal$, BO constructs a GP model $\mathcal{GP}_f$ to predict the outputs $\mathbf{y}_{f} = \left[ f \left(\mathbf{x}^{(1)}\right), ..., f \left(\mathbf{x}^{(n_d)}\right) \right]$. Similarly, a GP $\mathcal{GP}_{g_i}$ can be constructed to model each constraint, provided it can be measured, resulting in $\mathbf{y}_{g_i} = \left[ g_i \left(\mathbf{x}^{(1)}\right), ..., g_i \left(\mathbf{x}^{(n_d)}\right) \right]$. Consequently, the following optimization problem can be formulated:

\begin{equation}
    \begin{aligned}
     & \min_{\mathbf{x} \in \mathbb{R}^{n_x} } && \mathcal{A}_{\mathcal{GP}_f}(\mathbf{x}) \\
     & \text{s.t. } && \mathbf{g}_\mathcal{GP}(\mathbf{x}) \le \mathbf{0} \\
     \end{aligned}
\end{equation}

where $\mathcal{A}_{\mathcal{GP}_f}$ denotes the acquisition function using the GP model of the objective function (e.g., $\mu_f(\mathbf{x}) - \gamma \sigma^2_f(\mathbf{x})$), and $\mathbf{g}_\mathcal{GP}$ denotes the use of the GPs that model the constraints to ensure feasibility. One straightforward way to incorporate them is to use the mean of the GP for each constraint:

\begin{equation}
    \label{eq:constrained_BO}
    \begin{aligned}
     & \min_{\mathbf{x} \in \mathbb{R}^{n_x} } && \mathcal{A}_{\mathcal{GP}_f}(\mathbf{x}) \\
     & \text{s.t. } && \mu_{\mathcal{GP}_i}(\mathbf{x}) \le 0 \quad \text{for}~i=1,...n_g\\
     \end{aligned}
\end{equation}

While this can be effective, using only the mean of a GP defeats the purpose of this surrogate and it could be replaced by less computationally expensive surrogates such as RBFs. One way in which the variance term can be incorporated is as a backoff to add an extra layer of safety, for example:

\begin{equation}
    \begin{aligned}
     & \min_{\mathbf{x} \in \mathbb{R}^{n_x} } && \mathcal{A}_{\mathcal{GP}_f}(\mathbf{x}) \\
     & \text{s.t. } && \mu_{\mathcal{GP}_i}(\mathbf{x})+\sigma_{\mathcal{GP}_i}(\mathbf{x}) \le 0 \quad \text{for}~i=1,...n_g\\
     \end{aligned}
\end{equation}

Finally, if safe exploration is important a trust region can be incorporated \cite{chanonaRealtimeOptimizationMeets2021}.
In this book chapter, the implementation follows the version in Equation \ref{eq:constrained_BO} to maintain consistency with the rest of the methods described.
%\subsubsection{Alternative methods from BO literature}
%multiplying acquisition function by probability of constraint satisfaction:
%Barnett Chapter 11.2 - (11.6)

\subsection{Performance Assessment of Constrained Surrogate-based Optimization Algorithms}

For the constrained performance assessment of black-box problems the functions and algorithms are as follows:
\begin{align*}
    a &\in \mathbb{A}, & \mathbb{A} &= \left\{ \text{CBO, COBYLA, COBYQA, CUATRO} \right\} \\
    f &\in \mathbb{F}, & \mathbb{F} &= \left\{ \text{Constrained Rosenbrock, Constrained Quadratic, Constrained Matyas} \right\}.
\end{align*}

Therefore every algorithm $a \in \mathbb{A}$ is assessed on every constrained problem $f \in \mathbb{F}$ and its performance is compared relative to the other optimization algorithms on the same constrained black-box optimization problem.

\subsubsection{Mathematical Objective Functions}
The objective functions used in this subsection are again the ill-conditioned quadratic function and the Rosenbrock function presented in Section \ref{sec:umof}, as well as the Matyas function proposed by Matyas in 1965 \cite{matyasRandomOptimization1965}. The input dimension chosen for this section is $n_x = 2$. Subsequently, the objective functions and the constraints can be presented as:
\begin{itemize}
    \item Rosenbrock: 
    \begin{equation} \label{eq:Rosenbrock_con}
        \begin{aligned}[t]
            &\min &&f(x_1,x_2) = (1-x_1)^2 + 100(x_2-x_1^2)^2 \\
            & \text{s.t. } &&g(x_1, x_2) = x_1 + 1.27 - 2.83x_2 + 0.69x_2^2
        \end{aligned}
    \end{equation}
    \item Quadratic: 
    \begin{equation} \label{eq:Quadratic_con}
        \begin{aligned}[t]
            &\min &&f(x_1,x_2) = x_1^2 + 0.95x_1x_2 + 5.9x_2^2 \\
            & \text{s.t. } &&g(x_1, x_2) = 1.5x_1 + 0.6 - x_2
        \end{aligned}
    \end{equation}
    \item Matyas: 
    \begin{equation} \label{eq:Matyas_con}
        \begin{aligned}[t]
            &\min &&f(x_1,x_2) = 0.26(x_1^2 + x_2^2) - 0.48x_1x_2 \\
            & \text{s.t. } &&g(x_1, x_2) = 6.31x_1 + 3.60 - x_2
        \end{aligned}
    \end{equation}
\end{itemize}

It can be observed that when the algorithms commence evaluating alongside the constraint, there are very small constraint violations e.g. 0.0008. Therefore, the threshold of 0.001 was chosen to determine a constraint violation as such. This means that when the constraint is violated when $g(\mathbf{x})>0$ it is counted as a violation when $g(\mathbf{x})>0.001$ 

\newpage
\subsubsection{Results - Convergence Plots}
%\ifincludegraphs
\begin{figure}[H]
    \centering
    \begin{subfigure}{0.37\linewidth}
        \centering
            \includegraphics[width=\linewidth]{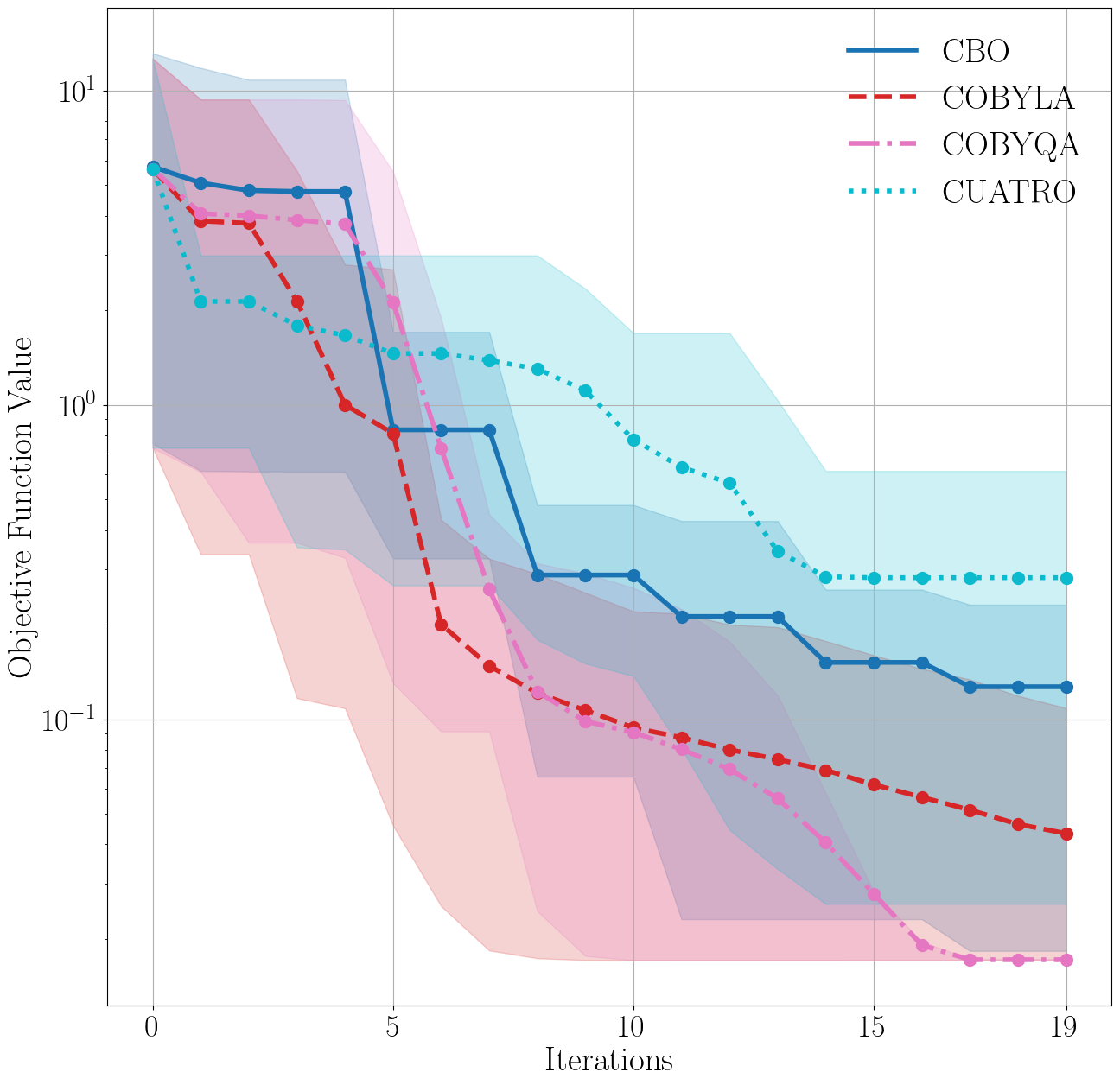}
            \caption{Matyas Objective Function}
    \end{subfigure}
    \quad
    \quad
    \begin{subfigure}{0.37\linewidth}
        \centering
            \includegraphics[width=\linewidth]{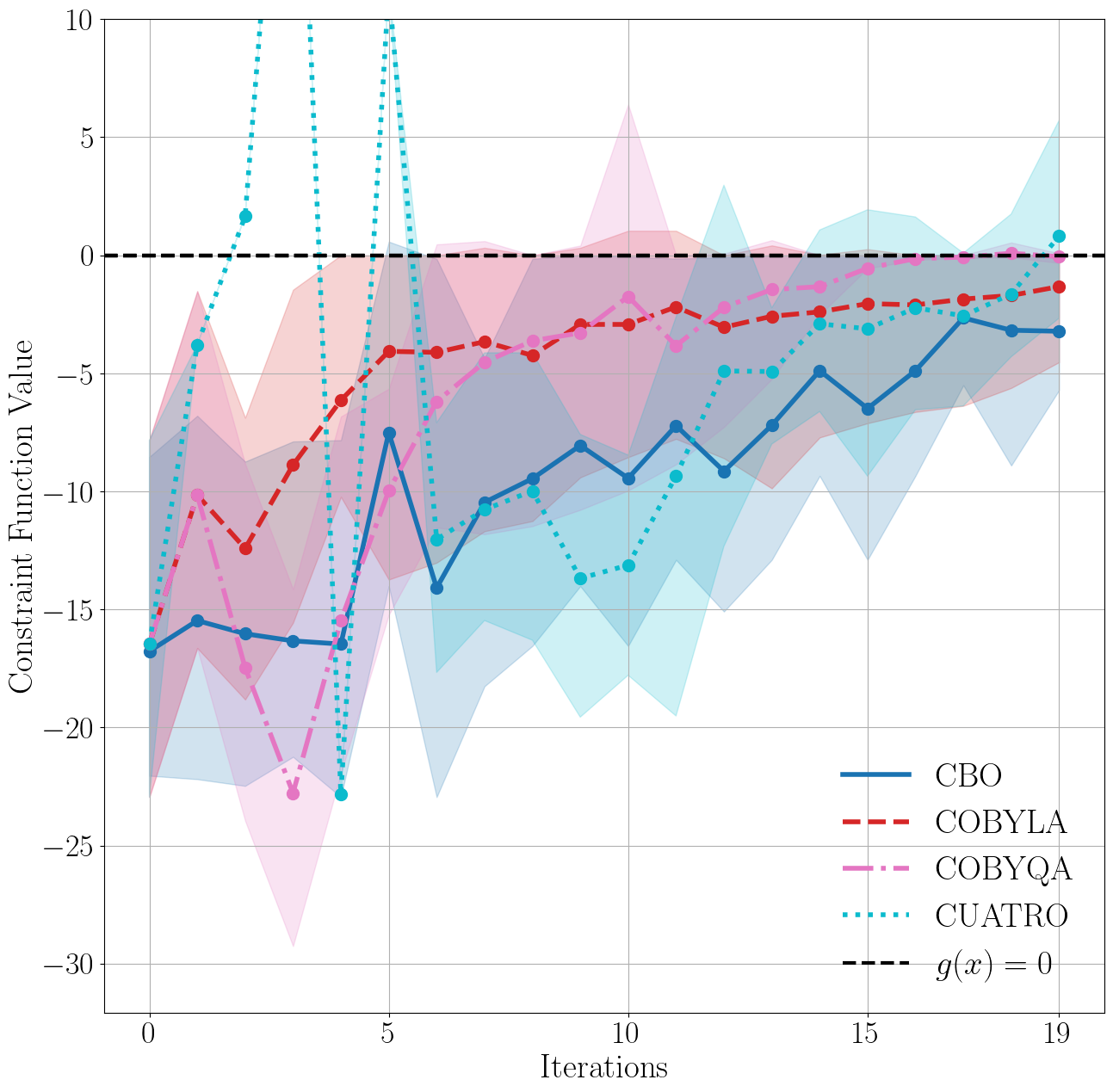}
            \caption{Matyas Constraint Function}            
    \end{subfigure}        
    \\
    \begin{subfigure}{0.37\linewidth}
        \centering
            \includegraphics[width=\linewidth]{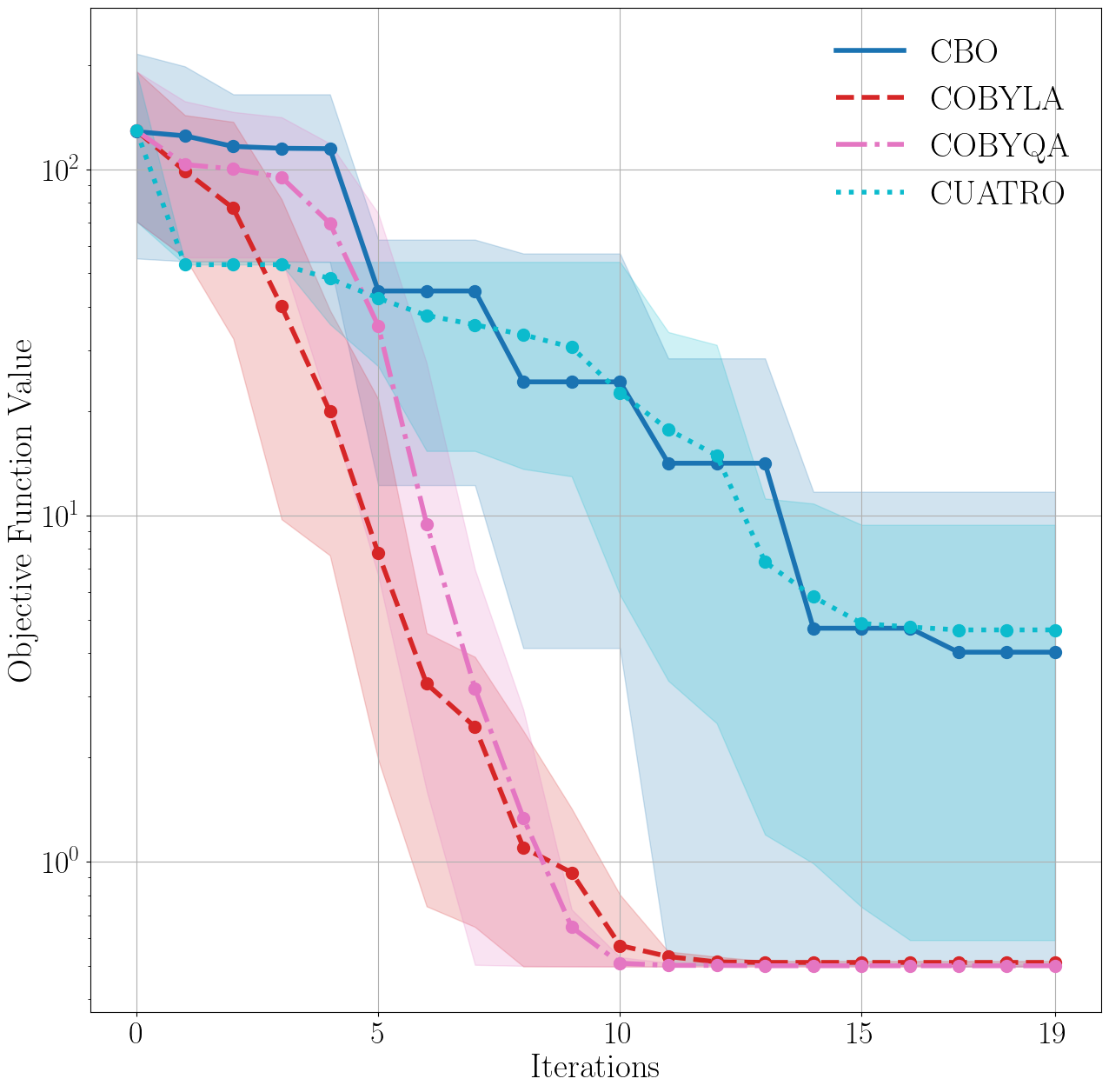}
            \caption{Quadratic Objective Function}
    \end{subfigure}   
    \quad
    \quad
    \begin{subfigure}{0.37\linewidth}
        \centering
            \includegraphics[width=\linewidth]{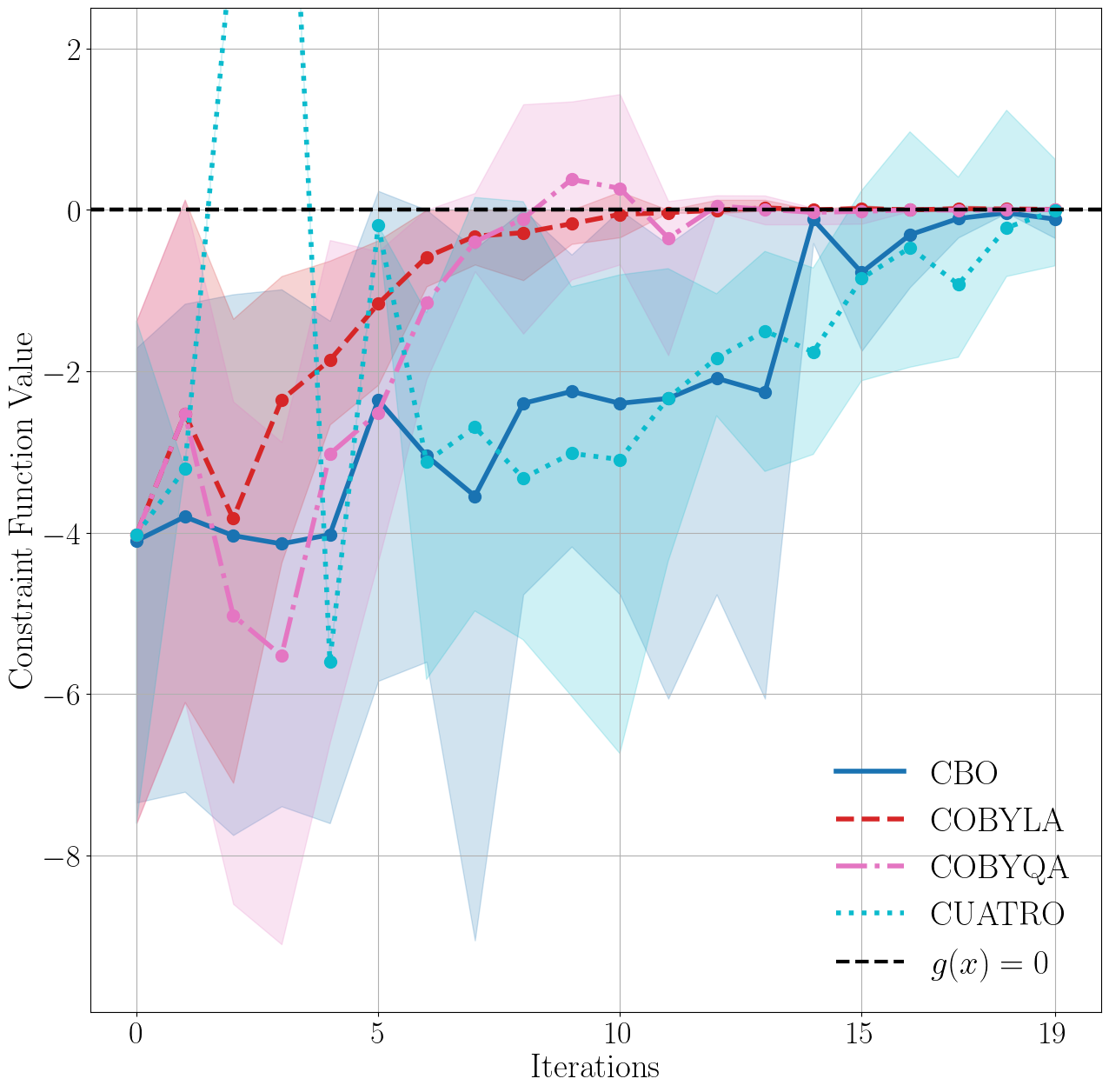}
            \caption{Quadratic Constraint Function}
    \end{subfigure}
    \\
    \begin{subfigure}{0.37\linewidth}
        \centering
            \includegraphics[width=\linewidth]{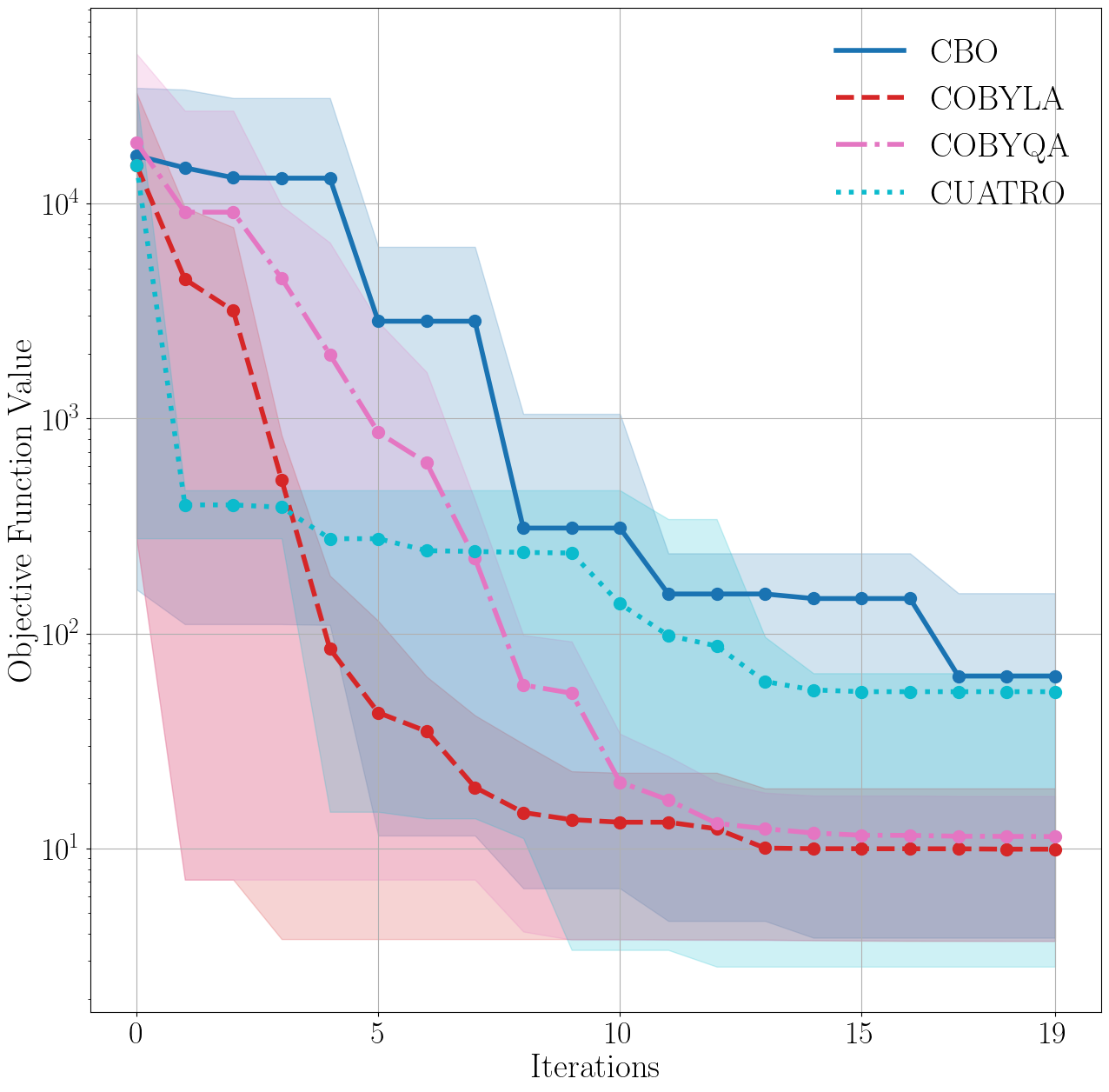}
            \caption{Rosenbrock Objective Function}
    \end{subfigure}
    \quad
    \quad
    \begin{subfigure}{0.37\linewidth}
        \centering
            \includegraphics[width=\linewidth]{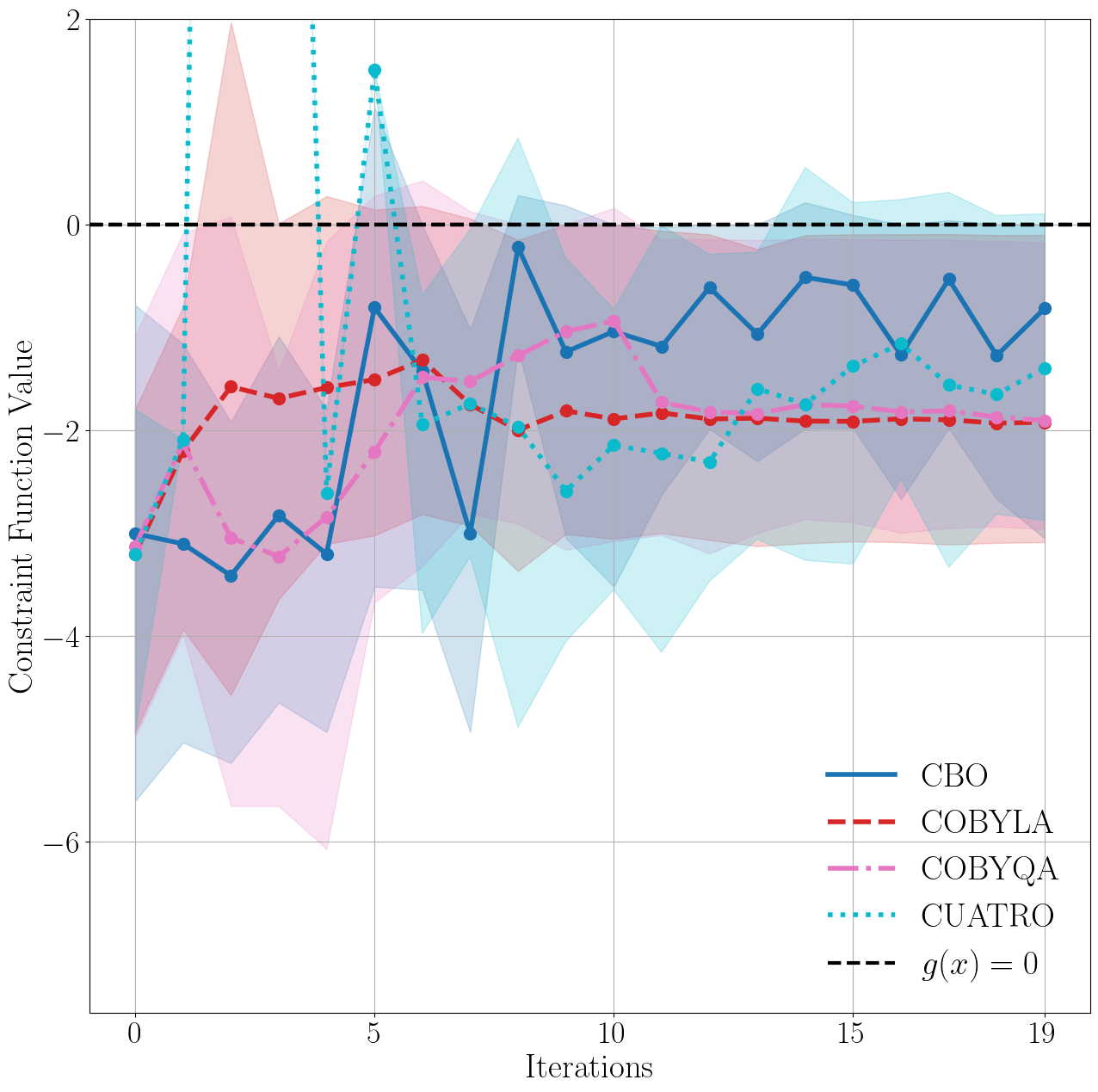}
            \caption{Rosenbrock Constraint Function}
    \end{subfigure}
    \caption{Convergence plots, showing mean objective function values (left) and constraint function values (right) and $90\%-10\%$-intervals enveloping trajectories over 10 repetitions from 10 different starting points with a budget (trajectory length) of 20 evaluations. The dashed line (right) indicates the boundary of constraint violation for $g(x)>0$. The dimensionality of the problem is 2 input dimensions and 1 output dimension for objective and constraint function.}  
    \label{fig:convergence-plots}
\end{figure}

%\fi
%%%%%%%%%%%%%%%%%%%%%%%%%%%%%%%%%%%%%%%%%%%%%%%%%%%%%%%%%%%%%%%%%%%%%%%%%%%%%%%%%%%%%%%%%%%%%%%%%%%%%%%%%%%%%%%%%%%%%%%%%%%%%%%%%%%%%%%%%%%%%%%%%%%
\newpage
\subsubsection{Results - Trajectory Plots}
\begin{figure}[H]
    \begin{subfigure}{0.33\linewidth}
        \centering
            \includegraphics[width=\linewidth]{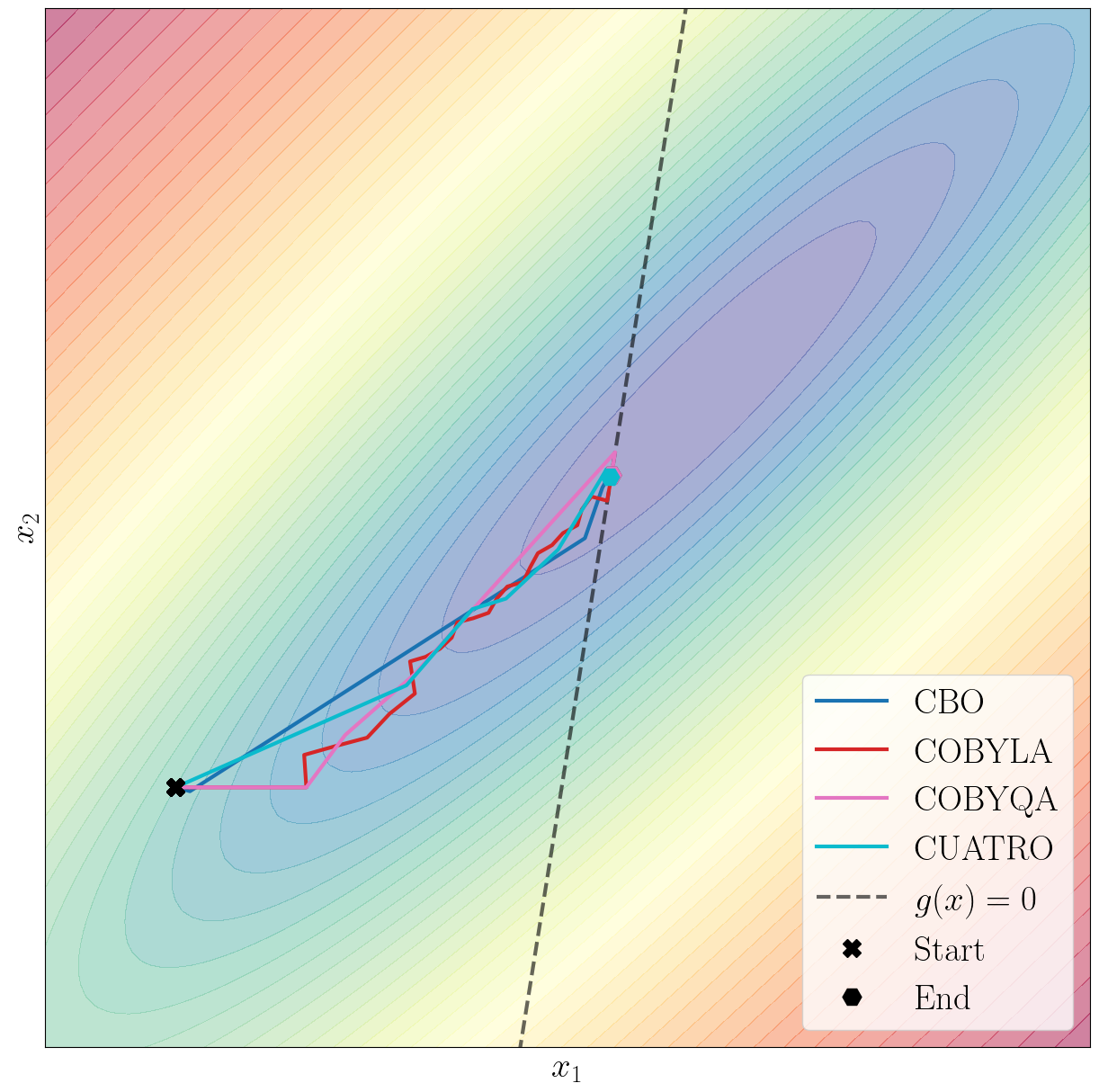}
            \caption{Constrained Matyas Objective Function}
    \end{subfigure}
    \begin{subfigure}{0.33\linewidth}
        \centering
            \includegraphics[width=\linewidth]{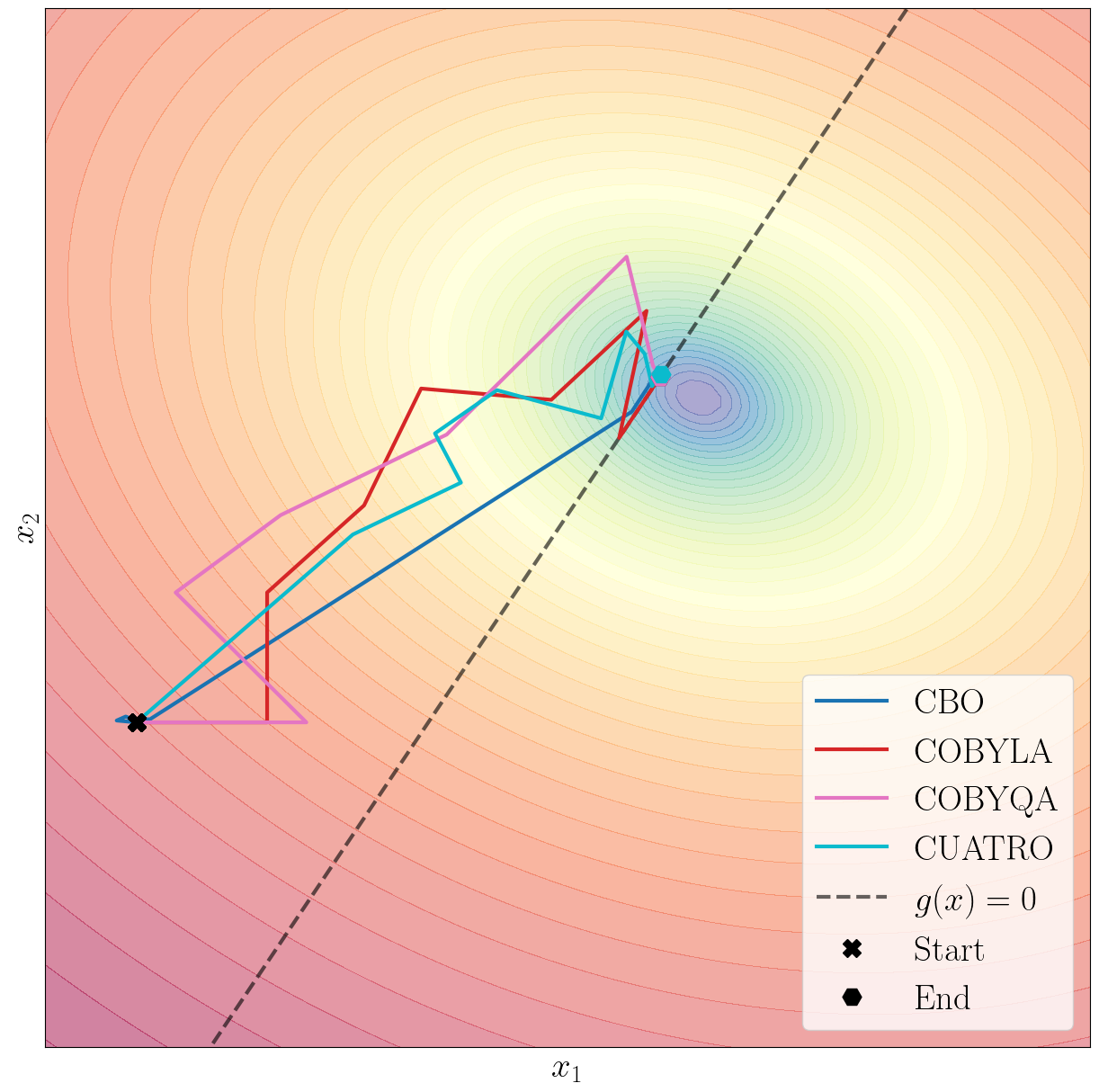}
            \caption{Constrained Quadratic Objective Function}            
    \end{subfigure}    
    \begin{subfigure}{0.33\linewidth}
        \centering
            \includegraphics[width=\linewidth]{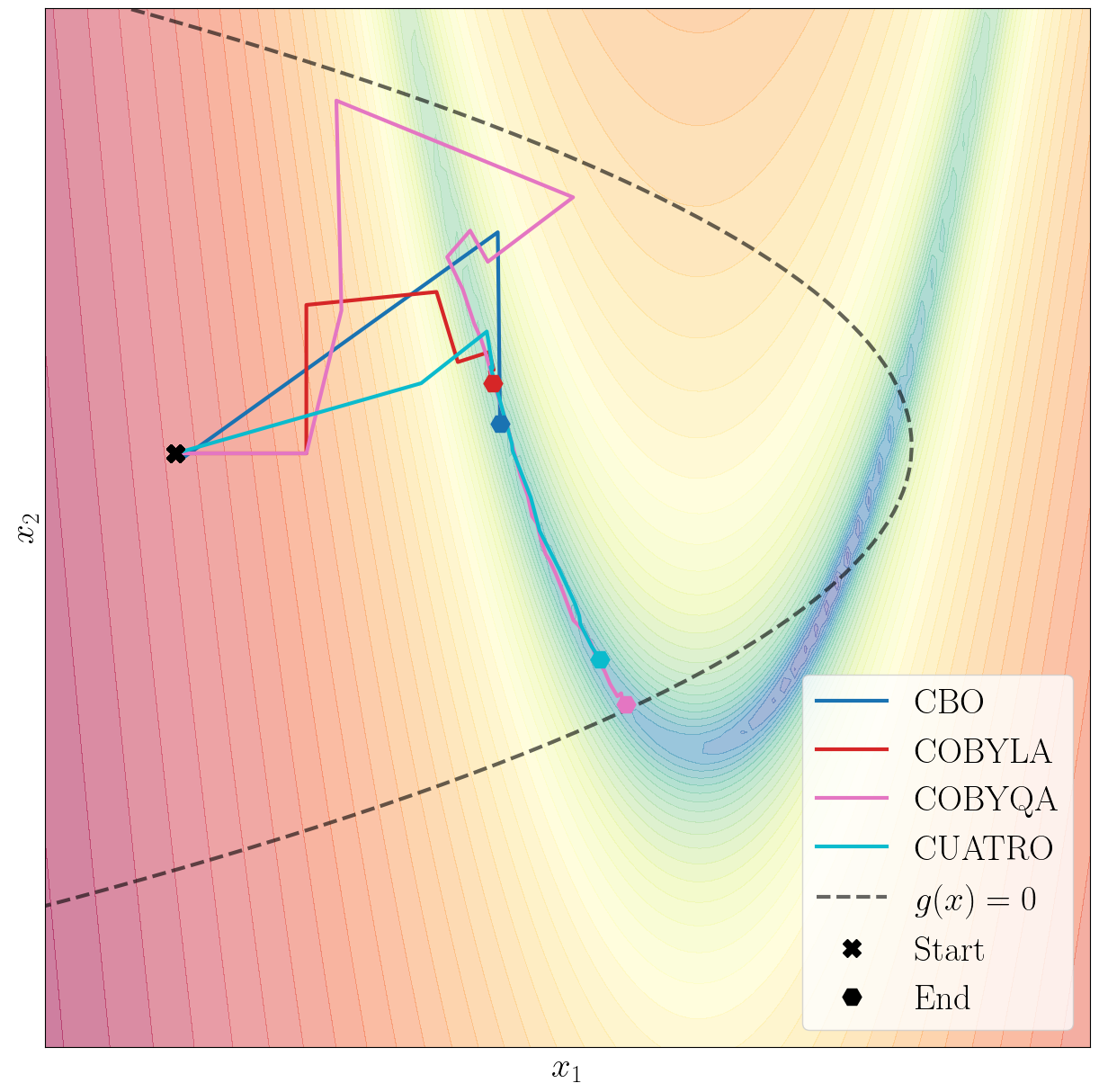}
            \caption{Constrained Rosenbrock Objective Function}            
    \end{subfigure}        
    \caption{2-Dimensional trajectory plots with constraint overlay for each constrained test function. The lines show connection of mean best-so-far evaluations of 10 repetitions per algorithm from a shared starting point and a budget of 60 evaluations}    
\end{figure}

\subsubsection{Results - Tables}

\begin{table}[H]
    \centering
    \captionof{table}{Benchmarking table for the constrained optimization}
    \begin{tabular}{lllll}
    \hline
          &\textbf{Rosenbrock} & \textbf{Quadratic} & \textbf{Matyas} & \textbf{Average}\\
    \hline
         \textbf{COBYQA} & 0.98 | 94.75\% | 0.29 & 1.00 | 88.75\% | 0.41 & 1.00 | 87.50\% | 1.26 & 0.99 | 90.33\% | 0.65 \\
         \textbf{COBYLA} & 1.00 |  94.39\% | 0.45 & 1.00 | 95.12\% | 0.10 & 0.92 |  92.20\% | 1.20 & 0.97 | 93.90\% | 0.58 \\
         \textbf{CBO} & 0.00 |  87.64\% | 0.25 & 0.04 |  97.09\% | 0.14 & 0.60 |  92.91\% | 0.46 & 0.21 | 92.55\% | 0.28 \\
         \textbf{CUATRO} & 0.49 |  87.66\% | 8.50 & 0.00 |  89.92\% | 3.05 & 0.00 |  74.24\% | 8.85 & 0.16 | 83.94\% | 6.80 \\
    \hline

    \end{tabular}
\end{table}
\vspace{-0.5cm}
\centerline{N.B. the cells are to be read as $p_a$ | constraint satisfaction | average constraint violation}
\vspace{0.5cm}
\subsubsection{Results and Discussion: Mathematical Constrained Functions}

For the constrained synthetic benchmarking, in addition to the trajectory plots and normalized algorithm scores, the constraint violations are considered through constraint violation trajectory plots. As in the unconstrained case, the algorithms' performance is observed for different test functions. It is important to note that all conclusions and results are relative to the specific set of algorithms included in this performance assessment. Algorithm-specific conclusions have been grouped and presented below.
\\[8pt]
\textbf{CBO}: The in-house implementation of CBO overall performed poorly, with the worst performance on the Rosenbrock function. The trajectory plots for the Matyas and Ill-Constrained Quadratic test functions show the largest standard deviation among all algorithms. Additionally, the objective function value remains roughly constant for the first five iterations across all test functions before improving, likely due to the algorithm using these iterations to construct its GP surrogate.
\\[8pt]
\textbf{COBYLA}: COBYLA performs best on the Ill-Constrained Quadratic and Rosenbrock functions, demonstrating very strong performance on the Matyas function as well.
\\[8pt]
\textbf{COBYQA}: COBYQA performs best on the Matyas function and shares the top performance with COBYLA on the Ill-Constrained Quadratic function, with identical trajectories in the quantitative ranking region. It also performs very well on the Rosenbrock function. The Ill-Constrained Quadratic constraint violation plot shows a small constraint violation by this algorithm.
\\[8pt]
\textbf{CUATRO}: CUATRO exhibits the poorest performance on the Ill-Constrained Quadratic and Matyas functions, with reasonable performance on the Rosenbrock function. Among all constraint violation plots, CUATRO demonstrates the most severe constraint violations. Furthermore, its first five iterations have a standard deviation of zero for all test functions, likely due to an in-algorithm imposed random sampling during these initial iterations.

\newpage
\section{Chemical Engineering Case Studies}\label{sec:d}
In addition to the benchmarking procedure based on mathematical functions, the performance of the algorithms was tested on two chemical engineering case studies. These case studies serve to evaluate the effectiveness of these algorithms in optimizing complex engineering systems. As the reader will discover in the results section, the performance of algorithms can vary significantly between mathematical case studies and real-world engineering problems. This discrepancy can be attributed to a variety of factors, including the fact that real-world problems generally present symmetries, and highly correlated variables which synthetic case studies generally do not. Examining the strengths and limitations of each algorithm in different contexts will provide insights into their suitability for various types of optimization problems and guide the selection of the most appropriate method for a given application.
\\[8pt]
Regarding the algorithms used for the unconstrained controller tuning case study, there are a few changes compared to those used for the unconstrained synthetic case study presented previously in Section \ref{sec:ub}. ENTMOOT and LSQM have proven to be too computationally expensive in higher dimensions, which makes them unsuitable for the controller tuning case study with 32 dimensions. Regarding the BO framework, the in-house implementation that was used in the synthetic case study has been swapped to make way for two state-of-the-art implementations by GPyOpt \cite{thegpyoptauthorsGPyOptBayesianOptimization2016} and TuRBO \cite{erikssonScalableGlobalOptimization2019}. Furthermore, CUATRO-pls was added to showcase the dimensionality reduction capabilities of this add-on to the base version of CUATRO.

\subsection{PID Controller Tuning Problem}

%\ifincludegraphs
\begin{wrapfigure}{l}{0.4\textwidth}
    \begin{center}
    \schemestart
        A \arrow{->[$r_A$]} B \arrow{->[$r_B$]} C
    \schemestop
    \end{center}
  \centering
\includegraphics[width=0.25\textwidth]{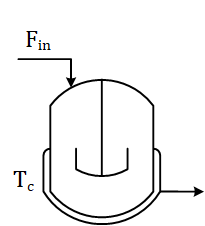}
  \caption{Diagram of a Continuous Stirred Tank Reactor (CSTR) with PID Control}
  \label{fig:cstr}
\end{wrapfigure}
%\fi

This 32-dimensional case study presents a classical chemical engineering problem involving the dynamic control of a \textit{Continuous Stirred Tank Reactor} (CSTR) equipped with a cooling jacket. In this system, the inlet flow rate and the cooling jacket temperature are the manipulated variables, and the reactor temperature is the control variable. Operating under a constant volume steady-state, the CSTR facilitates reactions where reactant A is converted into product B, which subsequently transforms into product C.
\\[8pt]
The objective is to maintain the control variable at its set point. This requires balancing the inherent nonlinearity and complexity of the system. To achieve this, a \textit{Proportional-Integral-Derivative} (PID) controller is employed to adjust the manipulated variables based on the deviation of the control variable from its setpoint. The performance of the control system is assessed by the objective function, defined as the system error. The error is a sum of the deviation of the control variable from its setpoint with a penalty for any changes in control action.
\\[8pt]
The optimization problem is high-dimensional, with 32 controller gain variables. The purpose of this engineering case study is to assess the potential of the benchmarked algorithms in optimizing high-dimensional, complex chemical processes.

\subsubsection{System Definition}
The rate of change of concentration of species A (\(C_A\)) in the reactor is given by:
\begin{equation}
\frac{dC_A}{dt} = \frac{F_{in}}{V} (C_{A_f} - C_A) - r_A
\end{equation}

where the reaction rate \(r_A\) for the reaction \(A \rightarrow B\) is defined as:
\begin{equation}
r_A = k_{0,AB} \exp\left(-\frac{E_{AB}}{RT}\right) C_A
\end{equation}

The rate of change of temperature (\(T\)) in the reactor is given by:
\begin{equation}
\frac{dT}{dt} = \frac{F_{in}}{V} (T_f - T) + \frac{\Delta H_{AB}}{\rho C_p} r_A + \frac{\Delta H_{BC}}{\rho C_p} r_B + \frac{UA}{V \rho C_p} (T_c - T)
\end{equation}

where the reaction rate \(r_B\) for the reaction \(B \rightarrow C\) is defined as:
\begin{equation}
r_B = k_{0,BC} \exp\left(-\frac{E_{BC}}{RT}\right) C_B
\end{equation}

\begin{itemize}
    \setlength{\itemsep}{0pt}   
    \item \( T_c \): Temperature of cooling jacket (K)
    \item \( F_{in} \): Inlet flow rate (\( \text{m}^3/\text{s} \))
    \item \( C_i \) (for \( i \in \{A, B, C\} \)): Concentration of species \( i \) (mol/m\(^3\))
    \item \( T \): Temperature in CSTR (K)
    \item \( T_f \): Feed temperature (K)
    \item \( C_{A_f} \): Feed concentration of A (mol/m\(^3\))
    \item \( V \): Reactor volume (m\(^3\))
    \item \( \rho \): Density of mixture (kg/m\(^3\))
    \item \( C_p \): Heat capacity of mixture (J/kg·K)
    \item \( UA \): Overall heat transfer coefficient (W/K)
    \item \( \Delta H_{j} \) (for \( j \in \{AB, BC\} \)): Heat of \( j \) reaction (J/mol)
    \item \( E_{j} \) (for \( j \in \{AB, BC\} \)): Activation energy for \( j \) reaction (J/mol)
    \item \( k_{0,j} \) (for \( j \in \{AB, BC\} \)): Pre-exponential factor for \( j \) reaction (1/s)
\end{itemize}

\subsubsection{PID Controller Benchmarking Results}
%\ifincludegraphs
\begin{figure}[H]
    \centering
    \includegraphics[width=0.8\textwidth]{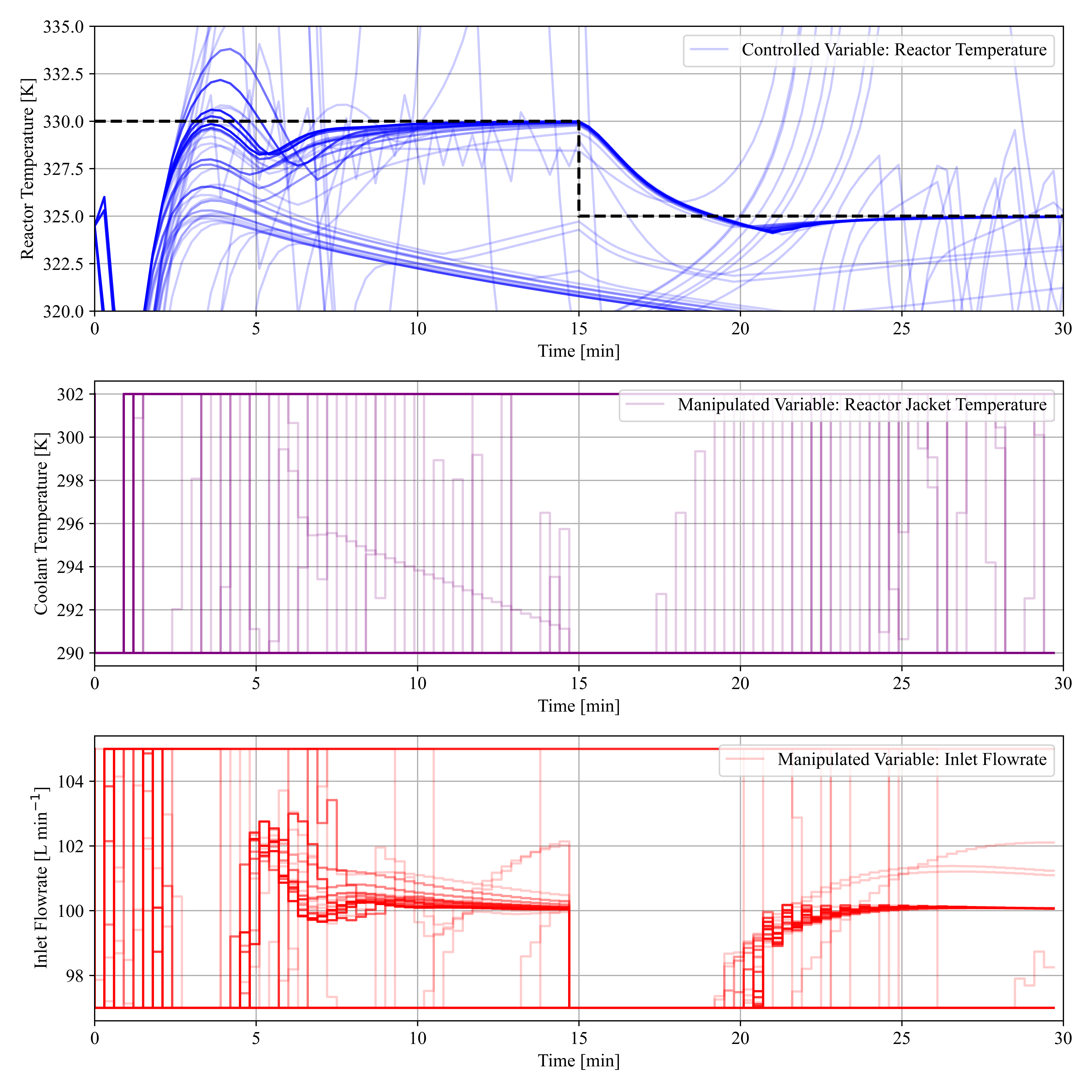}
    \caption{Example Training Trajectories of PID Controller}
    \label{fig:cstr_pid_traj}
\end{figure}
%\fi
\begin{figure}
    \centering
    \includegraphics[width=0.75\linewidth]{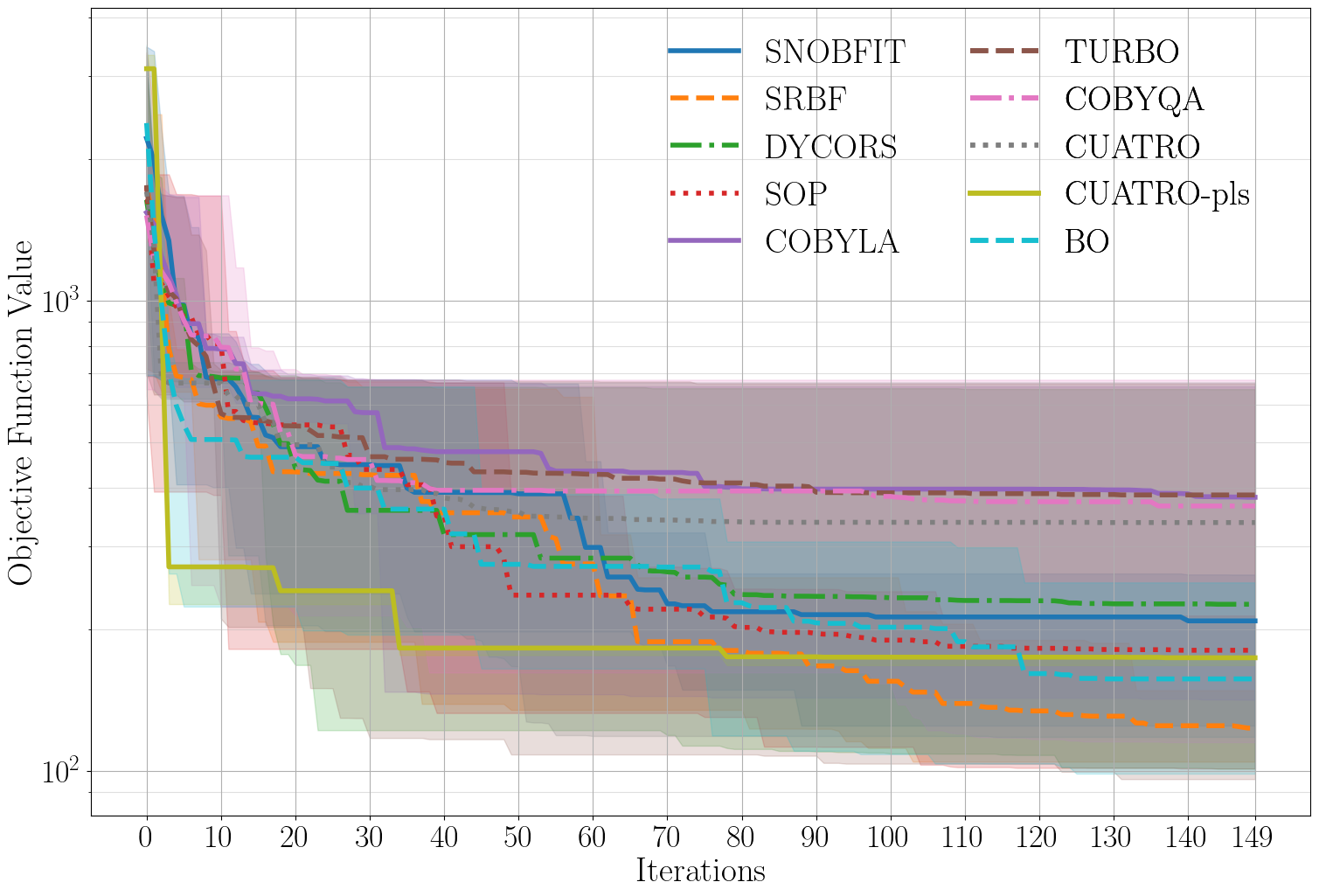}
    \caption{Trajectories for the CSTR-PID Case Study. The budget (trajectory length) is 150 evaluations. Thick lines represent the mean performance on each iteration over 10 repetitions and shaded regions indicate the $90\%$-$10\%$-interval enveloping trajectories for all repetitions.}
\end{figure}

\begin{table}[H]
    \centering
    \captionof{table}{Benchmarking Ranking for CSTR-PID Case Study} \label{tab:uncon_ALL}
    \begin{tabular}{lcc}
    \hline
        \textbf{Algorithm} & \textbf{Overall Performance (All)} & \textbf{Rank} \\
    \hline
         \textbf{CUATRO-pls} & 1.00 & 1 \\
         \textbf{SRBF} & 0.83 & 2  \\
         \textbf{BO} & 0.75 & 3  \\
         \textbf{SOP} & 0.73 & 4  \\
         \textbf{DYCORS} & 0.64 & 5  \\
         \textbf{SNOBFIT} & 0.61 & 6  \\
         \textbf{CUATRO} & 0.32 & 7  \\
         \textbf{COBYQA} & 0.17 & 8  \\
         \textbf{TURBO} & 0.08 & 9  \\
         \textbf{COBYLA} & 0.00 & 10  \\
    \hline
    \end{tabular}
\end{table}
In this study, the performance of several optimization algorithms on a PID controller tuning case was evaluated. The algorithms tested included CUATRO-pls, SRBF, COBYLA, and TuRBO. The results indicate that CUATRO-pls was the best-performing algorithm, followed by SRBF. COBYLA was the worst-performing algorithm, with TuRBO being the second worst.
\\[8pt]
To demonstrate the process of training the PID controller by continuously updating the controller gain variables, an exemplar set of training trajectory graphs has been included - these trajectories are derived from a CUATRO-pls optimization run. These graphs show the optimization process over time, where each algorithm was allowed 150 function evaluations. At each function evaluation, the algorithm updated the set of controller gains, resulting in different control trajectories.
\\[8pt]
There is one graph for each manipulated variable and one for the control variable. By overlaying these trajectories, a visual representation of the PID controller becoming more effective at maintaining the control variable at its setpoint is provided. The earlier, more poorly-tuned trajectories are shown in a fainter color, while the darker colors depict the later trajectories as the optimization progresses. It can be observed that the darker trajectories more closely follow the setpoint, indicating improved performance of the PID controller.

\subsection{Williams-Otto Benchmark Problem}
The Williams and Otto CSTR is a widely studied example, frequently used to benchmark algorithms. This process, depicted below, involves feeding the reactor with two pure component streams $F_a$ and $F_b$ for reactants $A$ and $B$ respectively. $A$ and $B$ react to form an intermediate product $C$, which further reacts with $B$ to yield the desired products, $P$ and $E$. A side reaction occurs between $C$ and $P$, resulting in the formation of $G$ which has no commercial value and is considered waste. The reactions and reactor scheme are displayed in Figure \ref{fig:WO_scheme}:
%\ifincludegraphs
%\begin{wrapfigure}{l}{0.4\textwidth}
%    \begin{center}
%    \begin{align*}
%        \schemestart
%        \chemfig{A} + \chemfig{B} \arrow{->} \chemfig{C}
%        \schemestop
%        \\
%        \schemestart
%        \chemfig{B} + \chemfig{C} \arrow{->} \chemfig{P} + \chemfig{E}
%        \schemestop
%        \\
%        \schemestart
%        \chemfig{C} + \chemfig{P} \arrow{->} \chemfig{G}
%        \schemestop
%    \end{align*}
%    \end{center}
%  \centering
%\includegraphics[width=0.25\textwidth]{images/Figure_60.png}
%  \caption{Diagram of a Continuous Stirred Tank Reactor (CSTR) with PID Control}
%  \label{fig:cstr}
%\end{wrapfigure}
%\fi
\begin{figure}[H] 
    \begin{subfigure}{0.5\linewidth}
        \begin{align*}
            A + B &\rightarrow C \\
            B + C &\rightarrow P + E \\
            C + P &\rightarrow G
        \end{align*}
        \caption{Chemical Reactions}
    \end{subfigure}
    \begin{subfigure}{0.3\linewidth}
        \centering
            \includegraphics[width=\linewidth]{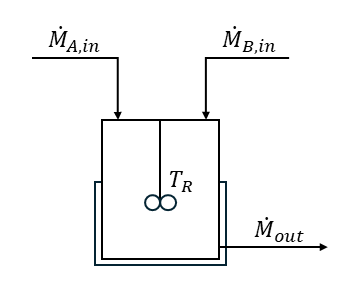}
            \caption{Reactor Scheme}            
    \end{subfigure}           
    \caption{Setting for the Williams-Otto benchmarking problem} 
\end{figure} \label{fig:WO_scheme}

The process is modeled in steady state using mass balance equations, with reactor temperature $T_R$ and component B flow rate $\dot{M}_{B,in}$ as controlled variables. The flow rate of reactant A, $\dot{M}_{A,in}$ and the total mass in the system are maintained at constant values. The objective is to maximize the profit-flow given by Equation \ref{eq:WO}:
\begin{equation} \label{eq:WO}
    \begin{aligned}[t]
        &\max &&f(T_R,\dot{M}_B) = price_R \dot{M}_R+price_E \dot{M}_E-cost_A \dot{M}_{A,in}-cost_B \dot{M}_{B,in} \\
        & \text{s.t. } &&g_1(w_{A,out}) \leq 0.12 \\
        & \quad &&g_2(w_{G,out}) \leq 0.08
    \end{aligned}
\end{equation}
with mass fractions $w_{A,out}$, and $w_{B,out}$ respectively. For brevity, the complete set of mass balance equations and kinetic rate equations for this reaction system, as detailed by \cite{mendozaAssessingReliabilityDifferent2016}, are not reproduced here.

\subsubsection{Williams-Otto Benchmarking Results}
\begin{figure}[H]
    \begin{subfigure}{0.5\linewidth}
        \centering
            \includegraphics[width=\linewidth]{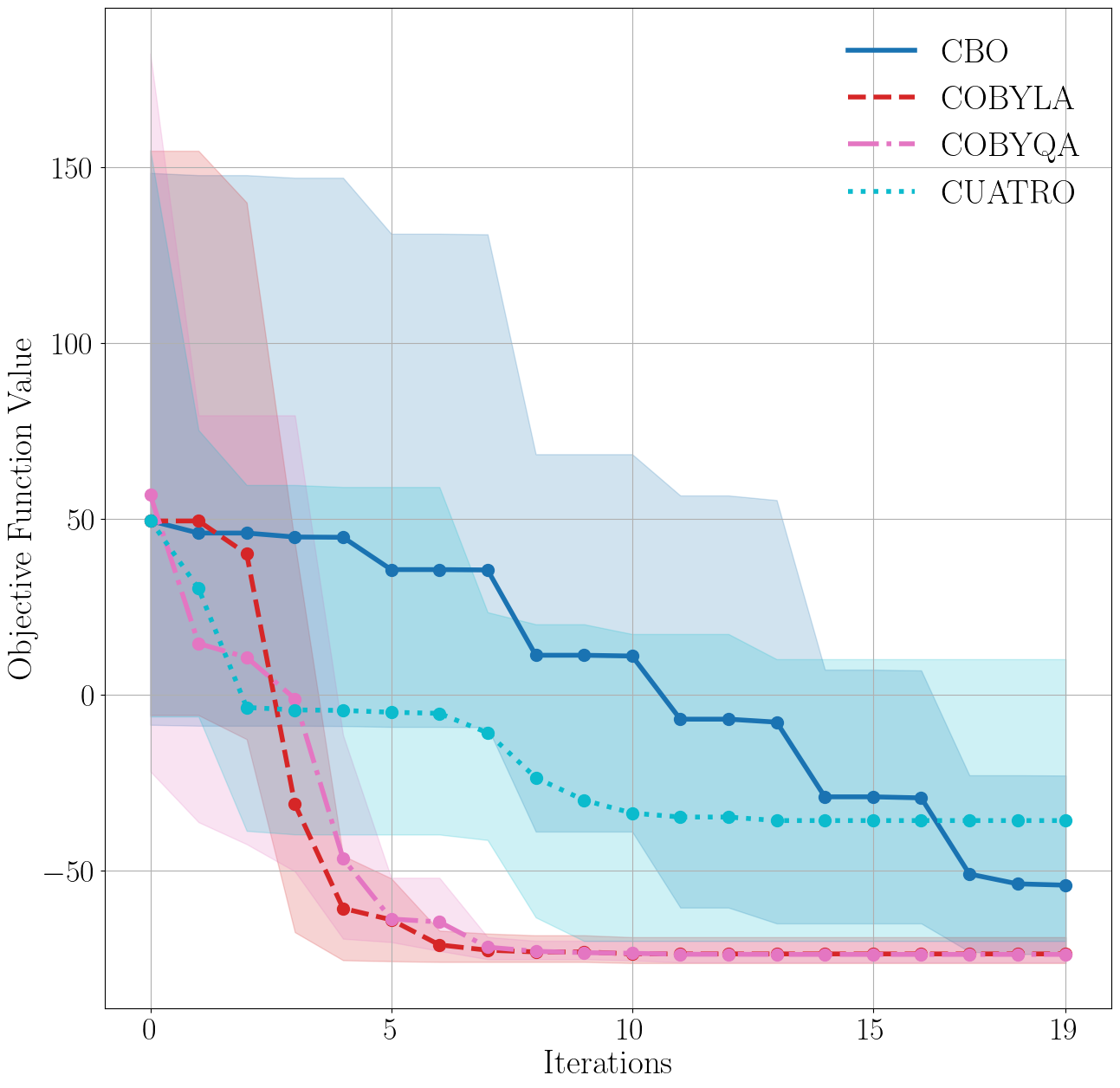}
            \caption{Convergence Plot}
    \end{subfigure}
    \begin{subfigure}{0.5\linewidth}
        \centering
            \includegraphics[width=\linewidth]{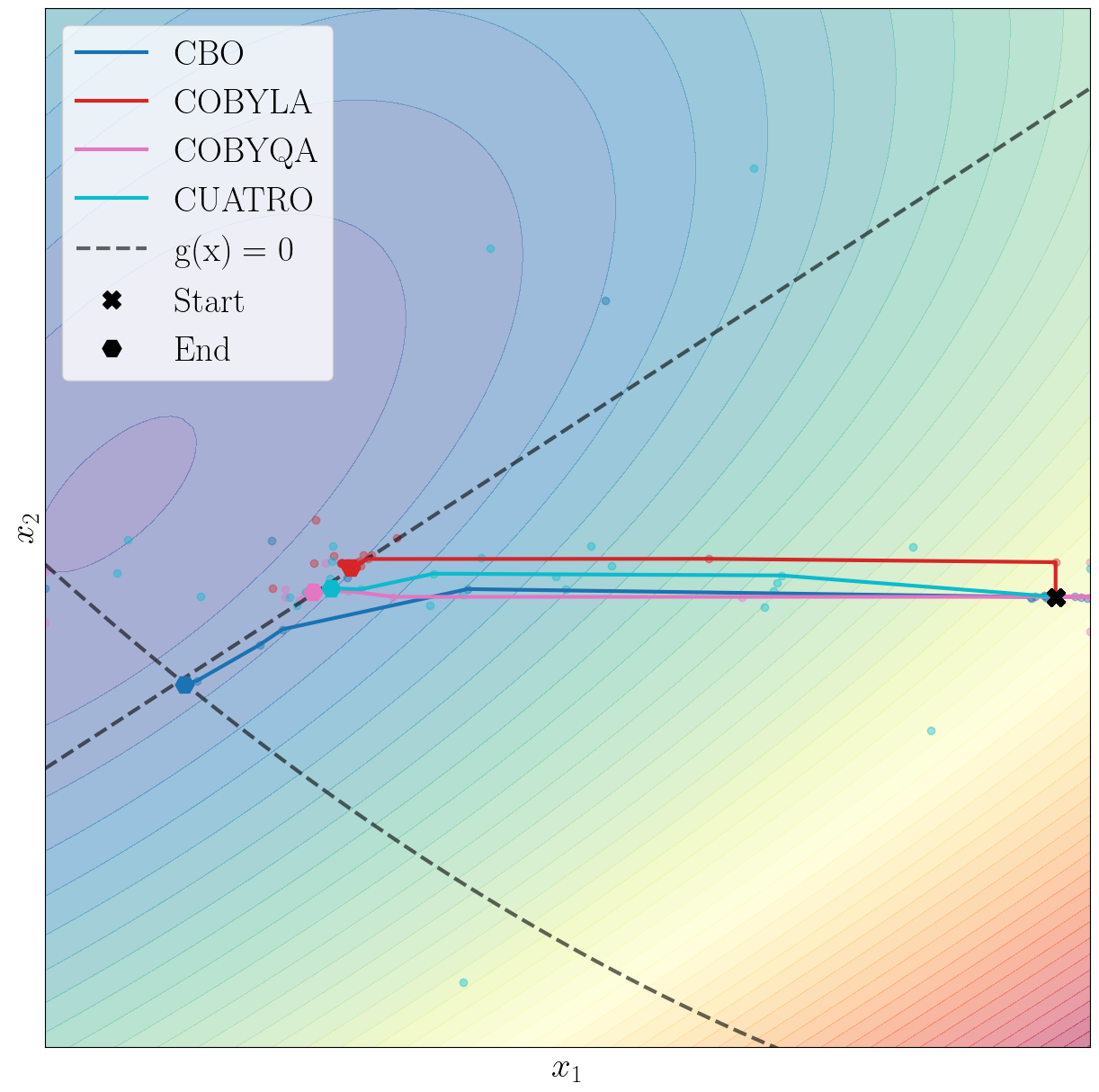}
            \caption{Trajectory Plot}            
    \end{subfigure}           
    \caption{Visualization for performance on Williams-Otto benchmarking problem: (a) Convergence plot, showing mean objective function values (thick lines) and $90\%$-$10\%$-intervals over 10 repetitions from 10 different starting points with a budget (trajectory length) of 20 evaluations. (b) Williams-Otto contour with dashed constraint-line and a single trajectory for each algorithm from an exemplary shared starting point. The lines show best-so-far evaluation positions and the scattered points show remaining evaluation positions. } 
\end{figure}
\begin{table}[H]
    \centering
    \captionof{table}{Benchmarking table for the constrained optimization of Williams-Otto}
    \begin{tabular}{lccc}
    \hline
          &\textbf{$p_a$} & \textbf{Feasible Samples} & \textbf{Mean Violation} \\
    \hline
         \textbf{COBYQA} & 1.00 & 81.75\% & 0.0062 \\
         \textbf{COBYLA} & 1.00 & 84.28\% & 0.0044 \\
         \textbf{CBO} & 0.00 & 90.86\% & 0.0085 \\
         \textbf{CUATRO} & 0.14 & 80.10\% & 0.0125 \\
    \hline
    \end{tabular}
\end{table}

Based on the performance scores, the best-performing algorithms are COBYQA and COBYLA with a perfect score of 1.0, followed by CUATRO with a score of 0.14. On the other end of the spectrum, the in-house implementation of Constrained BO is the worst-performing algorithm with a score of 0.0.

\section{Code}\label{sec:e}
%\todo[inline]{add github}
The source-code for this project including benchmarking routines and algorithms can be found at:
\begin{center}
 https://github.com/OptiMaL-PSE-Lab/DDO-4-ChemEng   
\end{center}

\newpage
\section{Concluding Remarks}\label{sec:f}

This chapter has presented the application and evaluation of various model-based derivative-free optimization algorithms for chemical processes, emphasizing the utility of data-driven approaches in enhancing optimization performance. The key insights from this study can be summarized as follows:

\subsection{Performance Comparison of Algorithms}

The performance evaluation of several optimization algorithms, including CUATRO-pls, SRBF, COBYLA, and TuRBO, highlighted significant differences in their effectiveness. CUATRO-pls emerged as the best-performing algorithm, particularly in the context of tuning a PID controller for a CSTR. SRBF followed closely, while COBYLA and TuRBO lagged behind, additionally, SNOBFIT performed much better in this case study. This ranking indicates that while traditional algorithms like COBYLA still hold value, newer and more sophisticated approaches such as CUATRO-pls and SRBF can offer substantial improvements in optimization tasks. Interestingly, some algorithms, such as CUATRO and SNOBFIT performed significantly differently in test functions compared to chemical engineering case studies, highlighting that assessing the performance of algorithms in engineering applications must be done in addition to a screening over their performance in traditional mathematical test functions. 

\subsection{Benchmarking on Mathematical Functions}

The study employed several well-known mathematical objective functions, including the Ackley, Levy, Rosenbrock, and a quadratic ill-conditioned function, to benchmark the algorithms. These functions, selected for their diverse characteristics—ranging from multi-modal to uni-modal—provided a testbed for evaluating algorithm performance. The results indicated that the algorithms' effectiveness varied with the complexity and nature of the objective functions. For example, CUATRO-pls performed exceptionally well on unimodal functions like Rosenbrock but showed mixed results on multimodal functions like Ackley and Levy, and similar conclusions can be drawn for other algorithms.

\subsection{General Observations and Future Work}

The benchmarking results highlighted the importance of algorithm selection based on the specific characteristics of the optimization problem at hand. Surrogate models, particularly those based on BO and RBFs, demonstrated robust performance across various scenarios. However, the study also highlighted areas for future exploration, such as the integration of noise in optimization processes and the application of these algorithms to more complex, real-world chemical engineering problems.
\\[8pt]
Finally, this study reaffirms the potential of surrogate-based optimization algorithms in enhancing the efficiency and accuracy of optimization in chemical processes for specific applications. The comparative analysis provides a reference for selecting appropriate algorithms, while the insights gained pave the way for future advancements and applications in the field. The continuous evolution of these algorithms promises further improvements, making DDO an indispensable tool in chemical engineering.
\newpage
\printbibliography
\end{document}